%% file: main.tex
\documentclass[manuscript,authordraft,review=false,timestamp=false,nonacm]{acmart}
\SetWatermarkText{}

\input{00_preamble}

\AtBeginDocument{%
  }

\title{A Modular Approach to Stochastic Optimisation for Inverse Problems Using the Core Imaging Library}

\input{version3/authors.tex}

\begin{abstract}

The Core Imaging Library (CIL) is an open-source versatile Python framework for solving inverse problems with special emphasis on imaging applications such as computed tomography (CT), using a plug-in architecture for data and operators, interfacing to toolboxes such as ASTRA, TIGRE and SIRF. A key component of CIL is its optimisation module enabling users to flexibly combine mathematical operators and functionals to form smooth and non-smooth optimisation problems and solve these with a range of first-order algorithms. The present work introduces an expansion of CIL with a new modular framework for stochastic optimisation, allowing researchers to easily use a variety of existing stochastic optimisation algorithms as well form new ones by combining modular building blocks.  
Users can flexibly configure algorithmic components, adapt to diverse problem structures, and experiment with various sampling and step size strategies. 
Rather than individual black-box implementations of each fixed algorithm with significant redundancies, our design is modular providing building blocks that can be flexibly combined to realise a wealth of algorithm instances.
The framework is particularly well-suited for large-scale applications, where stochastic methods offer notable computational advantages over deterministic approaches. To demonstrate its versatility and practical utility, we present experiments on real-world datasets from imaging inverse problems, such as X-Ray CT and Positron Emission Tomography (PET) reconstruction.

In summary, the presented software expansion aims to support the research community with a robust, extensible optimisation suite  for developing, testing, and benchmarking stochastic methods for inverse problems. 

\end{abstract}

\begin{document}

\maketitle

\input{version3/01_introduction}

\input{version3/01_stochastic_introduction}
\input{version3/03_implementation}

\input{version3/04_example_ct}

\input{version3/05_example_pet}
\input{version3/06_discussion}
\input{version3/07_conclusion}

\section*{Data availability statement}

The code to reproduce all results in this paper will be made available as open source upon acceptance of this manuscript. The battery dataset, Lithium-ion 10440 NMC532 cell, is available through the Battery Imaging Library \cite{Docherty2025} at \href{https://doi.org/10.5281/zenodo.17930456}{https://doi.org/10.5281/zenodo.17930456}. The PET dataset is available at \href{https://zenodo.org/records/1304454}{https://zenodo.org/records/1304454}.

\section*{Acknowledgments}

EPap acknowledges funding through the Innovate UK Analysis for Innovators (A4i) program for the project \emph{Denoising of chemical imaging and tomography data}, during which the initial design and implementation of the software were carried out. CD was supported by PET++, UKRI EPSRC EP/S026045/1. SP was supported in part by National Physical Laboratory through the National Measurement System of the Department for Science, Innovation and Technology (NPL grant PO496581) and by the Engineering and Physical Sciences Research Council (EPSRC) through an Industrial CASE studentship (EPSRC grant EP/W522077/1) and Mediso Medical Imaging Systems. 
The authors wish to acknowledge the support of the Collaborative Computation Project in Tomographic Imaging (via the EPSRC grant EP/T026677/1 and via STFC through the UKRI Digital Research Infrastructure programme), as well as the support of the Collaborative Computation Project in Synergistic Reconstruction for Biomedical Imaging, CCP SyneRBI (via  the EPSRC grant EP/T026693/1 and by UKRI STFC through the UKRI Digital Research Infrastructure programme).

This work made use of computational support by CoSeC, the Computational Science Centre for Research Communities, through CCP-SyneRBI and CCPi.
The authors wish to acknowledge the STFC cloud for providing the computational platform for this work.

The authors would also like to acknowledge Laura Murgatroyd, Matthias Ehrhardt, Tang Junqi,  Imraj Singh, Robert Twyman, Antony Vamvakeros, Casper da Costa-Luis, Daniel Deidda, Ashley Gillman, Evgueni Ovtchinnikov, Georg Schramm, Hannah Robarts, and Zeljko Kereta for their valuable discussions and feedback that helped shape this work.

\bibliographystyle{IEEEtran}
\bibliography{refs}

\end{document}

%% file: 00_preamble.tex
\usepackage{graphicx} 
\usepackage[utf8]{inputenc}
\usepackage{float}
\usepackage{graphicx}
\usepackage{textcomp}
\usepackage{xcolor}
\usepackage{tikz}
\usepackage{pgfplots}
\usepackage[export]{adjustbox}
\usetikzlibrary{spy}
\usepackage{graphicx}
\usepackage{tabularx}
\usepackage{xcolor}
\usepackage{amsmath,amsfonts, amsthm}
\usepackage{wrapfig}
\usepackage{bm}
\usepackage{tikz}
\usepackage{multirow}
\usepackage{booktabs}
\usepackage{blindtext}
\usepackage{multicol}
\usepackage{tikzscale}
\usepackage{numprint}
\npdecimalsign{.}
\nprounddigits{3} 
\usepackage{caption}
\usepackage{subcaption}
\usepackage{enumitem}
\usepackage{float}
\usepackage{pict2e}
\usepackage{siunitx}
\usepackage{gensymb}
\usepackage[english]{babel}
\usepackage{lipsum}
\usepackage{diagbox}
\usepackage{url}
\usepackage{hyperref}

\usepackage{array}   
\renewcommand{\arraystretch}{2.5}
\usepackage[most]{tcolorbox}
\definecolor{darkspringgreen}{rgb}{0.09, 0.45, 0.27}
\usepackage[cachedir=minted-cache]{minted}

\usepackage{mathtools}

\usepackage[percent]{overpic}
\usepackage{algorithm,algcompatible}

\DeclareMathOperator*{\argmin}{arg\,min}
\algnewcommand\INPUT{\item[\textbf{Input:}]}%
\algnewcommand\PARAMETER{\item[\textbf{Parameters:}]}%
\algnewcommand\OUTPUT{\item[\textbf{Output:}]}%

\usepackage{cleveref}

\usepackage{comment}
\usepackage{array}

\usepackage{pythonhighlight}
\newcommand{\code}[1]{\pyth{#1}\xspace}

\usepackage{algorithm}
\usepackage{algpseudocode}
\algnewcommand{\Inputs}[1]{%
  \State \textbf{Inputs:}
  \Statex \hspace*{\algorithmicindent}\parbox[t]{.8\linewidth}{\raggedright #1}
}
\algnewcommand{\Initialize}[1]{%
  \State \textbf{Initialize:}
  \Statex \hspace*{\algorithmicindent}\parbox[t]{.8\linewidth}{\raggedright #1}
}

\usepackage[most]{tcolorbox}
\definecolor{darkspringgreen}{rgb}{0.09, 0.45, 0.27}

\newcommand{\sysmat}{\bm{A}}
\newcommand{\image}{\bm{u}}

\newcommand{\bx}{\bm{x}}
\newcommand{\by}{\bm{y}}
\newcommand{\bz}{\bm{z}}

\newcommand{\bA}{\bm{A}}
\newcommand{\bb}{\bm{b}}

\DeclareMathAlphabet\mathbfcal{OMS}{cmsy}{b}{n}
\newcommand{\Kop}{\bm{K}}

\newcommand{\gradest}{\bm{G_k}}

\definecolor{applegreen}{rgb}{0.55, 0.71, 0.0}

\usepackage{comment}

\newcommand{\ie}{\textit{i}.\textit{e}., }

\let\oldnl\nl
\newcommand{\nonl}{\renewcommand{\nl}{\let\nl\oldnl}}

\usepackage{subcaption}

\usepackage{xr}
\makeatletter
\newcommand*{\addFileDependency}[1]{
  \typeout{(#1)}
  \@addtofilelist{#1}
  \IfFileExists{#1}{}{\typeout{No file #1.}}
}
\makeatother

%% file: version3/authors.tex
\author{Evangelos Papoutsellis}
\email{epapoutsellis@gmail.com}
\orcid{0000-0002-1820-9916}
\affiliation{
  \institution{Finden Ltd, Rutherford Appleton Laboratory}
  \city{Harwell Campus}
  \state{Oxfordshire}
  \country{UK}
}

\author{Margaret A. G. Duff}
\email{margaret.duff@stfc.ac.uk}
\orcid{0000-0003-1014-3147}
\affiliation{
  \institution{Science and Technology Facilities Council, Rutherford Appleton Laboratory}
  \city{Harwell Campus}
  \state{Oxfordshire}
  \country{UK}
}

\author{Jakob S. J\o{}rgensen}
\email{jakj@dtu.dk}
\orcid{0000-0001-9114-754X}
\affiliation{
  \institution{Department of Applied Mathematics and Computer Science, Technical University of Denmark}
  \city{Copenhagen}
  \country{Denmark}
}

\author{Sam Porter}
\email{sam.porter.18@ucl.ac.uk}
\orcid{0000-0002-1815-4716}
\affiliation{
  \institution{Institute of Nuclear Medicine, University College London}
  \city{London}
  \country{UK}
}
\author{Claire Delplancke }
\email{claire.delplancke@edf.fr}
\orcid{0000-0001-7483-0419}
\affiliation{
  \institution{EDF Lab Paris-Saclay}
  \city{Paris}
  \country{France}
}

\author{Gemma Fardell}
\email{gemma.fardell@stfc.ac.uk}
\orcid{0000-0003-2388-5211}
\affiliation{
  \institution{Science and Technology Facilities Council, Rutherford Appleton Laboratory}
  \city{Harwell Campus}
  \state{Oxfordshire}
  \country{UK}
}

\author{Edoardo Pasca}
\email{edoardo.pasca@stfc.ac.uk}
\orcid{0000-0001-6957-2160}
\affiliation{
  \institution{Science and Technology Facilities Council, Rutherford Appleton Laboratory}
  \city{Harwell Campus}
  \state{Oxfordshire}
  \country{UK}
}

\author{Kris Thielemans}
\email{k.thielemans@ucl.ac.uk}
\orcid{0000-0002-5514-199X}
\affiliation{
  \institution{Institute of Nuclear Medicine and UCL Hawkes Institute, University College London}
  \city{London}
  \country{UK}
}

%% file: version3/01_introduction.tex
\section{Introduction}

Imaging Inverse problems are a fundamental challenge in computational imaging, where the objective is to estimate an unknown quantity from indirect, incomplete and often noisy measurements. These arise in various fields including medical imaging, materials science, remote sensing and geoscience. For finite dimensional spaces $\mathbb{X}, \mathbb{Y}$ and a forward operator $\bm{A}:\mathbb{X}\rightarrow\mathbb{Y}$, many inverse problems are characterised by the following operator equation
\begin{equation}
\bm{A}\image\approx \bb,
\label{eq:general_operator_equation}
\end{equation}
where $\image\in\mathbb{X}$ is the unknown quantity and $\bb\in\mathbb{Y}$ is the measurement data. Often these problems are ill-posed, \ie they lack existence, uniqueness or stability of a solution. To address this, a common approach is to integrate prior information of the solution in the form of regularisation in a variational (optimisation-based) method.

Variational methods with regularisation aim to recover an approximated solution $\image\in\mathbb{X}$ of~\eqref{eq:general_operator_equation} by solving the following optimisation problem
\begin{equation}
\min_{\bm{u}\in\mathbb{X}},D(\bm{A}\bm{u}, \bb) + \alpha R(\bm{u}).
\label{eq:general_optimisation}
\end{equation}
Here, $D:\mathbb{Y}\times \mathbb{Y}\rightarrow\mathbb{R}_{\infty}:=\mathbb{R}\cup\{\infty\}$ denotes the fidelity term, which measures the discrepancy between the forward model applied to the unknown quantity and the measurement data. Depending on the noise statistics of the measurement data, $\bb$, the fidelity term often takes the form of a least squares distance, an absolute difference, or a Kullback-Leibler divergence for Gaussian, impulse~\cite{Nikolova2004}, and Poisson noise~\cite{Le2007}, respectively. The regularisation term, $R(\bm{u})$, promotes certain properties to the solution, $\bm{u}$, such as smoothness, sparsity or edge preservation and is usually weighted by a regularisation parameter, $\alpha>0$, acting as a balance between the two terms in~\eqref{eq:general_optimisation}. Popular choices for the regularisation term are Tikhonov~\cite{tikhonov1977solutions}, Total Variation (TV)~\cite{Rudin1992} and the high order extension, namely the Total Generalised Variation (TGV),~\cite{Bredies2010}, Total Nuclear Variation,~\cite{Holt2014} and more general tensor-based structure regularisation~\cite{Lefkimmiatis2013}.

Numerous algorithms have been proposed in the literature~\cite{ChambollePock2016, BenningBurger2018, Bredies2020, Ryu2022, Mueller2012} designed to solve~\eqref{eq:general_optimisation}, each tailored to specific imaging modalities, data fidelity requirements, and prior assumptions. Some of the present authors presented the Core Imaging Library (CIL) in~\cite{Jorgensen2021}, a versatile open-source software framework designed for imaging inverse problems with an emphasis on advanced tomographic imaging tasks. Besides important utilities related to pre-processing, post-processing and visualisation of imaging data encountered in computed tomography (CT) and related modalities, a flexible optimisation framework was designed to find solutions for inverse problems. The flexible \code{cil.optimisation} module, based on three main building blocks, namely \emph{Functions}, \emph{Operators}, and \emph{Algorithms}, allows the users to construct solutions with minimal effort for various optimisation problems tailored to specific imaging applications. Several algorithms have been implemented in CIL to solve minimisation problems in the form of~\eqref{eq:general_optimisation}. These are Gradient Descent (GD), Conjugate Gradient Least Squares (CGLS), as well as proximal splitting-based algorithms including Proximal Gradient Descent (PGD) (aliased to ISTA), Fast Iterative Shrinkage Thresholding Algorithm (FISTA)~\cite{Beck2009}, Primal-Dual Hybrid Gradient (PDHG) \cite{Chambolle2010} and Linearised Alternate Direction Method of Multipliers (LADMM)~\cite{OConnor2018}.

These algorithms, implemented in CIL, have been utilised in various applications with a strong emphasis on tomographic reconstruction problems on large and challenging data sets, including low-count, non-standard geometries, limited angle,  multispectral or multimodality data, including conventional X-ray CT in~\cite{Jorgensen2021}, dynamic X-ray CT,~\cite{Papoutsellis2021}, limited angle CT,~\cite{Jorgensen2023}, Hyperspectral and Neutron Tomography in~\cite{Warr2021, Evelina2021} and fusion plasma imaging, ~\cite{SimmendefeldtSchmidt2024}.   Moreover, the CIL optimisation offers a unique solution to the medical imaging community since it is fully compatible with the Synergistic Image Reconstruction Software (SIRF)~\cite{Ovtchinnikov2020}, open-source software for applications in medical imaging including Positron Emission Tomography (PET), Single Photon Emission Computed Tomography (SPECT) and Magnetic Resonance Imaging (MRI),~\cite{Brown2021, mcir}. Finally, the CIL library has also been applied beyond tomography, notably in synthetic aperture radar (SAR) imaging,~\cite{Watson2024}.

Despite the substantial advancements and successes of iterative regularisation algorithms for inverse problems in recent years, uptake of their use is limited because they can be slow to converge or have high memory demands. Traditional first-order solvers and their primal-dual variants are typically designed to operate in a deterministic mode, processing the entire dataset for every update, \ie one data pass for each iteration. This approach can be computationally intensive, especially as tomographic scanners continue to advance, producing increasingly larger and higher-resolution datasets. This challenge has motivated the creation of algorithms that provide faster convergence rates per data pass and reduce per-iteration computational costs, with many utilising stochastic optimisation techniques, see Ehrhardt et al.~\cite{Ehrhardt_2025} for a comprehensive review. Unlike traditional deterministic methods that process the entire dataset in every iteration, stochastic approaches randomly select and process a subset of the acquired data, for example, a small number of tomographic projection views. Similar to Stochastic Gradient Descent (SGD) methods widely used in machine learning, this strategy has the potential to reduce per-iteration computation, while still converging to the solution. These advantages position stochastic optimisation as a compelling solution for tackling the challenges posed by high resolution, large-scale datasets in modern inverse problems.

\subsection{Contribution}

In this paper, we extend the deterministic CIL optimisation framework by introducing a unique and user-friendly design that enables efficient implementation and testing of stochastic optimisation techniques.  The rapid advancement in optimisation research results in the frequent publication of new algorithms, making it challenging to evaluate their practical performance for specific imaging problems. Our software addresses this by providing users with an easy-to-use environment to test and validate these emerging algorithms, enabling them to assess their effectiveness on real-world problems, including X-ray CT, PET, SPECT and MR reconstruction. Furthermore, CIL version 24.3.0 provides additional utilities that enhance its functionality and flexibility, such as data splitting and sampling methods, preconditioning, line search methods for non-constant step size and callback methods to monitor and analyse the progress of algorithms. These features enable users to tailor the optimisation process to their specific needs, making it particularly effective for addressing complex inverse problems. Finally, the proposed modular design allows users to easily extend its functionality by incorporating novel optimisation methods, fostering innovation and adaptability in tackling a wide range of challenging problems.

The scope of this paper is not to provide an exhaustive literature review or comparison of algorithms  for~\eqref{eq:optimisation_three_terms}  but rather to introduce and demonstrate the superiority of stochastic optimisation alongside software solutions designed to support our community in tackling diverse optimisation problems. Our emphasis lies on delivering a framework that prioritises flexibility, versatility, and reproducibility, enabling users to efficiently prototype, improve, and extend optimisation methods for their applications. The experiments in the paper are designed to be indicative of the expected behaviour of these algorithms applied to real-world problems.
This paper describes the functionality introduced in CIL version 24.3.0.

\subsection{Related Work}

There are many open-source libraries that provide optimisation algorithms for inverse problems arising in imaging, each targeting different applications and communities. PyProximal~\cite{Ravasi2024} is a Python library for convex optimisation, developed as an integral part of the PyLops framework~\cite{Ravasi2020}, which is a library for matrix-free linear algebra and large-scale inverse problems. SCICO (Scientific Computational Imaging Code)~\cite{Balke2022} is a Python package for scientific imaging and inverse problems that uses JAX as its numerical backend and supports learned priors and deep-learning reconstruction methods. DeepInv~\cite{tachella2025deepinverse} is an open-source PyTorch library for imaging inverse problems that offers physics-based forward models, including for MRI, Tomography and Ptychography; learned reconstruction networks; and plug-and-play denoisers. PyXu~\cite{PyXu} is an open-source Python framework for computational imaging pipelines that emphasises modular operators and scalable execution on CPUs and GPUs. The Operator Discretization Library (ODL)~\cite{ODL} is a Python framework for prototyping inverse problems using abstract forward operators and iterative reconstruction methods. Although these libraries differ in emphasis - ranging from variational regularisation to deep learning to scalable operator algebra - they mainly expose deterministic optimisation algorithms such as ISTA, FISTA, ADMM and PDHG. To the best of our knowledge, our proposed software extension is the first to offer a unified and extensible interface for stochastic optimisation in this setting; none of the aforementioned libraries provides similar capabilities.

Our work introduces a novel software design to the open source optimisation ecosystem. Instead of simply adding stochastic optimisation algorithms to an existing module, we propose a unified optimisation interface that allows stochastic methods, such as data splitting, sampling-based updates, and variance-reduced gradient schemes, to be combined seamlessly with classical methods such as proximal or primal-dual algorithms.
As discussed later, the key distinction between deterministic and stochastic gradient-based algorithms lies in how the gradient of the objective function is evaluated. To leverage this, we encapsulate stochastic gradient estimation within dedicated functions. This design allows us to reuse CIL's existing deterministic gradient-based algorithms for stochastic optimisation without creating separate workflows. The same base algorithm can therefore execute either a standard deterministic method or its stochastic counterpart.
This approach provides a single, coherent interface for deterministic and stochastic optimisation in imaging inverse problems.

\subsection{Overview of Paper}

In \Cref{background-theory}, we introduce the key concept in our framework: creating stochastic algorithms by plugging stochastic gradient estimators into deterministic, gradient-based algorithms. In this section, we do this in a theoretical way, but also draw links to the software.  In \Cref{sec:implementation} we describe the framework again with a software walk-through of the implementation in CIL, providing readers with a step by step method to build up a stochastic algorithm. We then provide examples illustrating all aspects of CIL's stochastic optimisation framework through two case studies: a large-scale parallel-beam X-ray CT dataset in \Cref{sec:ct}, and a low-count PET dataset in \Cref{sec:pet}.

%% file: version3/01_stochastic_introduction.tex
\section{Theory and Implementation  \label{background-theory}}

The essence of our stochastic framework is the idea of plugging \textbf{stochastic gradient estimators} into \textbf{deterministic gradient-based algorithms}, providing a quick and easy way of combining algorithms to give a range of different methods. In this section, we first introduce the deterministic algorithms, then the idea of a stochastic gradient estimator and how this encompasses many stochastic algorithms in the literature.

Along the way, we will discuss other aspects of the deterministic algorithms that can be tuned: step sizes, preconditioners and momentum. We will also mention \textbf{partitioning} and \textbf{sampling} for the stochastic algorithms, but these will be covered in more detail in the future sections.

\subsection{Deterministic Optimisation Framework }
\label{sec:deterministic_opt}

For a detailed discussion of deterministic optimisation algorithms available in CIL, see references~\cite{Jorgensen2021,Papoutsellis2021}. In this section, we consider just those algorithms required for the stochastic framework and comparisons. We consider a general optimisation framework that rewrites ~\eqref{eq:general_optimisation}  as the sum of three terms:
\begin{align}
\bx^{*} & \in\argmin_{\bx\in\mathbb{X}} \large{F(\bx):= f(\bx) + g(\bx)+h(\Kop\bx)\large}.
\label{eq:optimisation_three_terms}
\end{align}
Here, $\Kop:\mathbb{X}\rightarrow \mathbb{Y}$ is a continuous linear operator with $$\|\Kop\| = \max\left\{\|\Kop\bx\|: \bx\in \mathbb{X},  \|\bx\|\leq1 \right\}<\infty,$$ and the functions $f, g:\mathbb{X}\rightarrow\mathbb{R}\cup\{\infty\}$ and $h:\mathbb{Y}\rightarrow\mathbb{R}\cup\{\infty\}$ are proper, convex and  lower-semicontinuous.  In addition, $f$ is $L$-smooth, i.e., a differentiable function with Lipschitz gradient $L>0$,
\begin{equation}
\|\nabla f(\bx) - \nabla f(\by)\|\leq L\|\bx-\by\|,\quad \forall \bx, \by \in \mathbb{X}.
\label{eq:Lsmooth}
\end{equation}
The functions $g$ and the convex conjugate of $h$ are assumed to be \textit{proximable}.
Recall that the convex conjugate of $h$ is defined as
\begin{equation}
h^{*}(\bm{y}^*)= \sup_{\bm{y}\in\mathbb{Y}}\{\langle \bm{y},\bm{y}^*\rangle - h(\bm{y})\}
\end{equation}  and the proximal operator of a function $g$ is
\begin{equation}
\mathrm{prox}_{\tau g}(\bz)  = \argmin_{\bx\in \mathbb{X}} \left\{\frac{1}{2}\|\bx-\bz\|^{2} + \tau g(\bx)\right\}
\label{eq:proxg}
\end{equation}
with $\tau>0$. A function is ``proximable" if~\eqref{eq:proxg} has either an analytic solution or can be solved efficiently, up to some precision, using an iterative process. The algorithms available in the latest version of CIL, and discussed in this paper, are written in pseudo-code, in Table~\ref{table:CIL_algorithms}.

\begin{table}[tb]
\noindent
\begin{minipage}[t]{0.45\textwidth}
\begin{algorithm}[H]
\begin{algorithmic}[1]
\State {\bf Parameters:} step-size $\gamma_k >0$, preconditioners $P_k$
\State {\bf Initialize:} $\bx_0 \in \mathbb{X}$
\For{$k=0,\ldots, K-1$}
\State $\bx_{k+1} = \bx_{k} - \gamma_k P_k\gradest(\bx_{k})$
\EndFor
\end{algorithmic}
\caption{GD}
\label{alg:GD}
\end{algorithm}
\end{minipage}\hfill
\begin{minipage}[t]{0.49\textwidth}
\centering
\begin{algorithm}[H]
\begin{algorithmic}[1]
\State {\bf Parameters:} step-sizes $\gamma_k >0$, preconditioners $P_k$
\State {\bf Initialize:} $\bx_0 \in \mathbb{X}$
\For{$k=0,\ldots, K-1$}
\State $\bx_{k+1} = \mathrm{prox}_{\gamma_k g}(\bx_{k} - \gamma_k P_k \gradest(\bx_{k}))$
\EndFor
\end{algorithmic}
\caption{PGD}
\label{alg:PGD}
\end{algorithm}
\end{minipage}
\begin{minipage}[t]{0.45\textwidth}
\begin{algorithm}[H]
\begin{algorithmic}[1]
\State {\bf Parameters:} step-sizes $\gamma_k >0$, momentum parameters $\beta_k> 0$, preconditioners $P_k$
\State {\bf Initialize:} $\bx_0 \in \mathbb{X}$
\For{$k=0,\ldots, K-1$}
\State $\bx_{k+1} = \mathrm{prox}_{\gamma_k g}(\by_{k} - \gamma_k P_k\gradest(\by_{k}))$
\State $\by_{k+1} =\bx_{k+1}+\beta_k\left(\bx_{k+1}-\bx_{k}\right)$
\EndFor
\end{algorithmic}
\caption{APGD}
\label{alg:FISTA}
\end{algorithm}
\end{minipage}\hfill
\begin{minipage}[t]{0.49\textwidth}
\centering
\begin{algorithm}[H]
\begin{algorithmic}[1]
\State {\bf Parameters:} stepsize $\sigma, \tau >0$
\State {\bf Initialize:} $\bx_0 \in \mathbb{X}$, $\by_0 \in \mathbb{Y}$
\For{$k=0,\ldots, K-1$}
\State $\bx_{k+1}=\mathrm{prox}_{\tau g}\left(\bx_k-\tau \gradest(\bx_k)-\tau \Kop^*\by_k\right)$
\State $\Bar{\bx}_{k+1} =2\bx_{k+1} -\bx_k +\tau\left(\gradest(\bx_k)-\gradest(\bx_{k+1})\right)$
\State $\by_{k+1} = \mathrm{prox}_{\sigma  h} \left(\by_k+\sigma \Kop\Bar{\bx}_k\right)$
\EndFor
\end{algorithmic}
\caption{PD3O}
\label{alg:PD3O}
\end{algorithm}
\end{minipage}

\caption{Deterministic optimisation algorithms in CIL relevant for stochastic optimisation. By choosing $\gradest = \nabla f$, the existing deterministic variant is achieved. Choosing a gradient estimator for $\gradest$, as described in the following sections, specific stochastic optimisation algorithms are achieved. Here, $P_k$ denotes a preconditioner and its use is illustrated later in Section~\ref{sec:pet}.
\label{table:CIL_algorithms}}
\end{table}

We consider two broad categories of algorithms for~\eqref{eq:optimisation_three_terms}: proximal gradient and primal-dual algorithms.

\paragraph{Proximal gradient algorithms}
When the composite term in~\eqref{eq:optimisation_three_terms} is missing, i.e., $h\circ\Kop=0$, we recover a simpler two-term optimisation problem solvable, for example, by the proximal gradient descent (PGD) algorithm (see algorithm~\ref{alg:PGD}). This flexible algorithm includes several other common algorithms as special cases, for example, if $g$ is the $\ell_{1}$ norm, i.e., $g=\|\cdot\|_{1}$, we recover the ISTA algorithm~\cite{Daubechies2004}. Note that in CIL, the names ISTA and PGD are used interchangeably; the algorithms are aliases of each other.  Additionally, if $g$ is an indicator function of a closed, non-empty, convex set, then the proximal of $g$ is the projection onto the set, and we recover the projected Gradient Descent (GD) algorithm. Furthermore, if $g=0$ as well, we obtain the simple case of one differentiable term solvable by GD (see algorithm~\ref{alg:GD}).

There is a range of parameter choices, which, when tuned correctly, can accelerate the gradient-based algorithms, including stochastic variants:
\begin{itemize}
\item \textbf{Step sizes} can be changed adaptively or remain constant across the iterations. By default, GD in CIL uses the Armijo rule, a backtracking line search algorithm to choose a step size that guarantees the iteration leads to a sufficient decrease in the optimisation objective (see algorithm 3.1 in~\cite{nocedal2006numerical}). CIL also has the Barzilai-Borwein step size method~\cite{bbstepsize} implemented. Custom decreasing step size rules are easy to implement in CIL.
\item \textbf{Momentum} based acceleration is implemented in Accelerated Proximal Gradient Descent (APGD, see algorithm~\ref{alg:FISTA}). Here, the user has a custom choice of momentum coefficient, $\beta_k$.  By default, the momentum parameter is calculated as $t_0=1$, $t_{k+1}=\frac{1+\sqrt{1+4t_k^2}}{2}$, $\beta_k=\frac{t_{k}-1}{t_{k+1}}$, giving the Fast Iterative Shrinkage-Thresholding Algorithm (FISTA)~\cite{beck2009fast}. This is Nesterov Momentum or Nesterov Accelerated Gradient method (NAG) and gives $\mathcal{O}\left(\frac{1}{k^2}\right)$ acceleration without any additional memory requirements in the case where $g=h=0$ ~\cite{Nesterov1983, Nesterov2004, Nesterov2004book, Nesterov2012, beck2009fast, Taylor2017}.  We refer the reader to a recent survey on accelerated methods for more details~\cite{Aspremont2021}.
\item \textbf{Preconditioning} the gradient estimator modifies the calculated gradient to improve convergence properties, especially in ill-conditioned problems. Preconditioning can be interpreted as applying a transformation to the gradient that accounts for the geometry of the problem, effectively rescaling directions in which the objective function changes more slowly. For example, in CIL, there are options for sensitivity preconditioning where each call to the preconditioner multiplies the gradient by ${(\bm{A}\bm{1})^{-1}}$,
where $\bm{1}$ is an object in the range of the operator filled with ones, see ~\cite{Kaufman1993, Clinthorne1993}.
The CIL framework is flexible to user-defined preconditioners, and the preconditioning strategies may include diagonal scaling, operator-based preconditioners, or learned preconditioners tailored to specific inverse problems, ~\cite{Qi2006}.

\end{itemize}

\paragraph{Primal-Dual Algorithms}
There are many primal-dual algorithms designed to solve~\eqref{eq:optimisation_three_terms}, each with its own set of strengths and weaknesses~\cite{Condat2013, Vu2011, Chen2016, Yan2018, Latafat2017, Salim2022}.

In the absence of a smooth function, i.e., $f=0$ in ~\eqref{eq:optimisation_three_terms}, we obtain the two-term problem suitable for the Primal-Dual Hybrid Gradient (PDHG) algorithm, also known in the literature as Chambolle-Pock, see~\cite{Chambolle2010, Esser2010, Fadili2011}. This is available in CIL, and we refer the reader to~\cite{Chambolle2010} for the algorithm.

For the three term problem~\eqref{eq:optimisation_three_terms}, the PD3O algorithm is implemented in CIL (see algorithm~\ref{alg:PD3O}). A recent study~\cite{Yan2024} provides convergence guarantees with relaxed step size conditions which are provably optimal. PD3O provides relaxed step size requirements, ease of tuning, low computational cost, and natural extensibility to stochastic optimisation. While we do not claim that this is the optimal algorithm for all use cases, its base implementation can serve as a valuable guide for implementing the above related algorithms and extensions, e.g., Bregman proximal operators~\cite{Jiang2022}.

\subsection{Stochastic Methods }

In this section, we extend the deterministic optimisation framework~\eqref{eq:optimisation_three_terms} to its stochastic counterpart.

In a first step, we assume that~\eqref{eq:optimisation_three_terms} can be \textbf{partitioned }into a finite sum structure:
\begin{equation}
\bx^{*} \in\argmin_{\bx\in\mathbb{X}} \sum_{i=0}^{N-1}f_i(\bx) + \sum_{m=0}^{M-1}g_m(\bx)+\sum_{t=0}^{T-1}h_t(\Kop_t\bx). \label{eq:stochastic_optimisation_three_terms}
\end{equation}

\noindent 
In particular, the finite dimensional space $\mathbb{Y}$ is now a Cartesian product of $(\mathbb{Y}_{t})_{t=0}^{T-1}$, i.e., $\mathbb{Y}:=\prod_{t=0}^{T-1}\mathbb{Y}t$ with elements $\by=(\by_0,…,\by_{T-1})$. Moreover, we define linear and bounded operators $\Kop_t:\mathbb{X}\rightarrow \mathbb{Y}_t$ where $\Kop_t x=(\Kop x)_tt$. Under these assumptions the composite term $h(\bx)$ is now represented as a sum of composite terms $\sum_{t=0}^{T-1}h_{t}(\Kop_{t}\bx)$. The sums of the functions $f_{i}:\mathbb{X}\rightarrow\mathrm{R}$ and $g_{m}:\mathbb{X}\rightarrow\mathrm{R}$ represent now the sum of smooth and proximable functions, respectively. As before, $f_{i}$, $i=0,\dots, N-1$ and $g_{m}$, $m=0, \dots, M-1$ are proper, closed and convex functions, and additionally, we assume that functions $f_{i}$ are $L_{i}$ smooth.

The finite sum design of~\eqref{eq:stochastic_optimisation_three_terms} is a key concept behind the genesis of stochastic optimisation, where update steps are computed using a \textbf{sample} of the above functions.

As a concrete example, in tomography reconstruction, see section~\ref{sec:ct} for more details, we can express the fidelity term, $f$, as a finite sum, by partitioning the acquisition data $\bb$ into a set of measured data $(\bb_{i})_{i=0}^{N-1}$, where each $\bb_i$ is often called a subset. A natural partition for tomography is to partition the data into subsets of projection angles.  We can define the corresponding forward operators, $\sysmat_{i}$. Then, for every $i\in{0, \dots, N-1}$, we  define $f_{i}(\bx) = \|\sysmat_{i}\bx-\bb_{i}\|^{2}$ and then

\begin{equation}
f(\bx) = \|\sysmat\bx - \bb\|^2 =
\left\|
\renewcommand{\arraystretch}{0.7}
\begin{pmatrix}
\sysmat_0 \bx - \bb_0 \\
\vdots \\
\sysmat_{N-1} \bx - \bb_{N-1}
\end{pmatrix}
\renewcommand{\arraystretch}{1.0}
\right\|^2
= \sum_{i=0}^{N-1} \|\sysmat_i \bx - \bb_i\|^2
= \sum_{i=0}^{N-1} f_i(\bx).
\label{eq:ls_finite_sum}
\end{equation}

Note, that in our framework, which is focused on inverse problems, we consider decomposing $f$ into a finite sum structure, $f = \sum_{i=0}^{N-1}f_i$ not, as in machine learning applications, where it is common to represent $f$ as an average of losses $f = \frac{1}{N}\sum_{i=0}^{N-1}f_i$~\cite{Bottou2018}.
Therefore, in CIL, slight modifications were sometimes required for the stochastic algorithms compared to the original literature. 

\subsubsection{Stochastic Gradient Descent}
\begin{table}[t]
\footnotesize
\begin{tabular}{p{2.8cm}|p{2.1cm}|p{7.5cm}|p{3cm}}
\hline
\textbf{Stochastic Gradient Estimator} &
\textbf{CIL Approx. Gradient Function} &
\textbf{Gradient calculation} &
\textbf{Initialisation} \\
\hline

Stochastic Gradient Descent (SGD) &
\texttt{SGFunction} &
\begin{minipage}[t]{7.5cm}
\raggedright
\[
\begin{aligned}
\gradest(\mathbf{x}) &= N\, \nabla f_{i_k}(\mathbf{x}) ,\\
\text{where } i_k & \text{ sampled from } \{0,1,\dots,N-1\}.
\end{aligned}
\]
\end{minipage}
&
Not applicable \\
\hline

Stochastic Average Gradient (SAG)~\cite{Schmidt2017} &
\texttt{SAGFunction} &
\begin{minipage}[t]{7.5cm}
\raggedright
\[
\begin{aligned}
\gradest(\mathbf{x}) &= 
\nabla f_{i_k}(\mathbf{x}) - \delta_{i_k}^k + \sum_{i=0}^{N-1}\delta_i^k ,\\[4pt]
\delta_j^{k+1} &=
\begin{cases}
\nabla f_j(\mathbf{x}), & \text{if } j = i_k, \\
\delta_j^k, & \text{otherwise},
\end{cases} \\[4pt]
\text{where } i_k & \text{ sampled from } \{0,1,\dots,N-1\}.
\end{aligned}
\]
\end{minipage}
&
Initialise $\delta_j^0 = 0$ or warm-start with $\nabla f_j(\mathbf{x}_0)$ for $j = 0,\dots, N-1$.\\
\hline

SAGA (SAG-amélioré)~\cite{Defazio2014} &
\texttt{SAGAFunction} &
\begin{minipage}[t]{7.5cm}
\raggedright
\[
\begin{aligned}
\gradest(\mathbf{x}) &=
N\bigl(\nabla f_{i_k}(\mathbf{x}) - \delta_{i_k}^k\bigr)
+ \sum_{i=0}^{N-1}\delta_i^k ,\\[4pt]
\delta_j^{k+1} &=
\begin{cases}
\nabla f_j(\mathbf{x}), & \text{if } j = i_k, \\
\delta_j^k, & \text{otherwise},
\end{cases} \\[4pt]
\text{where } i_k & \text{ sampled from } \{0,1,\dots,N-1\}.
\end{aligned}
\]
\end{minipage}
&
Initialise $\delta_j^0 = 0$ or warm-start with $\nabla f_j(\mathbf{x}_0)$ \\
\hline

Stochastic Variance Reduced Gradient (SVRG)~\cite{Johnson2013} &
\texttt{SVRGFunction} &
\begin{minipage}[t]{7.5cm}
\raggedright
\[
\begin{aligned}
\gradest(\mathbf{x}) &=
\begin{cases}
N\, \nabla f(\bar{\mathbf{x}}_k), & \text{if } k \bmod T = 0,\\[4pt]
N\bigl(\nabla f_{i_k}(\mathbf{x}) - \nabla f_{i_k}(\bar{\mathbf{x}}_k)\bigr)
+ \nabla f(\bar{\mathbf{x}}_k), & \text{otherwise},
\end{cases} \\[6pt]
\bar{\mathbf{x}}_k &=
\begin{cases}
\mathbf{x}, & \text{if } k \bmod T = 0,\\
\bar{\mathbf{x}}_{k-1}, & \text{otherwise},
\end{cases} \\[4pt]
\text{where } i_k & \text{ sampled from } \{0,1,\dots,N-1\}.
\end{aligned}
\]
\end{minipage}
&
Choose snapshot period $T$ \\
\hline

Loopless-SVRG (LSVRG)~\cite{Kovalev2020, Hofmann2015, Qian2021} &
\texttt{LSVRGFunction} &
\begin{minipage}[t]{7.5cm}
\raggedright
\[
\begin{aligned}
\gradest(\mathbf{x}) &=
\begin{cases}
N\, \nabla f(\bar{\mathbf{x}}_k), & \text{if } c < p,\\[4pt]
N\bigl(\nabla f_{i_k}(\mathbf{x}) - \nabla f_{i_k}(\bar{\mathbf{x}}_k)\bigr)
+ \nabla f(\bar{\mathbf{x}}_k), & \text{otherwise},
\end{cases} \\[6pt]
\bar{\mathbf{x}}_k &=
\begin{cases}
\mathbf{x}, & \text{if } c < p,\\
\bar{\mathbf{x}}_{k-1}, & \text{otherwise},
\end{cases} \\[4pt]
\text{where } i_k & \text{ sampled from } \{0,1,\dots,N-1\} \\
\text{and } c &\text{ sampled } U[0,1].
\end{aligned}
\]
\end{minipage}
&
Choose $p\in(0,1]$ (snapshot probability) \\
\hline

\end{tabular}

\caption{Stochastic gradient estimators in CIL.}
\label{tab:stochastic-estimators}
\end{table}

Consider the simplest case of~\eqref{eq:stochastic_optimisation_three_terms}, where $g_m, h_t \equiv 0$ and so
\begin{align}
x^* \in\argmin_{\bx\in\mathbb{X}} \sum_{i=0}^{N-1}f_i(\bx).
\end{align}
This is differentiable, so it can be solved using GD (algorithm~\ref{alg:GD}). For ``standard" gradient descent the \textit{full} gradient, over all functions $f_i$ and thus all of the data,  is calculated at each iteration, so in algorithm~\ref{alg:GD},  $\gradest(x) =\sum_{i=0}^{N-1} \nabla f_{i}(x)$.  The Stochastic Gradient Descent (SGD) algorithm proceeds as
\begin{equation}
\bx_{k+1} = \bx_{k} - \gamma\nabla f_{i_{k}}(\bx_{k}),
\label{eq:sgd_iteration}
\end{equation}
where $\gamma$ is a step size and $i_k$ is sampled from ${0,1,\dots,N-1}$. This is equivalent to algorithm~\ref{alg:GD} but  $\gradest(x)= N\nabla f_{i_{k}}(x)$ is an unbiased approximation of $\nabla f(x)$. This is the \textbf{key concept} in this paper and software: stochastic gradient algorithms can be implemented as deterministic algorithms with \textbf{a stochastic estimator of the gradient }plugged in.

The SGD algorithm is not a descent method, i.e., $f$ may not decrease in every iteration. However, in expectation, the iterates in~\eqref{eq:sgd_iteration} are moving in the direction of the negative full gradient. The key idea behind stochastic algorithms is that, although it may take more iterations to reach the same minimum, if the stochastic gradient estimator is computationally cheaper to calculate (in the best case, an $N$th of the computational cost), the stochastic algorithm will overall decrease the computational cost of the optimisation.

For a concrete example, consider the $f_i$ described in Equation~\ref{eq:ls_finite_sum}. In each iteration of a deterministic algorithm, calculating  $\gradest(x) =\sum_{i=0}^{N-1} \nabla f_{i}(x)$,  requires the computation of the forward operation $A_ix$ (and usually also its adjoint) for all $i \in \{0,…,N-1\}$ performing computations on all the data in each iteration, while the SG stochastic estimator, $\gradest(x)= \nabla f_{i_{k}}$, performs calculations on a subset of the data, using just one sampled $A_{i_k}$ (and corresponding data $b_{i_k}$).

To ensure a fair comparison between the deterministic and stochastic approaches, convergence is generally plotted and discussed in terms of ``data passes" rather than iterations, where one data pass corresponds seeing, in expectation, every subset of data. For SGD, there are in expectation $N$ iterations in one data pass. This is comparative to the ``epoch" commonly used in the Machine Learning community.

\subsubsection{Variance-Reduced Gradient Methods}

Table~\ref{tab:stochastic-estimators} gives all the stochastic estimators currently in CIL. The first row concerns the SGD estimator. As discussed above, using this stochastic estimator in algorithm~\ref{alg:GD} gives the iterations in equation~\ref{eq:sgd_iteration}.

There are multiple examples in the literature of stochastic methods which we can now write in the form of a stochastic gradient estimator  $\gradest$, which can be plugged into deterministic algorithms to create the stochastic algorithms.
With a combination of averaging and full gradient computations, these methods use a stochastic estimate of the full gradient $\gradest\approx\nabla f(x_{k})$, with a smaller error than in SGD and thus faster convergence rates and easier to tune step sizes compared to SGD. The drawback is that this improvement is realised with an increased computational cost and/or additional storage requirements compared to SGD. Here we discuss the methods currently in CIL, and listed in table \ref{tab:stochastic-estimators}.
\begin{itemize}
\item \textbf{SAG: }The first variance-reduced method is the Stochastic Average Gradient (SAG) method~\cite{Schmidt2017}, which is a stochastic extension of the Increment Aggregated Gradient method,~\cite{Blatt2007}. The SAG estimator, see Table~\ref{tab:stochastic-estimators}, requires storing a table of $N$ stochastic gradients, $\delta_i^{k}$, initialised to zero, for example, that are used in every iteration. In every step, only one position of the table is updated by computing the gradient of only one randomly sampled function. In this way, we have a reduced computational cost per iteration similar to the SGD algorithm.  SAG is the first algorithm that enjoys similar convergence rates to full gradient methods.  The initial convergence analysis in~\cite{Schmidt2017} considered only smooth objectives. Later, the analysis of SAG was extended to non-smooth objectives with proximal support~\cite{Driggs2022}.

\item \textbf{SAGA: }In~\cite{Defazio2014}, an ``improved" version of SAG is presented, known as SAGA, where the
additional `A` stands for the French \emph{am\'elior\'e}.
Compared to SAG and under our convention discussed above, the SAGA gradient estimator in table~\ref{tab:stochastic-estimators} has a weight $N$ in front of the difference between the stochastic gradient and the stored gradient $\delta^k_{i}$. The SAG calculation uses a biased estimator of the full gradient of $f$, whereas the change in scaling for SAGA makes it an unbiased stochastic estimator.

\item \textbf{SVRG: }Storing a table of gradients $\delta_i^{k} \in \mathbb{X}$ for all $ i\in\{0,1,\dots, N-1\}$ may require a large amount of storage for inverse imaging problems, e.g., for tomographic reconstruction applications~\cite{Karimi2016}. In~\cite{Johnson2013}, the authors introduce the Stochastic Variance Reduced Gradient (SVRG) algorithm, initially for smooth objectives and later for optimisation problems with proximal support in~\cite{Xiao2014}.  The SVRG algorithm stores a reference image, $\Bar{\bx}_k$, also known as a ``snapshot", and the full gradient computed at the reference, which are updated every $T$ iterations. When the gradient estimator for SVRG is called, one of two things happens: every $T$ iterations a snapshot is taken, and its full gradient is calculated, and in every other case, a gradient of just one subset is calculated and an approximate gradient computed using the information of the previously computed full gradient. Although memory requirements are decreased compared to SAG and SAGA, an additional parameter $T$ is introduced, which requires careful tuning. In terms of computational cost, in SVRG, two gradients are computed in every step in addition to the computation of the full gradient at every $T$ iterations.

\item \textbf{LSVRG: }A common variation on SVRG is the  Loopless-SVRG (LSVRG)~\cite{Kovalev2020, Hofmann2015, Qian2021} where the snapshot every $T$ iterations is replaced by a probability, $p$, that a snapshot is taken and a full gradient returned.  At every update, the reference point is updated with probability $p$ and the full gradient is computed. Otherwise, it remains unchanged, and we avoid computing the full gradient. 

\end{itemize}

There are other stochastic gradient estimators in the literature, for example, the  StochAstic Recursive grAdient algoritHm (SARAH) is a biased version of SVRG, which has also been extended to Loopless SARAH in~\cite{Bingcong2019} and with proximal support in~\cite{Nhan2019}. The latter algorithms are not implemented in the current version of CIL. However, our library provides the necessary building blocks to construct them with ease.

Another option for adding stochasticity is to apply random selection for the last term in~\eqref{eq:stochastic_optimisation_three_terms}, i.e., the sum of composite terms. For this, the optimisation problem can be written as
\begin{equation}
\bx^{*} \in\argmin_{\bx\in\mathbb{X}} f(\bx) + g(\bx)+\sum_{t=0}^{T-1}h_{t}(\Kop_{t}\bx). \label{eq:stochastic_optimisation_three_terms_c}
\end{equation}

The Stochastic Primal-Dual Hybrid-Gradient (SPDHG)~\cite{Chambolle_2018_SPDHG} algorithm adds stochasticity not through an approximate gradient but rather by decomposing the space $\mathbb{Y}$ into a product space $\prod_{t=0}^{T-1}\mathbb{Y}_i$ and updating only one subset of the dual variable at each iteration.  SPDHG is implemented in CIL, in the case of $f=0$, and we use it as a comparison method.

\subsection{The plug-and-play framework: Combining Stochastic Gradient Estimators with Deterministic Algorithms}

\begin{table}[h!]
\centering
\begin{tabular}{|>{\centering}m{3.0cm}|>{\centering}m{2.3cm}|>{\centering}m{2.6cm}|>{\centering}m{2.6cm}|>{\centering\arraybackslash}m{2.6cm}|}
\hline

    \diagbox[width=3.35cm, height=1.75cm]{\hspace{-0.1cm}\textbf{Functions}}{\vspace{0.5cm}\textbf{Algorithms}} 
    & \code{GD} (Alg. 
   ~\ref{alg:GD}) & \code{PGD} (Alg. 
   ~\ref{alg:PGD}) & \code{APGD} (Alg.~\ref{alg:FISTA}) & \code{PD3O} (Alg.~\ref{alg:PD3O})
    \\ \hline
    \code{SGFunction} & SGD~\cite{Robbins1951} & ProxSGD ~\cite{JMLR:v10:duchi09a} & Acc-ProxSGD ~\cite{pmlr-v124-zhou20a} & PD3O-SGD ~\cite{Salim2022}\\ \hline
    \code{SAGFunction} & SAG~\cite{Schmidt2017} & ProxSAG~\cite{Driggs2022} & Acc-ProxSAG ~\cite{Driggs2020} & PD3O-SAG ~\cite{Salim2022}\\ \hline
    \code{SAGAFunction} & SAGA~\cite{Defazio2014} & ProxSAGA~\cite{Defazio2014} & Acc-ProxSAGA ~\cite{Driggs2020} & PD3O-SAGA ~\cite{Salim2022}\\ \hline
    \code{SVRGFunction} & SVRG~\cite{Johnson2013} & ProxSVRG~\cite{Xiao2014} & Acc-ProxSVRG ~\cite{Driggs2020} & PD3O-SVRG ~\cite{Salim2022}\\ \hline
    \code{LSVRGFunction} & LSVRG~\cite{Kovalev2020} & ProxLSVRG~\cite{Qian2021} & Acc-ProxLSVRG ~\cite{Driggs2020} & PD3O-LSVRG ~\cite{Salim2022} \\ \hline
\end{tabular}
\caption{Overview of stochastic optimisation algorithms. Each row corresponds to a CIL \code{Function} class that implements a specific stochastic gradient approximation, and each column represents a specific deterministic CIL \code{Algorithm}.}
\label{table:CIL_PnP}

\end{table}

In the previous section, we discussed how GD, defined in algorithm~\ref{alg:GD},  combined with the stochastic estimators in~\ref{tab:stochastic-estimators}, creates SGD, SAG, SAGA,  SVRG and LSVRG algorithms.

In Table~\ref{table:CIL_PnP}, we present how different stochastic estimators, matched with the deterministic algorithms in Table ~\ref{table:CIL_algorithms}, can combine to give a range of stochastic algorithms.

Additionally, in the case of equation~\eqref{eq:stochastic_optimisation_three_terms} where $M=1$ and $h\equiv 0$ we can combine stochastic estimators with PGD (algorithm~\ref{alg:PGD}) and APGD (algorithm~\ref{alg:FISTA}) to give proximal stochastic methods and accelerated methods, with references given in the table.

Finally, we consider again equation~\eqref{eq:stochastic_optimisation_three_terms} when the composite term $h(\Kop\bx)$ is present. For simplicity, we consider the case where $M=T=1$. We can expand the PD3O algorithm (algorithm~\ref{alg:PD3O}) by again replacing the gradient with a stochastic gradient estimator.
In~\cite{Yurtsever2016}, the three operator splitting method is presented for the case of $\Kop=\mathbb{I}$ and extended for the $\Kop\neq\mathbb{I}$ in~\cite{Zhao2019} for stochastic gradients. Recently, in~\cite{Salim2022}, PD3O, alongside other algorithms,  are used with stochastic variance reduced estimators such as SAGA and SVRG.

%% file: version3/03_implementation.tex
\section{Software Solution Step-by-step} \label{sec:implementation}

This section gives a walkthrough of our software solution.   The next two sections then demonstrate its application across a range of imaging tasks on real-world datasets, validating the software's versatility and practical relevance. In the supplementary material section, we provide an additional example on a non-imaging task.

A stochastic algorithm described in the previous section can be fully defined once we determine all of the following:

\begin{enumerate}
    \item \textbf{Data Partitioner}: How to partition the dataset in subsets;
    \item \textbf{Sampler}: Sampling strategy for the stochastic gradient estimator to select the subset to operate on;
    \item \textbf{The stochastic gradient estimator, $\gradest$};
    \item \textbf{Deterministic Algorithm}: a gradient type deterministic algorithm (GD, PGD or APGD) or the PD3O algorithm; 
    \item \textbf{Optionally} define pre-conditioning, and step size strategy and for APGD choose the acceleration technique such as momentum.
\end{enumerate}

Our software solution enables the definition of all the 5 points above, providing new CIL functionality for flexible data partitioning for CT and general sampling strategies, together with a flexible and extensible framework to integrate stochastic gradient estimators into deterministic algorithms, enabling the construction of stochastic optimisation methods.

\begin{enumerate}
\item \textbf{Data Partitioner:}
 In the presented stochastic framework, we focus on the splitting of the differentiable term (the first term) in~\eqref{eq:optimisation_three_terms} into a sum of $N$  differentiable terms $\{f_i\}$, see the first sum in~\eqref{eq:stochastic_optimisation_three_terms}. For the stochastic algorithms to provide a meaningful speed up, each $f_i$ must operate on just a subset of the data at each iteration, and thus partitioning $f$ into a finite sum requires us to partition the data. 
 
 Exactly how this is done depends on the specific data domain and the form of $f$.
 
 We will present the X-ray CT data partitioning, implemented in CIL, in section~\ref{subsec:XCT-partitioning} and PET data partitioning to a lesser extent in section~\ref{subsec:PET-partitioning}.

 For the remainder of this section, we define the list of functions obtained by the data partitioner:

 \begin{center}
\begin{center}
\begin{tcolorbox}[
    enhanced,
    attach boxed title to top center={yshift=-2mm},
    colback=darkspringgreen!20,
    colframe=darkspringgreen,
    colbacktitle=darkspringgreen,
    title= Partitioned functions,
    text width = 9.0cm,
    fonttitle=\bfseries\color{white},
    boxed title style={size=small,colframe=darkspringgreen,sharp corners},
    sharp corners,
]
\begin{minted}{python}
list_fi = [f_1, ..., f_N]
\end{minted}
\end{tcolorbox}
\end{center}

\end{center}

\item \textbf{Sampler:}

All the stochastic gradient estimators, defined in table~\ref{tab:stochastic-estimators}, for each evaluation require an $i_k$ sampled from $\{0,1,...,N-1\}$. In CIL, we have defined a \code{Sampler} class responsible for generating indices $i_k$ at each iteration.

In the code snippet below, we present several sampling strategies available to users. 
Scatter plots and histograms illustrating the sampled indices across iterations for each strategy are shown in Figure~\ref{fig:scatter-hist-all-methods}.

\begin{center}
\begin{center}
\begin{tcolorbox}[
    enhanced,
    attach boxed title to top center={yshift=-2mm},
    colback=darkspringgreen!20,
    colframe=darkspringgreen,
    colbacktitle=darkspringgreen,
    title= Sampling Strategies,
    text width = 11.0cm,
    fonttitle=\bfseries\color{white},
    boxed title style={size=small,colframe=darkspringgreen,sharp corners},
    sharp corners,
]
\begin{minted}{python}
N = 30 # The number of indices to sample from
sampler_seq = Sampler.sequential(N)
sampler_s06 = Sampler.staggered(N, stride=6)
sampler_s12 = Sampler.staggered(N, stride=12) 
sampler_rwr = Sampler.random_with_replacement(N, seed=6)
sampler_rnu = Sampler.random_with_replacement(N, seed=6,
                                prob= [0.9/20]*20 + [0.1/10]*10)

sampler_rwo = Sampler.random_without_replacement(N, seed=6)
sampler_her = Sampler.herman_meyer(N)
def my_sampling(iteration_number):
    return np.floor(N*np.sin((np.pi/(2*N))*iteration_number)**2)
sampler_cus = Sampler.from_function(num_indices=N, function=my_sampling)
\end{minted}
\end{tcolorbox}
\end{center}

\end{center}

\begin{figure}[h!]
\centering
\includegraphics[width=\linewidth]{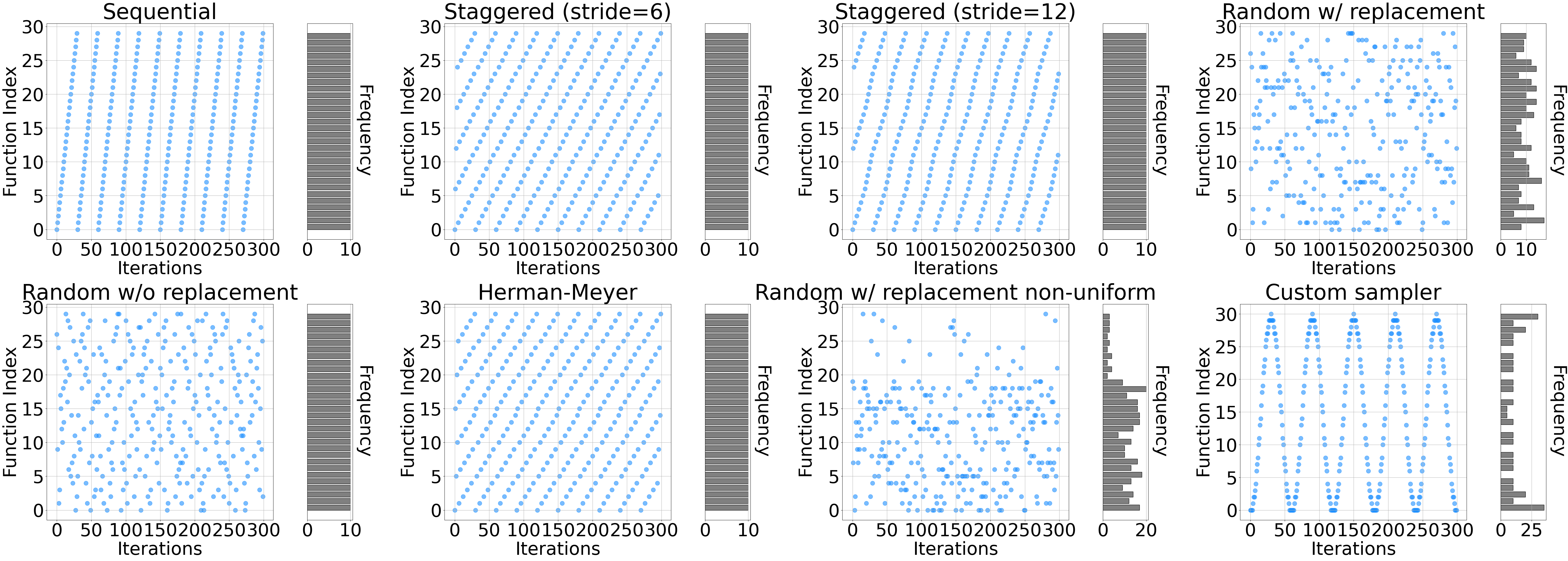}
\caption{Plots to show different sampling strategies for 30 indices over 300 calls to the sampler (10 data passes). Next to the scatter plots, the frequency of each index over the 300 calls is plotted. }
\label{fig:scatter-hist-all-methods}
\end{figure}

These examples include random sampling with or without replacement, implemented via \code{random\_with\_replacement} and \code{random\_without\_replacement}, respectively. In the case of sampling with replacement, non-uniform probabilities can be specified when initialising the sampler, provided that the probabilities are positive and sum to one.
For sequential sampling (\code{sequential}), indices are selected in order. In staggered sampling (\code{staggered}), a variation of sequential sampling, indices are selected cyclically using a predefined stride.
The Herman-Meyer sampling method, proposed in~\cite{Herman1993} and used in~\cite{ImrajHerman}, is primarily used in tomography applications. It requires the number of functions $f_i$ to be non-prime, as it is based on prime number decomposition.

In most cases, the sampler is passed straight to the stochastic gradient function or algorithm. However, the user can sample from it directly using the \texttt{next()}\footnote{\url{https://docs.python.org/3/library/functions.html\#next}} method. 
 
Our \code{Sampler} class shares similarities with sampler modules found in other open source libraries. For example, in PyTorch, the \code{sampler} is used within the \code{DataLoader} class to enable efficient batching and shuffling of data. Similarly,  JaxOpt~\cite{Blondel2021}, an open source optimisation library for JAX, provides a data iterator that is directly attached to the corresponding algorithm.
In contrast, our approach decouples the \code{Sampler} from both the \code{Algorithm} class and the input data. Instead, the \code{Sampler} is attached directly to the stochastic gradient estimator, offering greater flexibility and modularity.

\item \textbf{The Stochastic Gradient Estimator, $\gradest$
\label{subsec:ApproximateGradientSumFunction}:} 
The next step is to choose a stochastic gradient estimator to approximate the gradient of the function $f$, the sum of the partitioned functions $\{f_i\}$.  Rather than creating separate CIL \texttt{Algorithm} classes for each method listed in Table~\ref{tab:stochastic-estimators}, we introduce a new abstract base class: \texttt{ApproximateGradientSumFunction}, a subclass of the CIL \texttt{SumFunction} class. This class acts as a core component for stochastic optimisation in CIL, supporting the implementation of various algorithms that rely on different gradient approximation strategies.

The class is initialised with a list of functions, $f_i$, and an instance of a stochastic sampler, \texttt{sampler\_xyz}, which could be any of the ones described above or a custom user-defined sampler.  When the function is called, at point $\bx$,  it returns the sum over all $f_i(\bx)$. Its key feature is the \texttt{gradient} class method, which calls an abstract \texttt{approximate\_gradient} method.  Derived classes, which behave like CIL \code{Function} classes, can then implement the stochastic approximate gradients defined in Table~\ref{tab:stochastic-estimators} in the \texttt{approximate\_gradient} method.  

For the user, this machinery is hidden, and they simply pass their Python list \code{list_fi} of functions $\{f_0, f_1, \dots, f_{N-1}\}$, and chosen stochastic \texttt{sampler} to their chosen stochastic gradient estimator, e.g. for SG
\begin{center}
\begin{tcolorbox}[
    enhanced,
    attach boxed title to top center={yshift=-2mm},
    colback=darkspringgreen!20,
    colframe=darkspringgreen,
    colbacktitle=darkspringgreen,
    title=Initialising the \code{SGFunction},
    text width = 11.5cm,
    fonttitle=\bfseries\color{white},
    boxed title style={size=small,colframe=darkspringgreen,sharp corners},
    sharp corners,
]
\begin{minted}{python}
f_SG = SGFunction(list_fi, sampler=sampler_xyz)
\end{minted}
\end{tcolorbox}
\end{center}

\item \textbf{Deterministic Algorithm:}
Depending on the properties of the objective function defined in~\eqref{eq:stochastic_optimisation_three_terms}, users can choose from:
\begin{itemize}
    \item Algorithm~\ref{alg:GD}, GD,  for smooth objectives,
    \item Algorithms~\ref{alg:PGD} and \ref{alg:FISTA} PGD, and APGD for objectives that can be split into a smooth term and a proximable term, 
    \item Algorithm~\ref{alg:PD3O}, PD3O, objectives that can be split into a smooth term, a proximable term and an additional composite term with a proximable convex conjugate.
\end{itemize} 

\end{enumerate}

In the two code snippets below, we begin with a CIL function \texttt{f}, which has already been partitioned into a list of $N$ CIL functions $\{f_i\}$. For the deterministic algorithm, we choose GD with a fixed step size \texttt{gamma}, which depends on the Lipschitz constant of the function of \texttt{f}. For the stochastic algorithm, we first set up the approximate gradient function, taking in the list of functions and any of the samplers demonstrated above. This is passed to the SGD algorithm alongside a fixed step size \texttt{gamma}, which depends on the Lipschitz constant of the approximate gradient function \texttt{f\_SG}, i.e., $\gamma = \frac{1}{\sum_{i=0}^{N-1}L_{i}}$. Both methods are then run for $K$ iterations in the same manner. Both algorithms are by default initialised with zero images, but custom initialisations can be used.

\begin{center}
\begin{tcolorbox}[
    enhanced,
    attach boxed title to top center={yshift=-2mm},
    colback=darkspringgreen!20,
    colframe=darkspringgreen,
    colbacktitle=darkspringgreen,
    title=Gradient Descent (GD),
    text width = 11.5cm,
    fonttitle=\bfseries\color{white},
    boxed title style={size=small,colframe=darkspringgreen,sharp corners},
    sharp corners,
]
\begin{minted}{python}
gamma = 1/f.L
gd = GD(f=f, step_size=gamma)
gd.run(K)
\end{minted}
\end{tcolorbox}
\end{center}

\begin{center}
\begin{tcolorbox}[
    enhanced,
    attach boxed title to top center={yshift=-2mm},
    colback=darkspringgreen!20,
    colframe=darkspringgreen,
    colbacktitle=darkspringgreen,
    title=$\textbf{SGD} \mbox{ := }$ \code{GD} + \code{SGFunction},
    text width = 11.5cm,
    fonttitle=\bfseries\color{white},
    boxed title style={size=small,colframe=darkspringgreen,sharp corners},
    sharp corners,
]
\begin{minted}{python}
f_SG = SGFunction(list_fi, sampler=sampler_xyz)
gamma = 1/f_SG.L
sgd = GD(f=f_SG, step_size=gamma)
sgd.run(K)
\end{minted}
\end{tcolorbox}
\end{center}

With minor adjustments, such as changing the stochastic gradient estimator or the sampling strategy, new methods and algorithms can be tested with minimal effort. This highlights the plug-and-play nature of our software, which supports flexible experimentation with minimal code changes.

%% file: version3/04_example_ct.tex
\section{Stochastic Optimisation for  Parallel-beam CT reconstruction} \label{sec:ct}

In this section, we demonstrate the application of our novel stochastic optimisation framework to a real-world imaging task: Computed Tomography (CT) reconstruction. Specifically, we utilise a real 3D CT dataset acquired with parallel beam geometry. While parallel beam CT data is traditionally reconstructed slice by slice, we choose to reconstruct the full volume in one go to leverage 3D total variation (TV) regularisation across slices. This approach provides additional regularisation that would not be achieved if slices were treated independently.

Real-world 3D CT datasets present unique challenges due to their high resolution and large data volumes, necessitating efficient optimisation techniques to ensure accurate and scalable reconstruction. To address this, we introduce specialised utilities for partitioning CT datasets, facilitating a user-friendly workflow for large-scale tomographic reconstruction.

In this section, we first illustrate how to partition data for a CT use case and then show how to create the required approximate gradient sum functions in CIL. In this use case, we aim to minimise an objective that includes a non-differentiable term, here the TV regularisation. We therefore use these approximate gradient sum functions within proximal or primal-dual methods, which are able to handle the non-differentiable component.
We then demonstrate the performance of these stochastic algorithms, providing comparisons between a range of stochastic and deterministic methods for reconstructing real 3D data. We also include comparisons across different numbers of subsets used in the stochastic methods.

\subsection{Dataset}

We use the X-ray CT measurements of a 10440 NMC532 Li ion battery~\cite{Vamvakeros2021} that were performed at beamline I12 of the Diamond Light Source using a 100 keV monochromatic X-ray beam~\cite{Drakopoulos2015}. The dataset is open-sourced and is available through the Battery Imaging Library (BIL), \cite{Docherty2025}. Radiographs were recorded with an X-ray imaging camera (PCO.edge), and the pixel size was 7.91 $\mu$m. Each scan consisted of 1800 projections (radiographs) covering an angular range of 0$\degree$--180$\degree$ with an acquisition time per frame of 8 ms. The size of the detector is 1560$\times$1560. Flat-field (50) and dark current (50) images were also collected before the CT measurements and were used to normalise the acquired radiographs before the tomographic reconstruction. In addition, we applied pre-processing routines,~\cite{Vo2021}, to remove vertical stripes in the sinogram space, which typically result in ring artefacts in the reconstruction space. For our numerical experiments presented below, we use a cropped dataset of 900 projection angles and 300 vertical slices.

\begin{figure}[tb]
\centering
\includegraphics[width=0.6\linewidth]{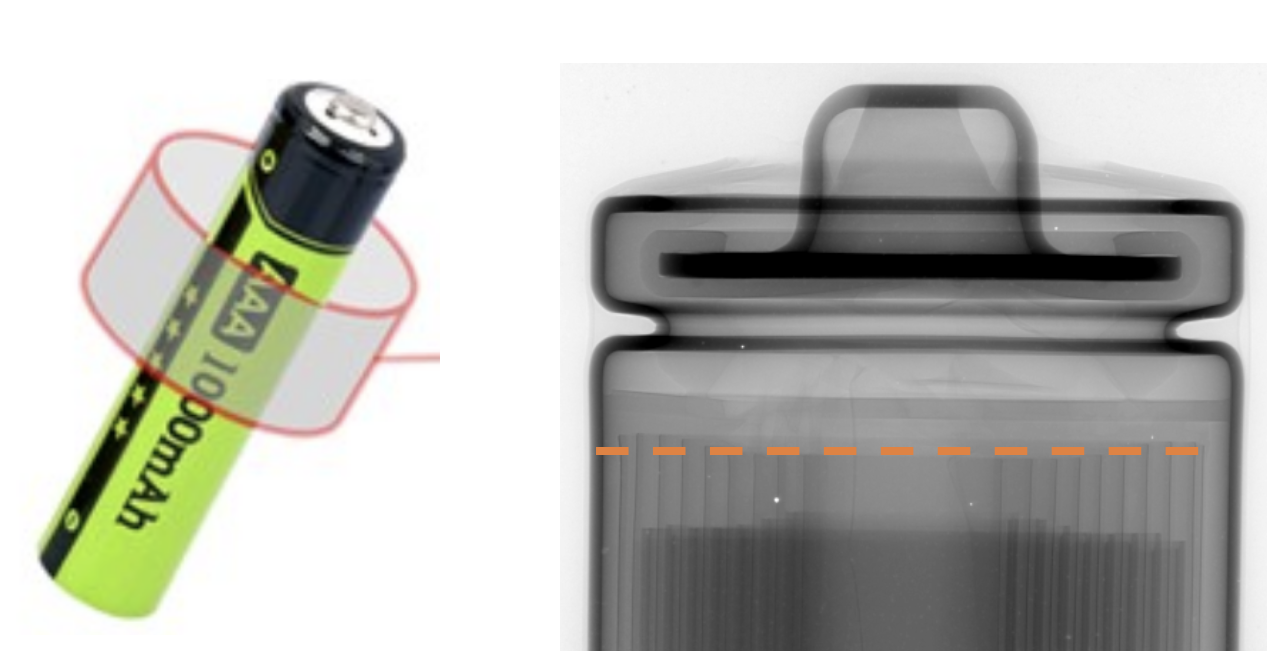}
\caption{Lithium-ion 10440 NMC532 battery, discharged at Capacity rate of 20 (left), and the radiograph for the top part of the battery (right).}
\label{fig:CT_battery_radiograph}
\end{figure}

\subsection{Methods}

The optimisation problem that we consider for the CT reconstruction is Total Variation regularisation with a non-negativity constraint. The acquired data is denoted by $\bb$, see Figure~\ref{fig:battery_sinogram}, and $\bA$ is a discretised version of the Radon transform obtained via either the ASTRA and TIGRE libraries with GPU acceleration, which we have wrapped into the \code{ProjectionOperator} class in CIL.

\begin{equation}
\argmin_{\bx}\frac{1}{2}\|\bA\bx-\bb\|_2^{2} + \alpha \mathrm{TV}(\bx) + \mathbb{I}_{\{\bx>0\}}(\bx).
\label{eq:CT_deter_tv_recon}
\end{equation}

Splitting the least squares function into a sum, as in~\eqref{eq:ls_finite_sum}, gives us the form
\begin{equation}
\argmin_{\bx}\frac{1}{2}\sum_{i=0}^{N-1}\|\bA_{i}\bx-\bb_{i}\|_2^{2} + \alpha \mathrm{TV}(\bx) + \mathbb{I}_{\{\bx>0\}}(\bx).
\label{eq:CT_stoch_tv_recon}
\end{equation}

The regularisation parameter $\alpha$ is chosen to be 0.75 for all experiments.

\subsubsection{Partitioning of Tomographic Datasets}\label{subsec:XCT-partitioning}\noindent  A key part of the stochastic framework, presented in this paper, is a partitioning strategy for tomography data that takes into account the underlying CT geometry. In the base CIL \texttt{AcquisitionData} class, we introduced a \texttt{}{partition} method that partitioned the data into a list of acquisition data, returning a \code{BlockDataContainer}. The data is split into a user-defined number of partitions, according to projection angles, with different angles being assigned to different groups according to a user-defined mode.  This code snippet demonstrates the user code for splitting the acquisition data, \code{b}, into a list of data \code{list_bi} of length $N$, with angles assigned to the $N$ partitions in a staggered way ( data from projection angles $l, l+N, l+2N, l+3N,…$ assigned to the $l$th partition.
\begin{center}
\begin{tcolorbox}[
enhanced,
attach boxed title to top center={yshift=-2mm},
colback=darkspringgreen!20,
colframe=darkspringgreen,
colbacktitle=darkspringgreen,
title=Data Partitioning Method,
text width = 9cm,
fonttitle=\bfseries\color{white},
boxed title style={size=small,colframe=darkspringgreen,sharp corners},
sharp corners,
]
\begin{minted}{python}
list_bi = b.partition(N, mode='staggered')
\end{minted}
\end{tcolorbox}
\end{center}

\begin{figure}[tb]
\vspace{0.5cm}
\begin{subfigure}[t]{0.45\textwidth}
\centering
\includegraphics[width=\textwidth]{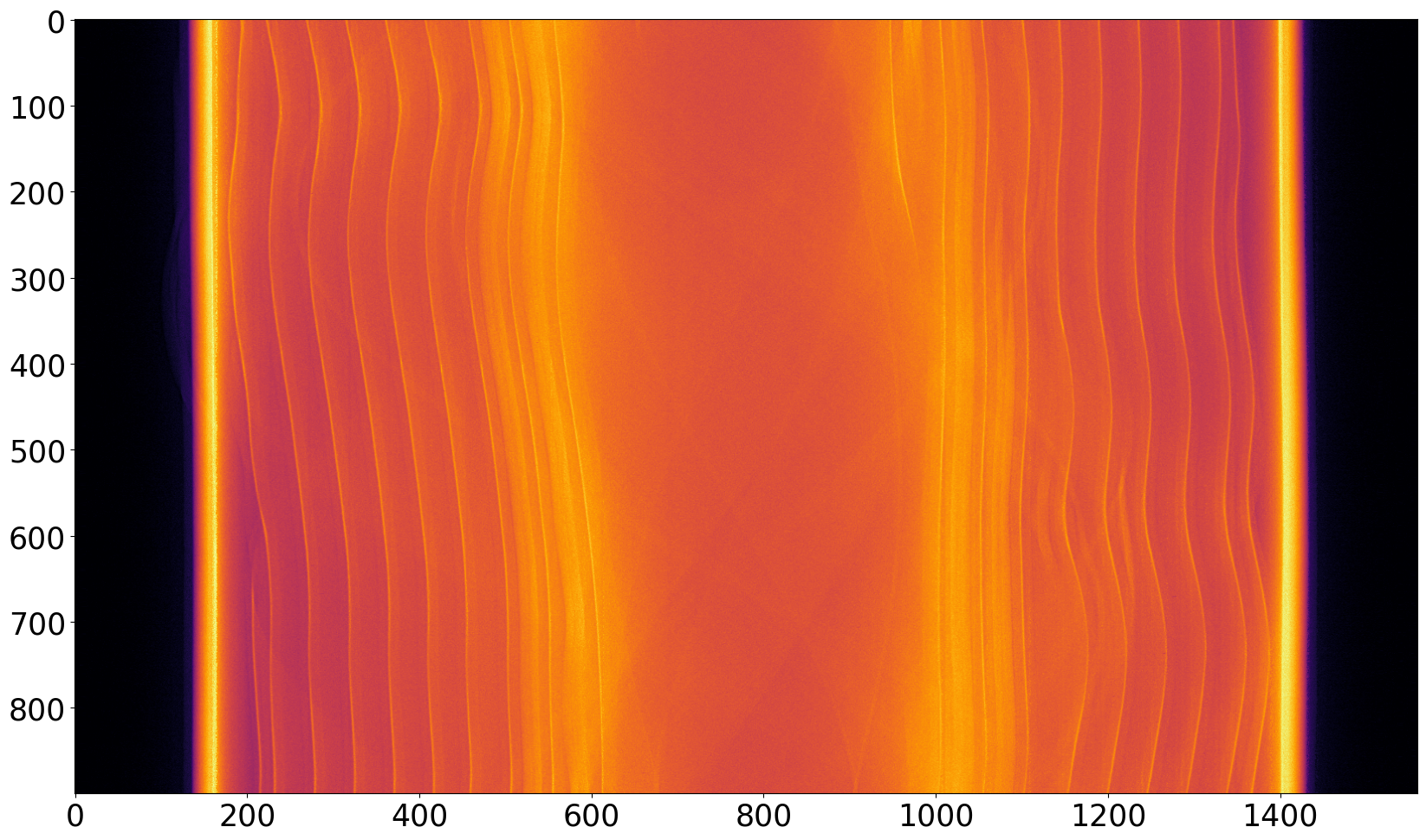}
\caption{Full sinogram}
\label{fig:battery_sinogram}
\end{subfigure}
\begin{subfigure}[t]{0.45\textwidth}
\centering
\includegraphics[width=\textwidth]{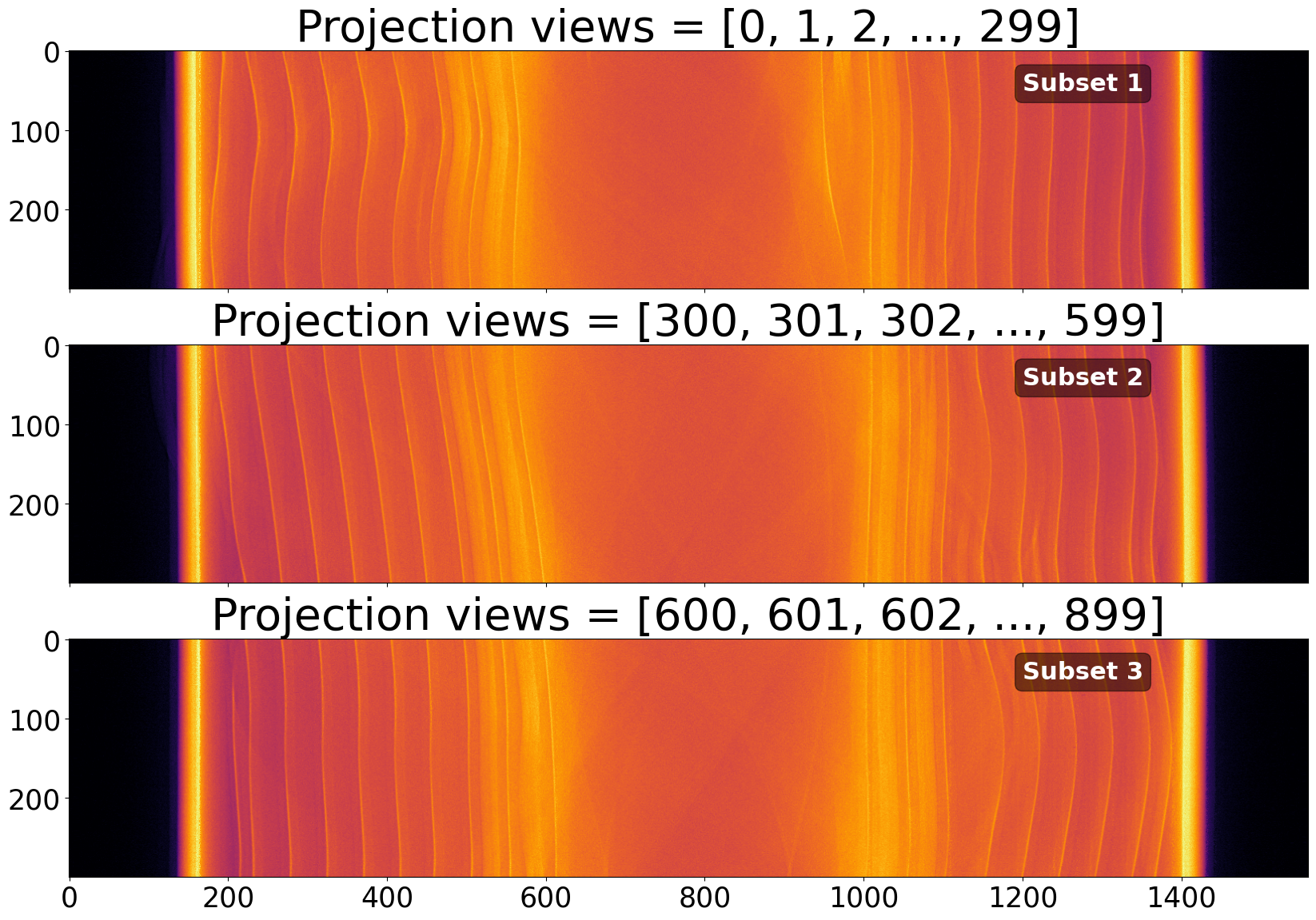}
\caption{Sequential Sampling}
\label{fig:battery_sequential_sinogram}
\end{subfigure}
\begin{subfigure}[t]{0.45\textwidth}
\centering
\includegraphics[width=\textwidth]{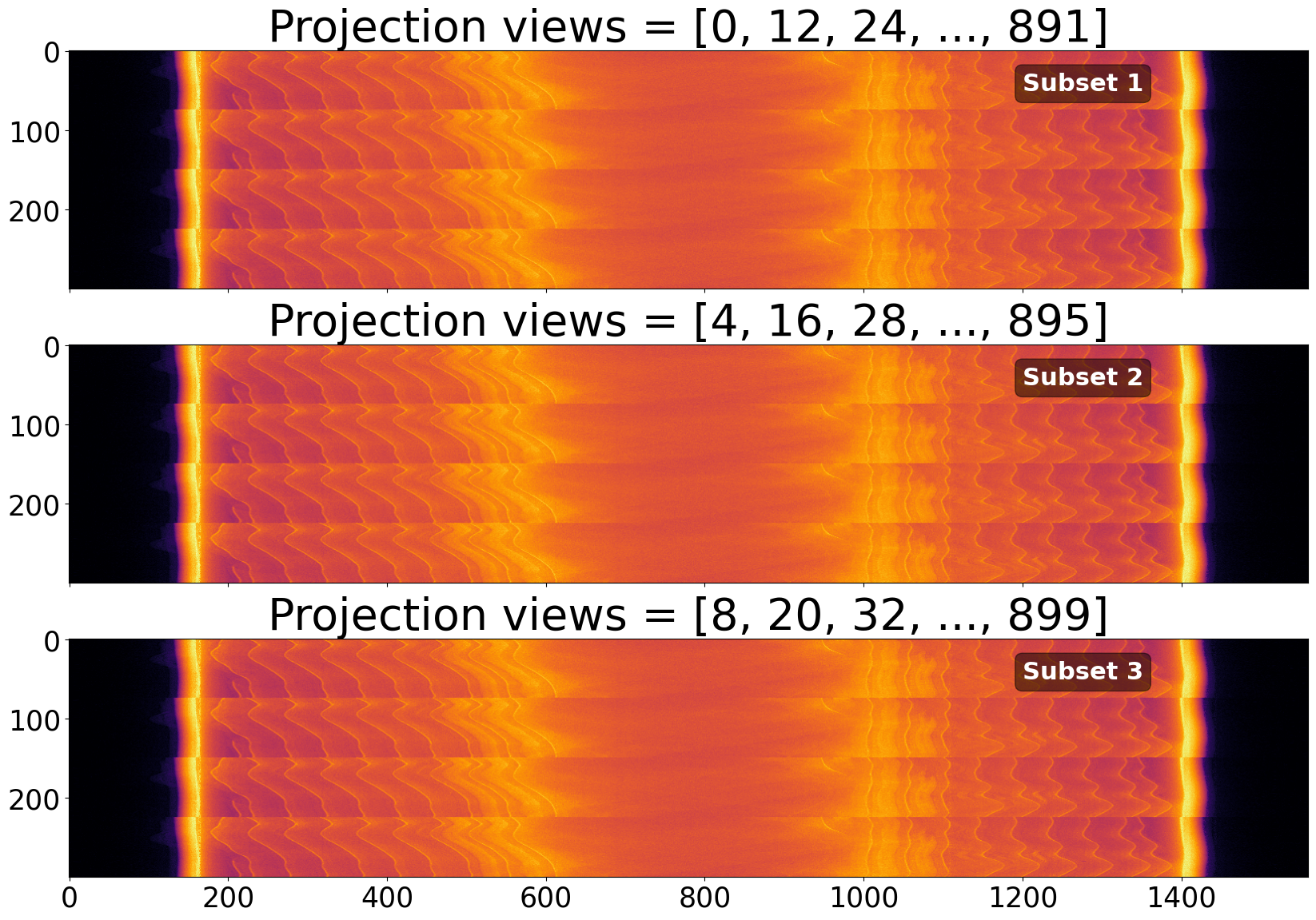}
\caption{Staggered Sampling}
\label{fig:battery_sequential_sinogram}
\end{subfigure}
\begin{subfigure}[t]{0.45\textwidth}
\centering
\includegraphics[width=\textwidth]{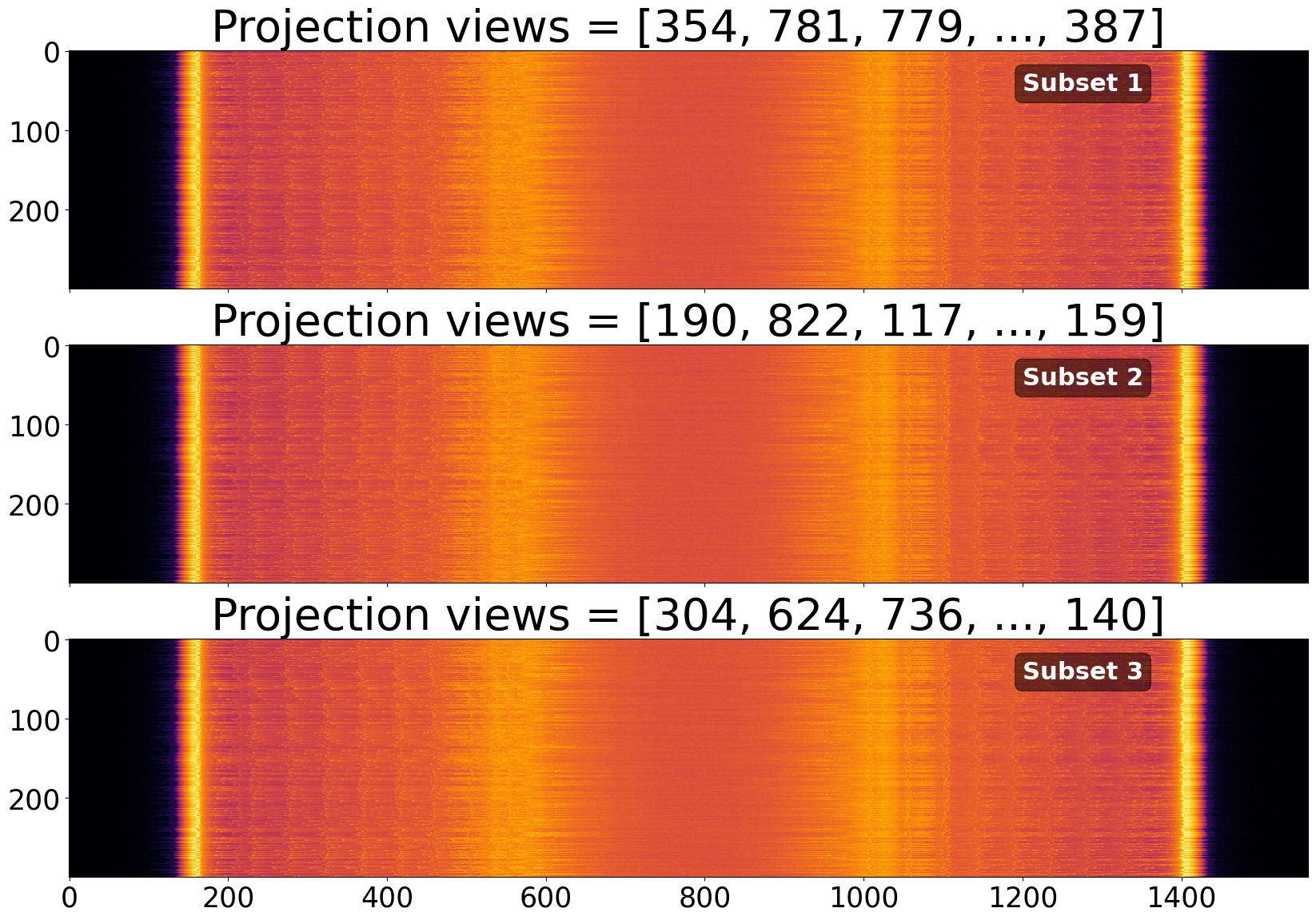}
\caption{Random Sampling}
\label{fig:battery_sequential_sinogram}
\end{subfigure}
\caption{Different data partitioning modes for the Lithium-ion battery sinogram. In each example, we split the data into 3 subsets of 300.  We choose to show the partitioned datasets as 2D sinograms, but remind readers that this is 3D data.  }
\label{fig:CT_partition_mode}
\end{figure}

In Figure~\ref{fig:CT_partition_mode}, we present different modes for partitioning tomography data. We plot in 2D the projection views in each partition to demonstrate \texttt{staggered}, \texttt{sequential} and \texttt{random} partitioning.  Note that, due to the symmetry of our sample, the individual partitioned sinograms for staggered and random all look similar to each other; each partition has a broad range of information, whereas for the sequentially partitioned data, each partition is quite different. In general, it is advisable to partition the data such that each partition contains information about all of the phantom, i.e. staggered or random.

We define the corresponding projection operators $(\bA_{i})_{i=0}^{N-1}$ with the geometric information associated with each $\bb_{i}$. The CIL projection operators can take a \code{BlockDataContainer} acquisition geometry and return a \code{BlockOperator}.  In addition, we construct a list of functions $(f_{i})_{i=0}^{N-1}$ that represent the finite sum term in~\eqref{eq:CT_stoch_tv_recon}.

\begin{center}
\begin{tcolorbox}[
enhanced,
attach boxed title to top center={yshift=-2mm},
colback=darkspringgreen!20,
colframe=darkspringgreen,
colbacktitle=darkspringgreen,
title=Family of $f_{i}$,
text width = 12cm,
fonttitle=\bfseries\color{white},
boxed title style={size=small,colframe=darkspringgreen,sharp corners},
sharp corners,
]
\begin{minted}{python}
Ai = ProjectionOperator(image_geometry, list_bi.geometry)
list_fi = []
for i in range(N):
    list_fi.append(LeastSquares(A=Ai[i], b=list_bi[i], c=0.5)
\end{minted}
\end{tcolorbox}
\end{center}

This leaves \texttt{Ai} as a CIL BlockOperator, \texttt{list\_fi} as a list of functions, and \texttt{list\_bi} as a BlockDataContainer.

\subsubsection{Experimental Set-up and Parameter Choices}

There is no ground truth image for this real-world dataset, so instead we construct an optimal solution for the deterministic optimisation problem~\eqref{eq:CT_deter_tv_recon}.

\begin{figure}
\centering
\includegraphics[width=0.9\linewidth]{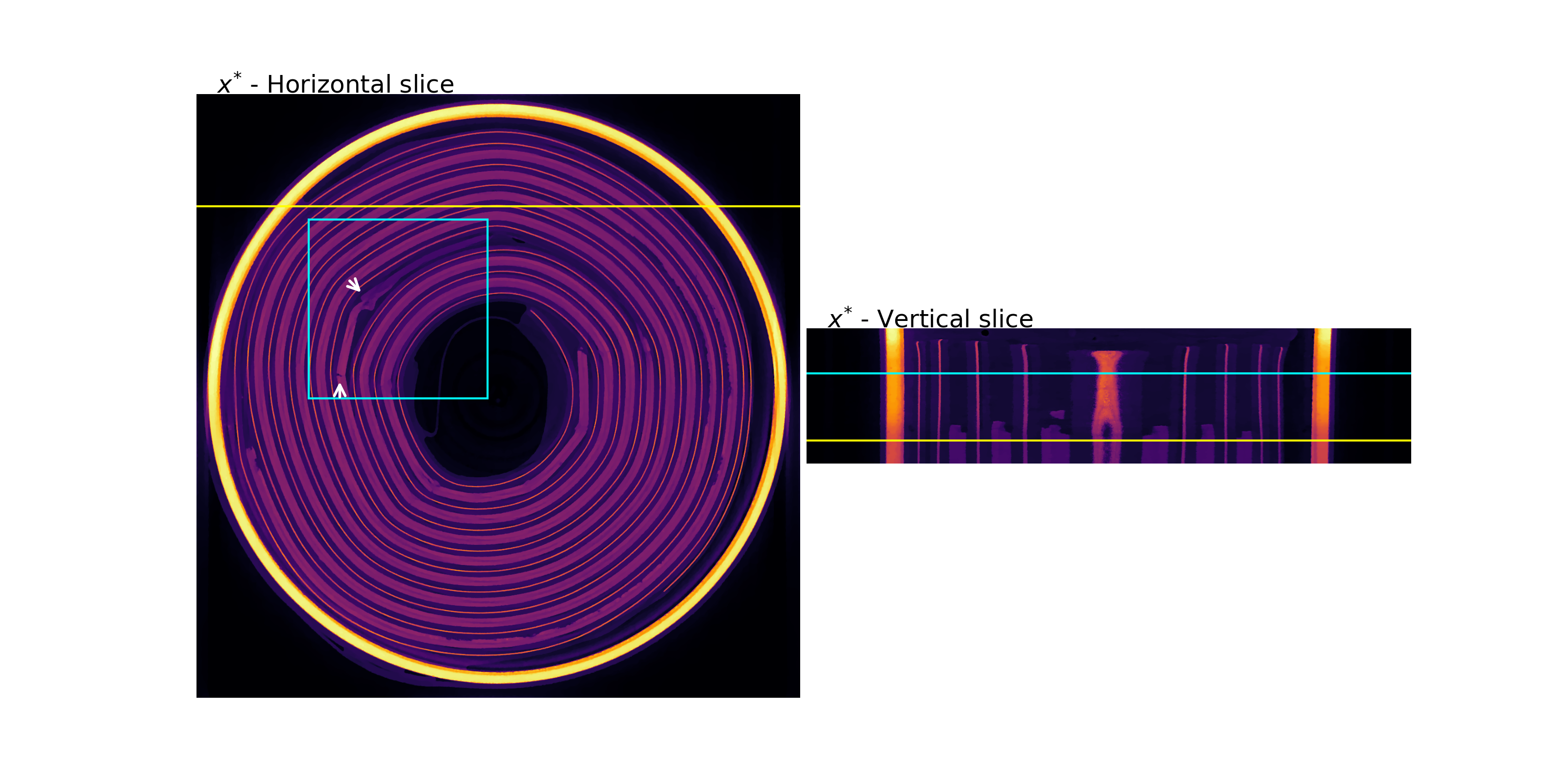}
\caption{The optimal 3D reconstruction of the CT battery dataset, used for comparisons. The blue box on the horizontal slice and blue line on the vertical slice are used for comparisons in figures~\ref{fig:50_subsets_horizontal_epochs} and~\ref{fig:50_subsets_vertical_epochs}. The yellow lines highlight the links between the two images - i.e. the yellow line on the left image is the vertical slice on the right image.   }
\label{fig:CT_optimal_fig}
\end{figure}

For the optimal solution, we use the APGD algorithm with Nesterov Momentum, algorithm~\ref{alg:FISTA}, where the TV proximal is calculated using the Fast Gradient Projection (FGP) algorithm~\cite{Beck2009} applied on the dual TV denoising problem~\cite{Chambolle2004} implemented with GPU acceleration in the CCPi regularisation toolkit~\cite{Kazantsev}. We do 300 inner iterations of the FGP algorithm to reduce any bias caused by this implicit calculation of the proximal.  We run the (outer) APGD algorithm for 5000 iterations. Each iteration corresponds to a full pass over the data (a ``data pass"), resulting in 5000 data passes in total to ensure convergence. We choose the step size to be $\frac{1}{L}$, see Figure~\ref{fig:CT_optimal_fig} for horizontal and vertical slices of the optimal reconstruction.
We compare this reference optimisation algorithm with a range of different deterministic and stochastic optimisation algorithms from the proposed framework to minimise~\eqref{eq:CT_deter_tv_recon} or~\eqref{eq:CT_stoch_tv_recon}, respectively. The parameter choices presented below are selected to yield stable, representative performance and to illustrate the capabilities of our new software design; a fully optimised comparison would require fine-tuning partition sizes, sampling strategies, step sizes, update frequencies and probabilities, which is beyond the scope of this paper.

For the stochastic comparisons, following the steps above, we consider:
\begin{enumerate}
\item \textbf{Partitioning}: We compare two types of partitioning \texttt{staggered} and \texttt{sequential}. We also compare a range of numbers of partitions: 10, 30 and 50 partitions.
\item \textbf{Sampler}: For all stochastic estimators, we use a   \code{random_with_replacement} sampler to provide fair comparison and to match convergence proofs in the literature.
\item \textbf{Stochastic gradient estimator}: We choose the following  stochastic gradient estimators, initialisations and update choices:
\begin{itemize}
\item ProxSAGA: we additionally ``warm start" the SAGA approximation by calculating the gradient of each function $f_i$ at the initial point.
\item ProxSVRG:  we use a full gradient update and snapshot every $2N$.
\item ProxLSVRG:  at each iteration, the probability of a full gradient update is $1/N$.
\end{itemize}
\item \textbf{Deterministic algorithm}: We choose the PGD algorithm, algorithm~\ref{alg:PGD}, and PD3O algorithm, algorithm~\ref{alg:PD3O},  for our deterministic algorithm, as these both allow for non-differential terms, such as the TV norm or the non-negativity constraint. For the PGD algorithm, we must choose a step size and use $\frac{1}{3\tilde{L}}$, $\frac{1}{\tilde{L}}$ and $\frac{1}{\tilde{L}}$ for the SAGA, SVRG and LSVRG stochastic gradient estimators, respectively,  where $\tilde{L}$ is the Lipschitz constant of $f_0$ multiplied by the number of functions. For PD3O, the step sizes are chosen as per the defaults in the CIL PD3O algorithm. Specifically,  with nomenclature as in~\cite{yan2018new}, the primal step size $\gamma = \frac{0.99*2.0}{L}$, where $L$ is the Lipschitz constant of $f$ and dual step size $\delta = \frac{L}{\|\Kop\|^2}$, where $\Kop$ is the operator which is composed with the function $h$ in PD3O, in our case the gradient operator multiplied by the regularisation parameter.
\end{enumerate}

Putting this all together, the following code blocks show a PGD-SAGA and PD3O-SGD setup. In the first code block, we first choose the  \code{random_with_replacement} sampler, as we do for each stochastic approximation. The second and third lines define a stochastic approximate gradient function with SAGA approximation and ``warm start" the SAGA approximation by calculating the gradient of each function $f_i$ at the initial point.  Each algorithm, deterministic and stochastic, is initialised with a volume of zeros. The code block then sets up the TV function $G$, with non-negativity and 100 inner iterations. The step size for SAGA is $\frac{1}{3\tilde{L}}$. Finally, in the last two lines, we set up PGD with the stochastic approximate gradient function and run for a number of iterations, using a callback to calculate the mean squared error at each data pass.
\begin{center}
\begin{tcolorbox}[
enhanced,
attach boxed title to top center={yshift=-2mm},
colback=darkspringgreen!20,
colframe=darkspringgreen,
colbacktitle=darkspringgreen,
title=ProxSAGA set-up,
text width = 0.8\linewidth,
fonttitle=\bfseries\color{white},
boxed title style={size=small,colframe=darkspringgreen,sharp corners},
sharp corners,
]
\begin{minted}{python}
sampler = Sampler.random_with_replacement(len(list_fi), seed=40)
F_saga = SAGAFunction(list_fi, sampler=sampler)
F_saga.warm_start_approximate_gradients(initial)
G = alpha * FGP_TV(max_iteration=100, device='gpu', nonnegativity=True)
step_size  = 1/(3*F_saga.L)
prox_saga = PGD(initial=initial, f=F_saga , g=G, step_size=step_size,
    update_objective_interval=N)
prox_saga.run(iterations=K, callbacks=[metric_callback])
\end{minted}
\end{tcolorbox}
\end{center}

There is a similar approach for the PD3O-SGD setup, although instead of using the TV function from the regularisation toolkit, we explicitly define the gradient operator and \code{MixedL21Norm} and pass this to the algorithm.
\begin{center}
\begin{tcolorbox}[
enhanced,
attach boxed title to top center={yshift=-2mm},
colback=darkspringgreen!20,
colframe=darkspringgreen,
colbacktitle=darkspringgreen,
title=PD3O-SGD set-up,
text width = 0.8\linewidth,
fonttitle=\bfseries\color{white},
boxed title style={size=small,colframe=darkspringgreen,sharp corners},
sharp corners,
]
\begin{minted}{python}
sampler = Sampler.random_with_replacement(len(list_fi), seed=40)
F_sgd = SGFunction(list_fi, sampler=sampler)
G = IndicatorBox(lower=0)
operator = alpha*GradientOperator(image_geometry)
H = L1Norm()
pd3o_sgd = PD3O(initial=initial, f=F_sgd, g=G, h=H, operator=operator,
    update_objective_interval=N)
pd3o_sgd.run(iterations=K, callbacks=[metric_callback])
\end{minted}
\end{tcolorbox}
\end{center}

We also compare the stochastic algorithms against some deterministic algorithms:
\begin{itemize}
\item APGD - although the optimal solution was computed with APGD, for a more realistic APGD run, we compute just 50 iterations, which is equivalent to 50 data passes, with the TV proximal calculated with just 100 inner iterations. We take a step size of $1/{L}$, where $L$ is the Lipschitz constant of $f$.
\item PD3O: we use PD3O, set up as per the stochastic version but with a deterministic \code{f}. Again, we use the default step sizes in CIL.

\end{itemize}

To compare the results, we consider two metrics, the first is the optimisation objective, \eqref{eq:CT_deter_tv_recon} or~\eqref{eq:CT_stoch_tv_recon}, for deterministic and stochastic,  respectively. The second is a normalised root mean squared error between the reconstructed volume, $\bx$, and the optimal volume, $\bx^*$,
\begin{equation}
NRMSE(x) = \dfrac{\|\bx - \bx^{*}\|_2}{\|\bx^{*}\|_2}.
\end{equation}

These metrics are calculated and saved each iteration for the deterministic methods, and every data pass for the stochastic methods.

\subsection{Results}

\subsubsection{Number of Subsets}
\begin{figure}
\begin{subfigure}[t]{0.475\textwidth}
\centering
\includegraphics[width=\linewidth]{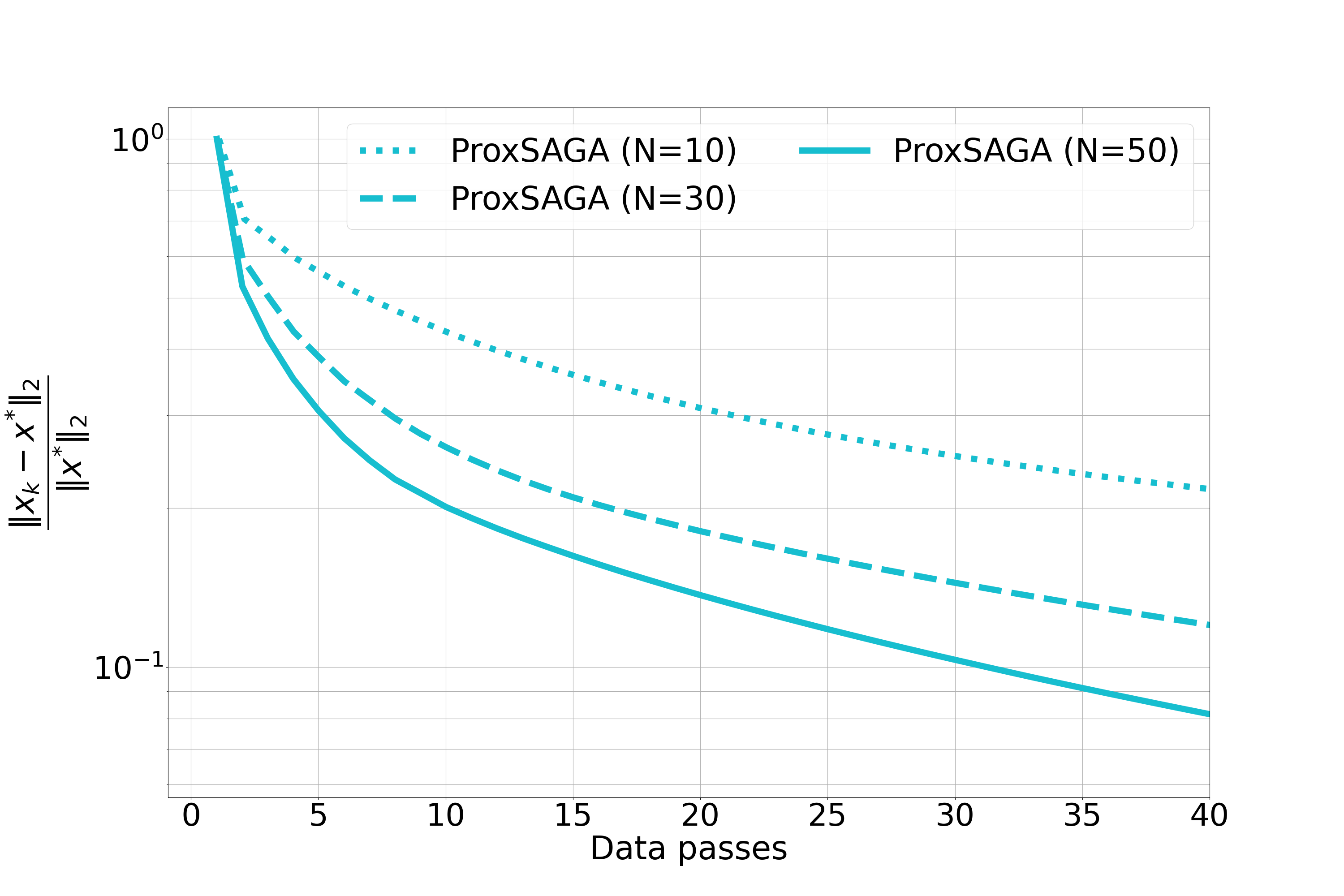}
\caption{ProxSAGA with 10, 30 and 50 subsets }
\label{fig:ct-saga-subsets}
\end{subfigure}
\begin{subfigure}[t]{0.47\textwidth}
\centering
\includegraphics[width=\linewidth]{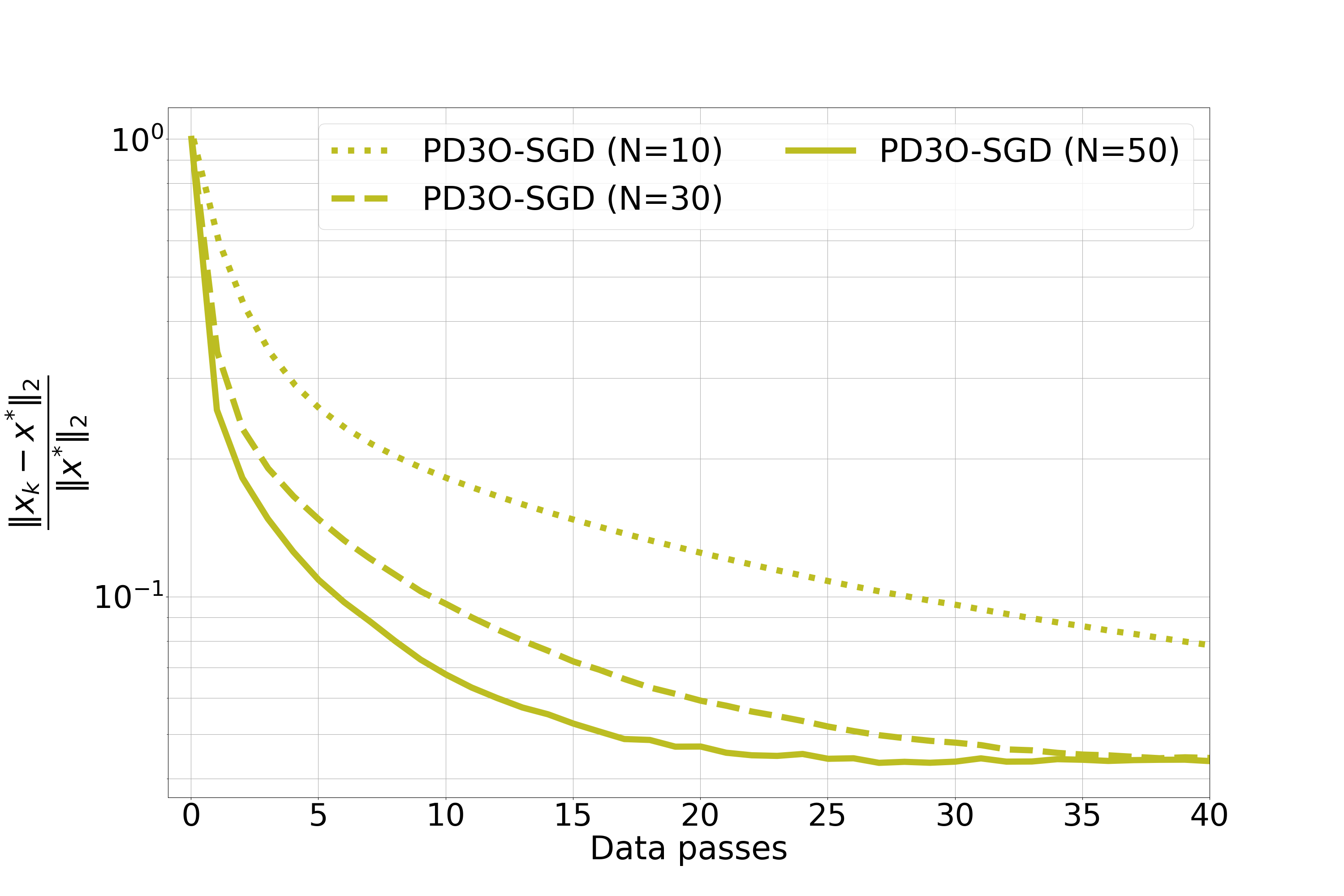}
\caption{PD3O-SGD with 10, 30 and 50 subsets}
\label{fig:ct-pd3o-subsets}
\end{subfigure}
\caption{Comparing the normalised mean squared error of SAGA and PD3O-SGD algorithms against data passes.}
\end{figure}

We begin by comparing various choices for the number of subsets used in the stochastic algorithms. Figures~\ref{fig:ct-saga-subsets} and~\ref{fig:ct-pd3o-subsets} present the normalised mean squared error against data passes for both the SAGA and PD3O-SGD stochastic methods. The results demonstrate that increasing the number of subsets enhances the convergence speed for these methods, measured in data passes. This trend is also observed in other methods, although these results are omitted here for brevity. In the remainder of this section, we analyse the case with 50 subsets.

\subsubsection{Primal-Dual Algorithms}
\begin{figure}
\centering
\begin{subfigure}[t]{0.47\textwidth}
\centering
\includegraphics[width=\linewidth]{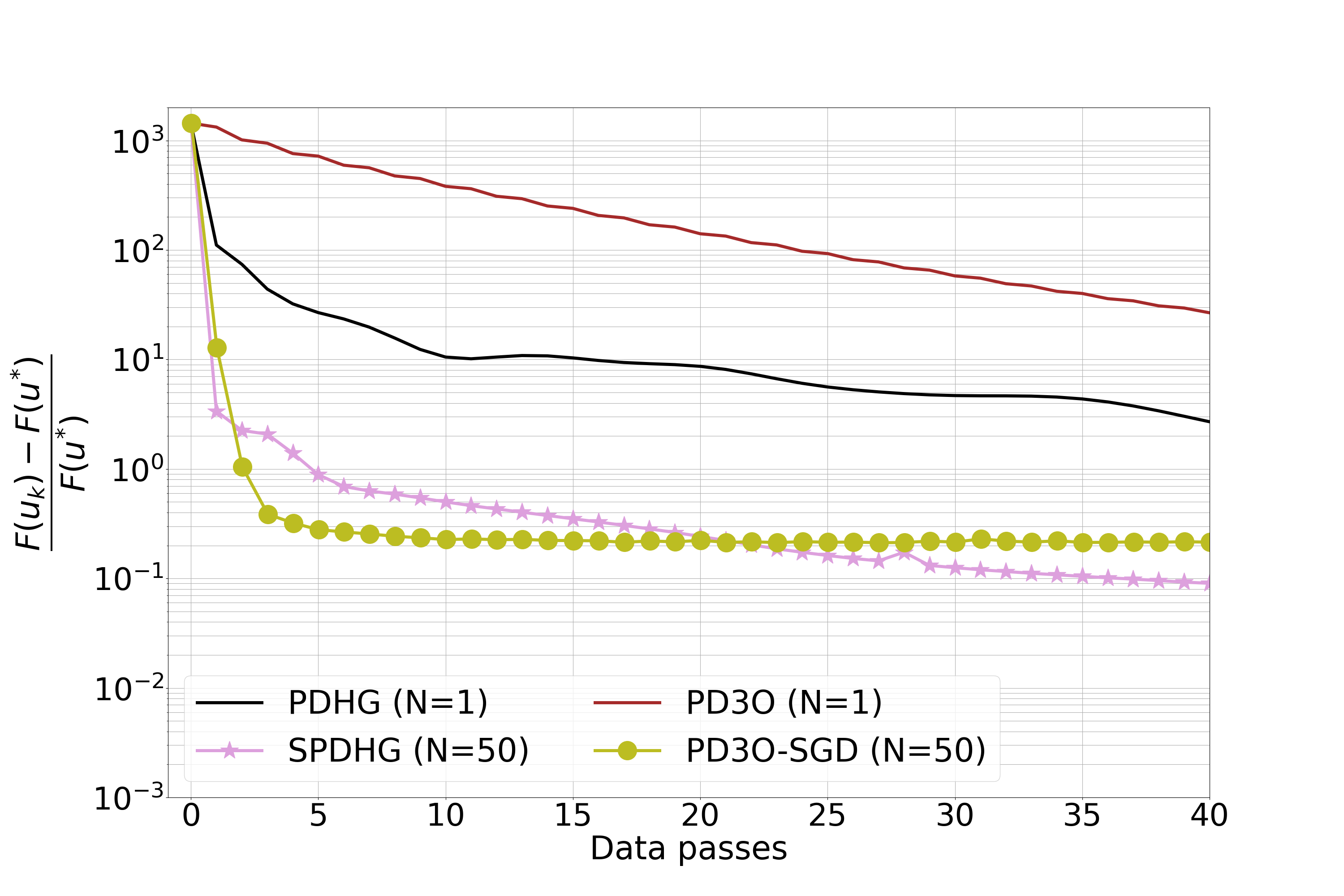}
\caption{Primal-dual algorithms comparing objective value}
\label{fig:primal_dual_obj}
\end{subfigure}
\begin{subfigure}[t]{0.47\textwidth}
\centering
\includegraphics[width=\linewidth]{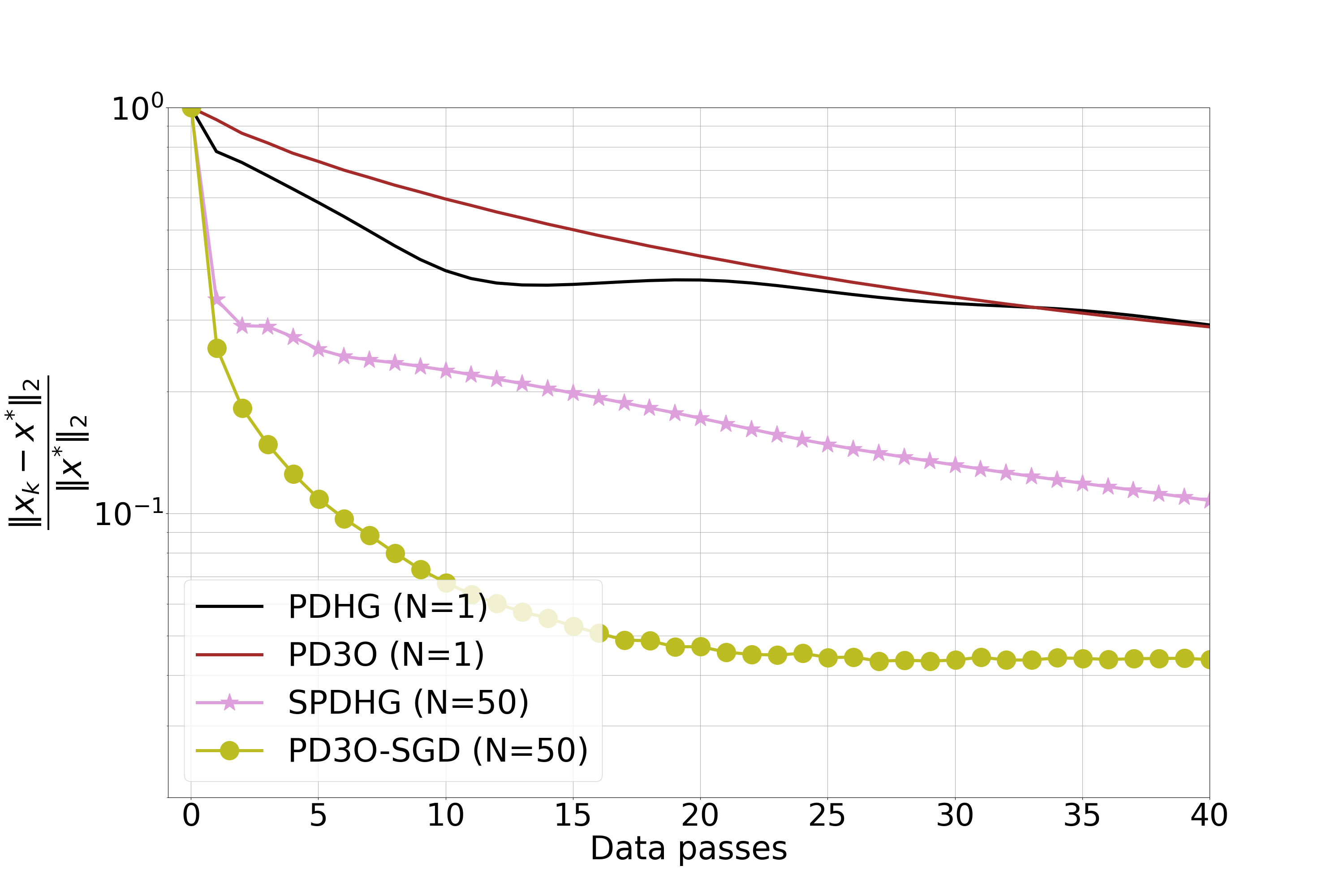}
\caption{Primal-dual algorithms comparing normalised mean squared error}
\label{fig:primal_dual_nrmse}
\end{subfigure}
\begin{subfigure}[t]{0.47\textwidth}
\centering
\includegraphics[width=\linewidth]{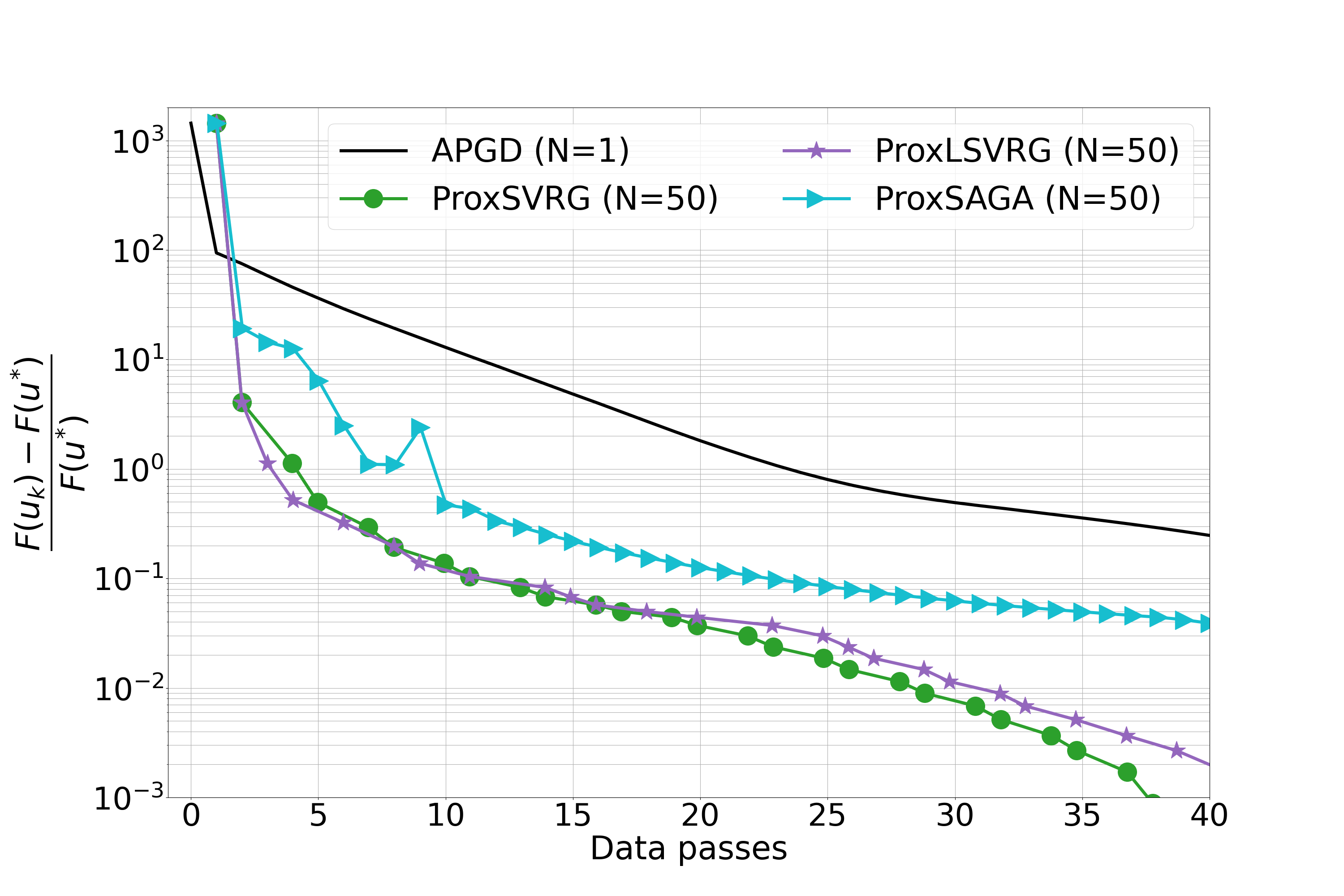}
\caption{Proximal gradient descent algorithms comparing objective value}
\label{fig:gradient_obj}
\end{subfigure}
\begin{subfigure}[t]{0.47\textwidth}
\centering
\includegraphics[width=\linewidth]{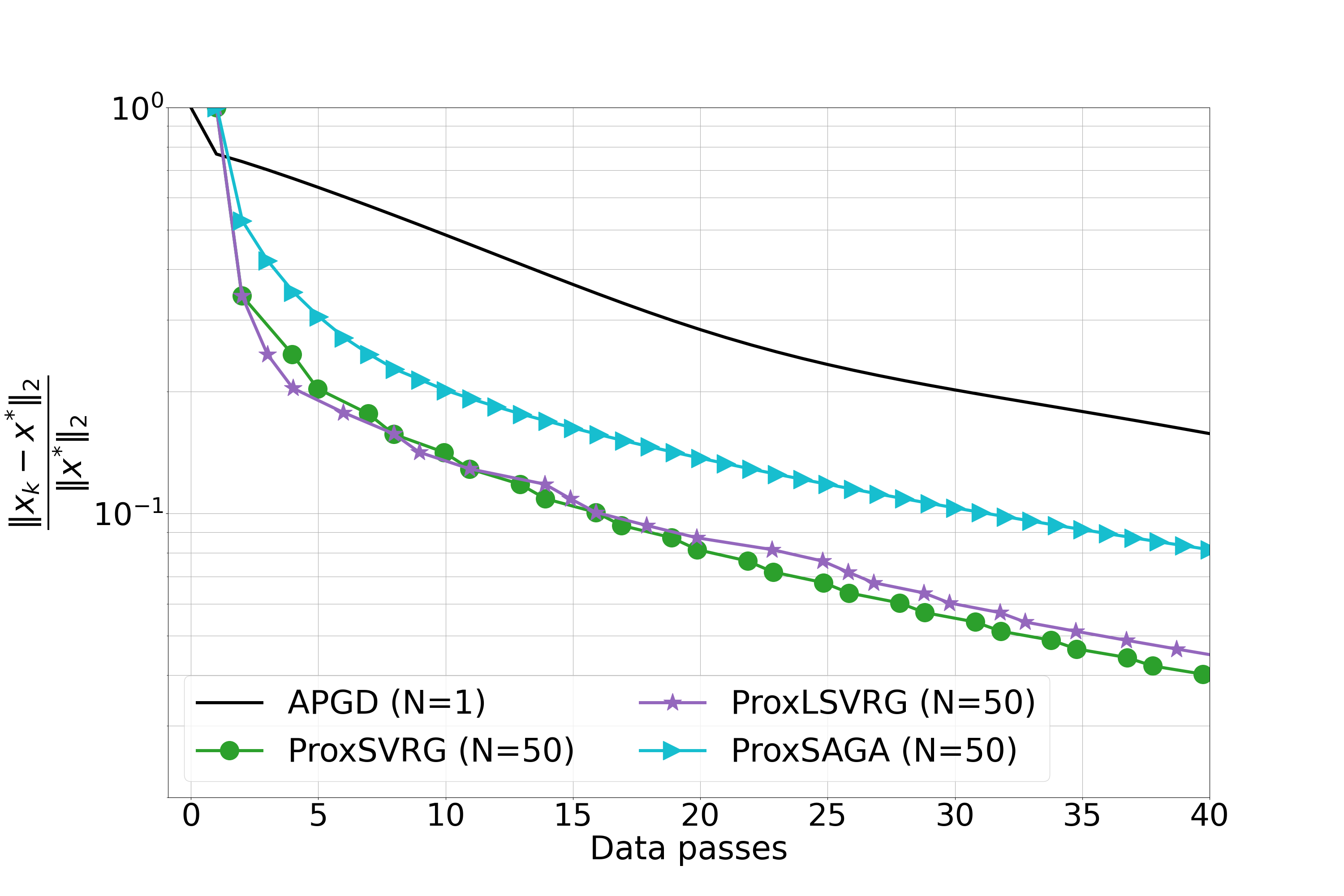}
\caption{Proximal gradient descent algorithms comparing normalised mean squared error}
\label{fig:gradient_nrmse}
\end{subfigure}
\caption{Comparing the normalised mean squared error and objective values of the proximal gradient descent and primal-dual algorithms against data passes. For the stochastic algorithms, we consider just 50 subsets. }
\label{fig:CT-epoch-50-subsets}
\end{figure}

In Figure~\ref{fig:CT-epoch-50-subsets} and in particular sub-figures~\ref{fig:primal_dual_obj} and~\ref{fig:primal_dual_nrmse}. We see that where the deterministic algorithms PDHG and PD3O show similar convergence rates on the mean squared error, PDHG seems to converge in fewer data passes for the objective value. The PDHG and PD3O differences illustrate that convergence in objective value may not tell the whole story, and, looking at convergence in other metrics, where they are available, is valuable. We see the stochastic algorithms SPDHG and PD3O-SGD converge faster than both deterministic algorithms in both metrics. In particular, PD3O-SGD does particularly well in the NRMSE metric.

In Figure~\ref{fig:CT-epoch-50-subsets}, particularly in sub-figures~\ref{fig:primal_dual_obj} and~\ref{fig:primal_dual_nrmse}, we observe that the deterministic algorithms PDHG and PD3O exhibit similar convergence rates for the mean squared error. However, PDHG appears to converge in fewer data passes for the objective value, emphasising the fact that these metrics measure different characteristics. The stochastic algorithms SPDHG and PD3O-SGD demonstrate faster convergence compared to both deterministic algorithms across both metrics. Notably, PD3O-SGD performs exceptionally well in the NRMSE metric.

\subsubsection{Proximal Gradient Descent Algorithms}

Once again, in Figure~\ref{fig:CT-epoch-50-subsets}, particularly in sub-figures~\ref{fig:gradient_obj} and~\ref{fig:gradient_nrmse}, we compare the deterministic APGD algorithm with ProxSVRG, ProxSAGA, and ProxLSVRG, all using 50 subsets. By 20 data passes, both ProxSVRG and ProxLSVRG achieve results within 10\% of the optimal solution, whereas the APGD algorithm does not reach this level within the 40 data passes shown.

\begin{figure}
\begin{subfigure}[t]{\textwidth}
\centering
\includegraphics[width=\linewidth]{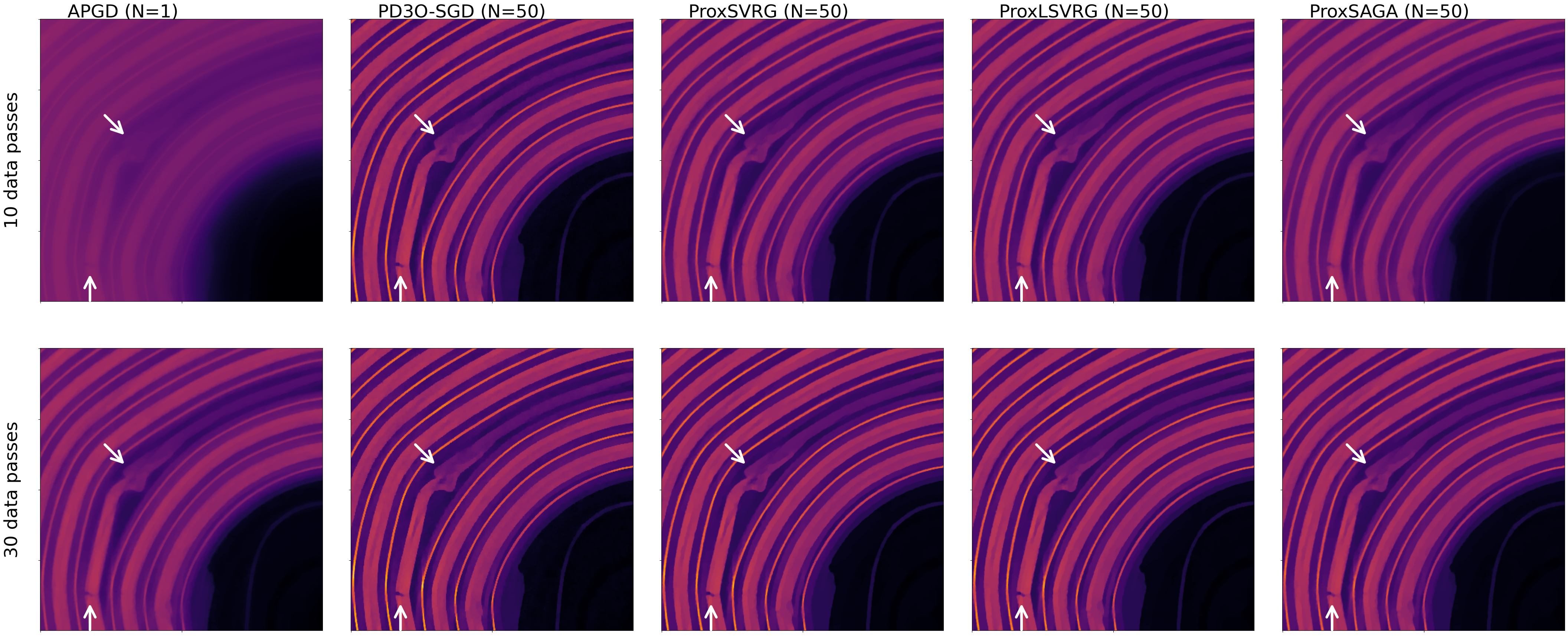}
\caption{Area of interest in the horizontal slice, zoom area given in Figure~\ref{fig:CT_optimal_fig}.}
\label{fig:50_subsets_horizontal_epochs}
\end{subfigure}
\begin{subfigure}[t]{\textwidth}
\centering
\includegraphics[width=\linewidth]{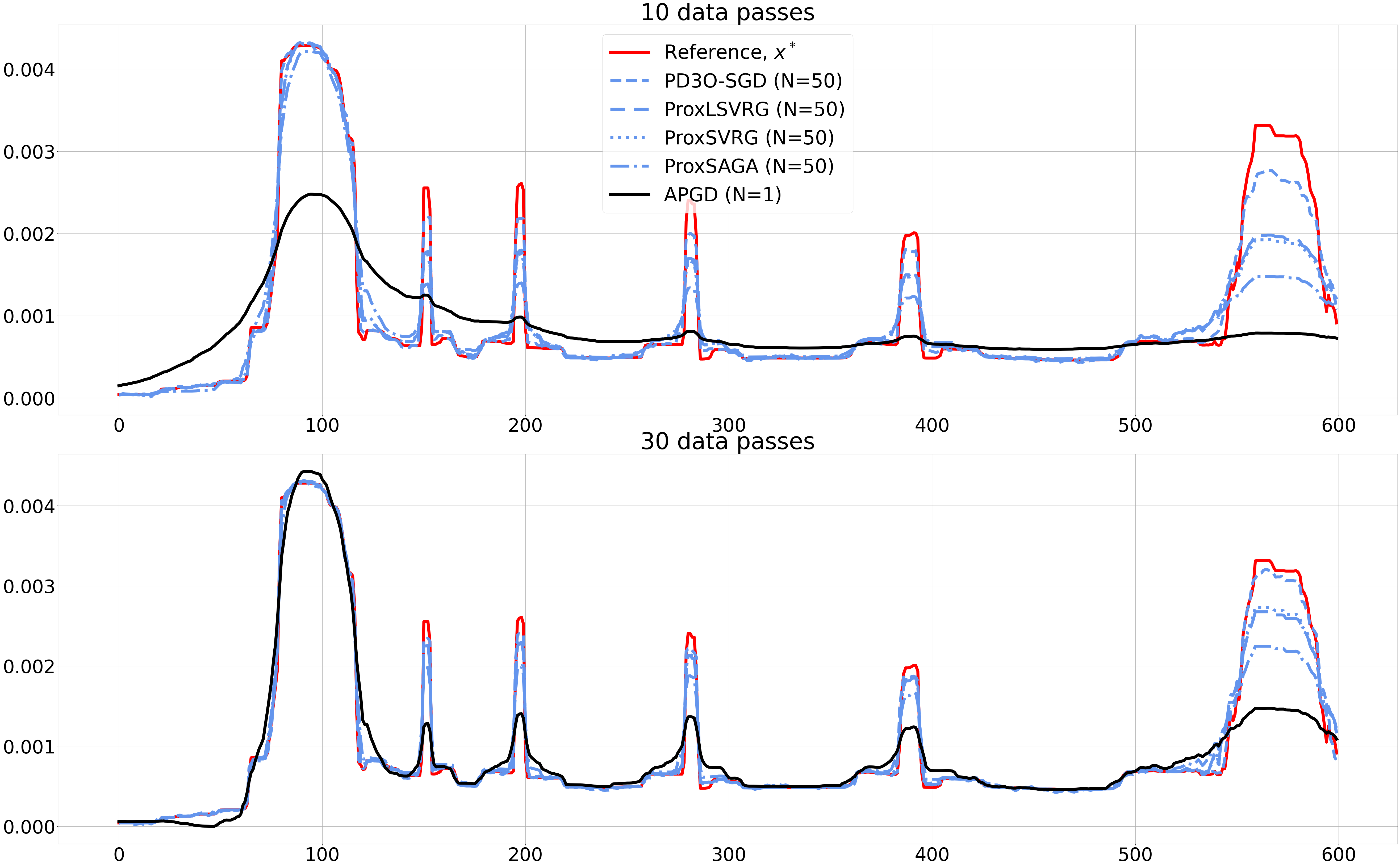}
\caption{Line plot of a cut through the vertical slice shown in Figure~\ref{fig:CT_optimal_fig}}
\label{fig:50_subsets_vertical_epochs}
\end{subfigure}
\caption{Plots comparing the reconstructed volume for the different algorithms at different numbers of data passes.}
\end{figure}

To visualise this, we compare ProxSVRG, ProxLSVRG, ProxSAGA and PD3O-SGD with the deterministic APGD algorithms at 10, 20 and 30 data passes. In Figure~\ref{fig:50_subsets_horizontal_epochs}, we see 10 and 30 data passes for a region of interest in the horizontal slice.  We see at 10 data passes that the internal structure has not been reconstructed for the APGD algorithm, and, for example, the truncated thick layer in the upper left of the region of interest is not visible. Whereas in all the stochastic images the structures are visible, and most clear in the PD3O-SGD reconstruction. By 30 data passes, the internal structure is visible in the APGD reconstruction, but even there, it is more blurred than the stochastic reconstructions. In Figure~\ref{fig:50_subsets_vertical_epochs}, we look at a line plot of a cross-section of the vertical slice (see the line in Figure~\ref{fig:CT_optimal_fig}). We see clearly that the APGD reconstruction has not yet managed to reconstruct the central peak at all, compared to the PD3O-SGD algorithm, which has managed to reconstruct all of the peaks of the optimal solution.

\subsubsection{Conclusion}
This case study demonstrates the stochastic framework on a CT use case, including the CIL CT data partitioner. We see that the deterministic algorithms are outperformed by their stochastic counterparts. We see that the choice of the number of subsets has a significant effect, and that a larger number of subsets appears to have better convergence values, in terms of data passes. We see the flexibility of the stochastic framework, testing a range of stochastic algorithms with only small changes in the code.

%% file: version3/05_example_pet.tex
\section{Case study: Positron Emission Tomography}\label{sec:pet}

In this section, we focus on a Positron Emission Tomography (PET) case, thanks to the integration of CIL with the Synergistic Image Reconstruction Framework (SIRF)~\cite{Ovtchinnikov2020}. 
The integration with SIRF has been described in~\cite{Jorgensen2021,Brown2021}, and it is based on Python's duck-typing, which enables to use CIL algorithms on objects that behave like its \code{Function}, \code{Operator} or \code{DataContainer}, thereby extending the original capabilities of CIL and demonstrating its adaptability when used in conjunction with third-party libraries. 

Furthermore, we demonstrate that not only is it easy to set up a stochastic algorithm to reconstruct a PET dataset, but also that appropriate tuning of 
the algorithm can have a large impact on its convergence speed. We will limit our attention to two SAGA algorithms with different preconditioning and step size rules.

\subsection{Dataset}

\begin{figure}[tb]
    \centering
    \begin{subfigure}[t]{0.3\textwidth}
        \centering
        \includegraphics[width=\textwidth]{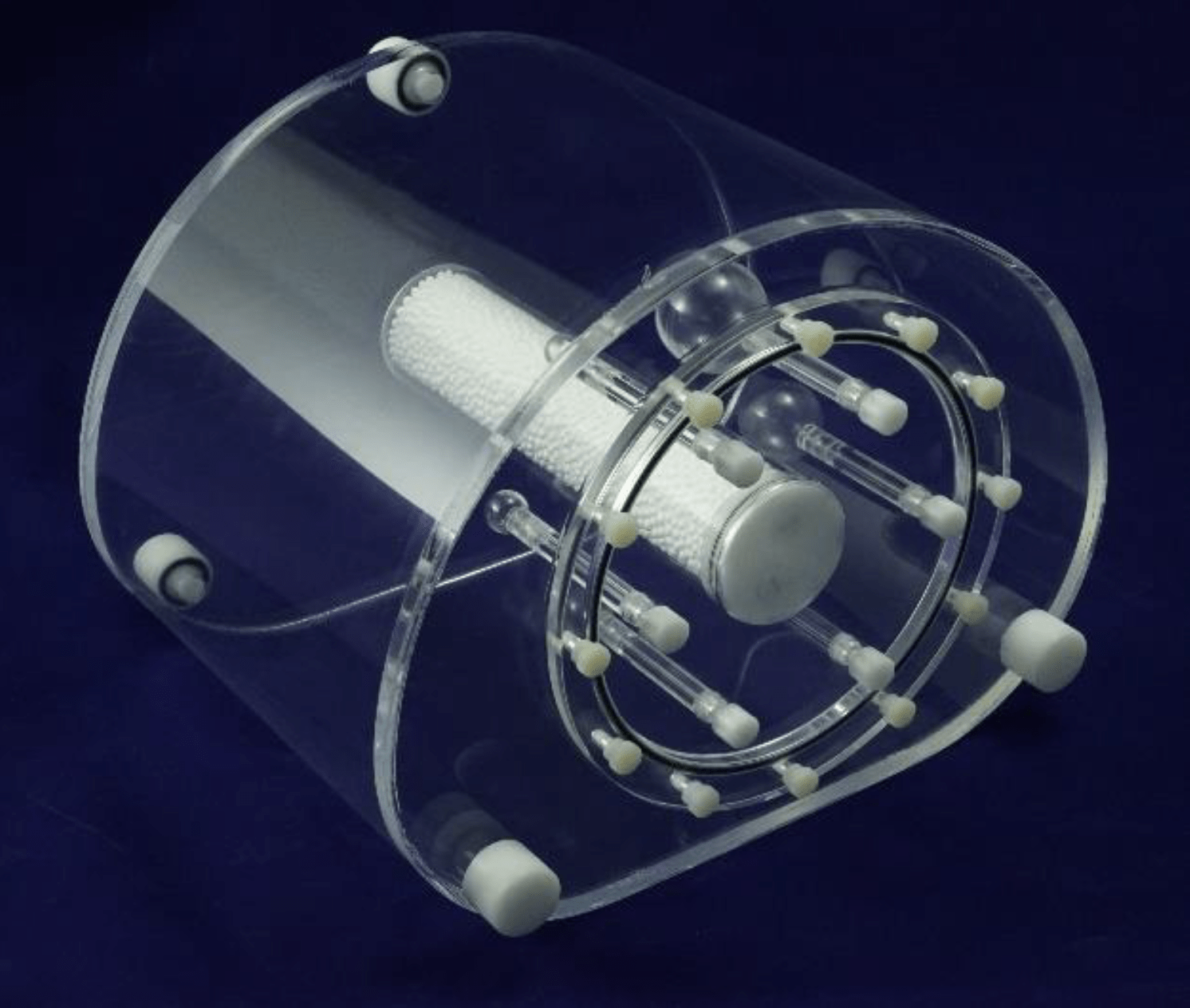}
        \label{fig:sub1}
    \end{subfigure}
    \begin{subfigure}[t]{0.3\textwidth}
            \centering
            \includegraphics[width=\textwidth]{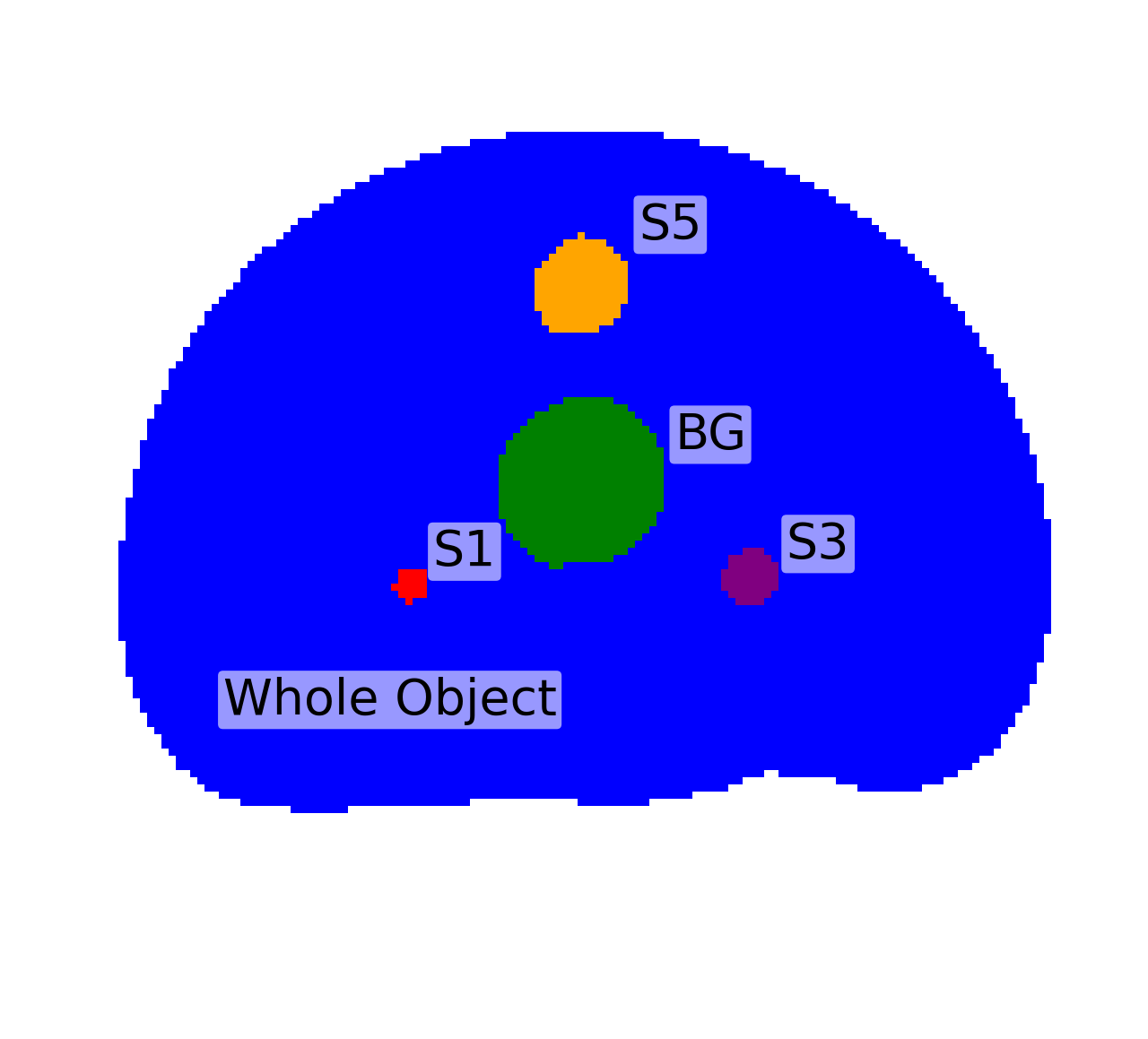}
            \label{fig:sub2}
        \end{subfigure}
    \begin{subfigure}[t]{0.3\textwidth}
        \centering
        \includegraphics[width=\textwidth]{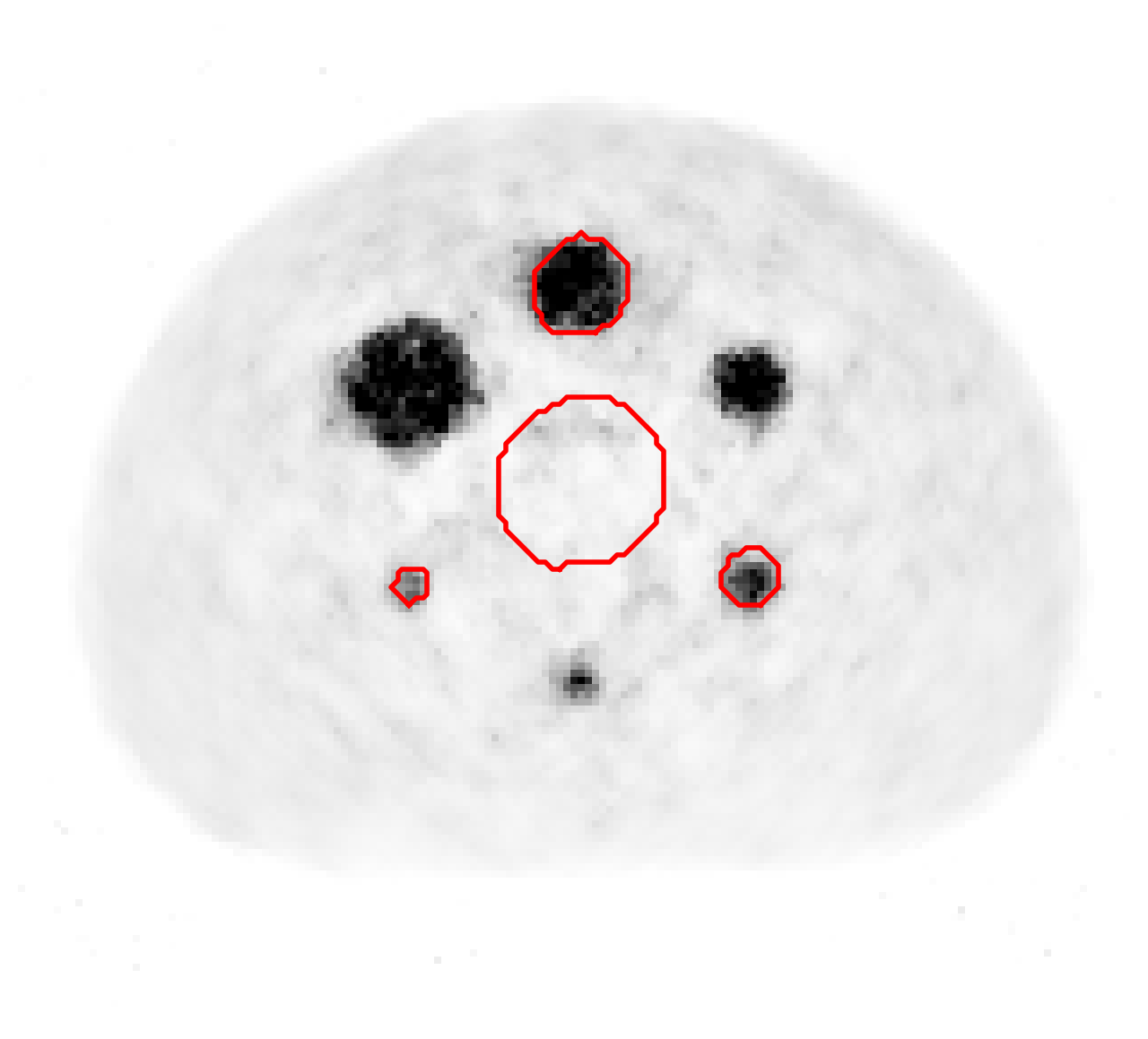}
        \label{fig:sub3}
    \end{subfigure}    
        \caption{Left: a picture of the NEMA phantom used in the PET comparison (Image reproduced with permission from Data Spectrum). Centre: a slice describing the volumes of interest (VOI) used for the evaluation in the PET Rapid Image Reconstruction Challenge~\cite{PETRIC_arxiv}: S1, S3, S5, BG and Whole Object indicate respectively: sphere with 10mm internal diameter, sphere with 17mm internal, sphere with 28 mm internal diameter, background, the whole object of interest. Out of these, we used the whole object and S1. Right: a slice of the reference reconstruction provided by the PETRIC organisers with highlighted S1, S3, S4 and BG.}        
    \label{fig:nema_pet}
\end{figure}

In order to demonstrate our new stochastic optimisation framework for PET, we use the National Electrical Manufacturers Association Image Quality (NEMA IQ) Phantom~\cite{NEMAstandard}, which is a standard phantom for testing PET scanner and reconstruction performance, see Figure~\ref{fig:nema_pet}. It consists of a Perspex container with inserts of different-sized spheres, some filled with liquid with a higher radioactivity concentration than the background, others with cold (non-radioactive) water. This allows assessment of resolution and quantification. 
We use a 12-min PET dataset~\cite{PETRIC_NEMAIQlowcountsDataset} of the NEMA IQ phantom that was acquired on a Siemens Biograph mMR PET/MR scanner at the Institute of Nuclear Medicine, UCLH, London. 
The data used in this section, the optimisation problem and the data partitioner, described below, were developed for the PET Rapid Image Reconstruction Challenge (PETRIC),~\cite{PETRIC_arxiv}, which was held in the summer of 2024 to foster research on fast image reconstruction algorithms for PET reconstruction. 

\subsection{Methods}

\subsubsection{Objective Function\label{sec:pet-objective}}

In this PET case study, we consider a maximum a-posteriori (MAP) estimation problem incorporating a smoothed Relative Difference Penalty (RDP) ~\cite{Nuyts2002} and a non-negativity constraint. The objective is to recover a non-negative image $\bx$
by solving the following optimisation problem ~\cite{Qi2006}:
 \begin{equation}
 \argmin_{\bx} \sum \bA\bx - \bb \mathrm{log}(\bA\bx + \bm{\eta}) + \alpha\mathrm{RDP}(\bx) + \mathbb{I}_{\{\bx>0\}}(\bx)
 \label{PET_RDP}
 \end{equation}
 where $\bA$ is the operator modelling the system and $\bb$ is the acquired data, and $\bm{\eta}$ is the background data, e.g., random and scatter events. 

The first term in~\eqref{PET_RDP} is the Kullback-Leibler data-fidelity term, which arises from the negative log-likelihood of the Poisson distribution used in the MAP formulation. The RDP is  defined as 
 \begin{equation}
     \mathrm{RDP}(\bx) := \frac{1}{2}\sum_{i=1}^{N_{v}}\sum_{\ell\in \mathcal{N}_{i}}w_{ij}\kappa_i\kappa_j\frac{(\bx_{i}-\bx_{\ell})^{2}}{(\bx_{i}+\bx_{\ell}) + \omega\|\bx_{i}-\bx_{\ell}\| + \epsilon}
 \label{RDP}
 \end{equation}
 where $\mathcal{N}_{i}$ is a $3\times3\times3$ neighbourhood around the ith image voxel, $\omega=2$ is a parameter used to control edge preservation, $w_{ij}$
is a weight factor (here taken as ``horizontal" voxel-size divided by Euclidean distance between the $i$-th and $j$-th voxels) and $\epsilon$ is a small constant to make the penalty differentiable. The parameter $\bm{\kappa}$ is an image to give voxel-dependent weights. In our setup, $\bm{\kappa}$ is precomputed as the square root of the negative row-sum of the Hessian of the log-likelihood evaluated at an initial reconstruction with the standard Ordered Subset Expectation Maximisation (OSEM) algorithm (see Equation 25 in~\cite{Tsai2020}).

\subsubsection{Data Partitioning}\label{subsec:PET-partitioning} 

For the stochastic optimisation algorithms we express~\eqref{PET_RDP}, similar to the CT case study, as a finite-sum optimisation problem. We begin by partitioning the PET data, following a similar procedure outlined in the previous sections. Given the number of subsets $N$, we partition our acquired data $\bb$ and the noisy background data $\bm{\eta}$ to $\bb_{i}$ and $\bm{\eta_{i}}$ respectively for $i=0,1,\cdots,N-1$. For each $i$, we define the corresponding operators $\bA_{i}$. Using the above partitioned data $\bb_{i}$, $\bm{\eta_{i}}$ and operators $\bA_{i}$, we define the following stochastic optimisation objective:
\begin{align}
&\argmin_{\bx}  \sum_{i=0}^{N-1} f_i(\bx)+g(\bx) 
\\ f_i(\bx) &:= \sum_i\bA_{i}(\bx) - \bb_{i}\mathrm{log}(\bA_{i}\bx + \bm{\eta}_{i}) + \frac{\beta}{N} \mathrm{RDP}(\bx)\\
g(\bx)&:=\mathbb{I}_{\{\bx>0\}}(\bx)
\end{align}

The PET data \code{partitioner}, which is available in the SIRF-contribs repository~\cite{PET_partitioner}, follows a similar interface to the CIL data partitioner described above. It takes as arguments the experimental PET data, background corrections  (i.e. scatter and randoms), and the bin efficiencies (including attenuation and detection efficiencies), as well as a mode for partitioning, which could be ``sequential", ``staggered" or ``random permutation". 

In contrast to CIL, the data \code{partitioner} not only partitions the data but also returns a list of SIRF objective functions for each subset, which we again call \code{list_fi}. This is because the construction of the $A_i$ operators in PET depends on the subset via the bin efficiencies and background.  These SIRF objective functions behave like CIL \code{Function} and can be passed to CIL algorithms. Note that the returned SIRF objective function objects are meant for maximisation, so they need to be multiplied by -1 to be used with CIL optimisation algorithms, which perform minimisation.

 We use a staggered partition scheme, as in the CT case.  We partition the data into 21, 42 and 63 subsets, each a divisor of the total number of views in the dataset, which is 252.

\subsubsection{Experimental Setup and Parameters}

In this section, we show how different step sizes and preconditioning schemes affect the performance of the same SAGA algorithm, obtained by plugging the \code{SAGAFunction} stochastic estimator into the PGD algorithm.

\paragraph{Step Size Strategies}
\label{sec:step-size-strategies}
Along with the stochastic framework, we include new functionality to support non-constant step sizes for first-order optimisation problems. Specifically, we have implemented methods such as Armijo~\cite{armijo1966minimization} and Barzilai-Borwein~\cite{BARZILAI1988}. 
Moreover, users are given the flexibility to define their own step size strategies. Subclassing the base abstract class \code{StepSizeRule}, users need to override the \code{get_step_size} abstract method, which is called during the \code{update} step of the algorithm classes \code{GD}, \code{PGD}, \code{APGD}, providing access to all attributes of the running algorithm instance. This way, users can customise and implement various step size rules.

The first algorithm, referred to as SAGA-1, uses a decreasing step size rule defined by $\gamma_k = \gamma_0 / (1 + \beta k)$, with $\gamma_0$ an initial step size, $\beta$ a decay parameter and $k$ the current iteration number; here $\beta=0.02$ and $\gamma_0=1$. This can be constructed in CIL by defining a \code{DecreasingStepSize} class, which inherits from the \code{StepSizeRule} base abstract class as below.

\begin{center}
\begin{tcolorbox}[
    enhanced,
    attach boxed title to top center={yshift=-2mm},
    colback=darkspringgreen!20,
    colframe=darkspringgreen,
    colbacktitle=darkspringgreen,
    title=Decreasing Stepsize,
    text width = 9.5cm,
    fonttitle=\bfseries\color{white},
    boxed title style={size=small,colframe=darkspringgreen,sharp corners},
    sharp corners,
]
\begin{minted}{python}
class DecreasingStepSize(StepSizeRule):

  def __init__(self, gamma0, beta=0.02):
        self.gamma0 = gamma0
        self.beta = beta

    def get_step_size(self, algorithm):
        step_size = self.gamma0 /
             (1.0 + self.beta * algorithm.iteration)
        return step_size

step_size_rule = DecreasingStepSize(gamma0=initial_step_size)
\end{minted}
\end{tcolorbox}
\end{center}

\paragraph{Preconditioners}
We introduce a new preconditioning feature via the \code{Preconditioner} class, which can be used to manipulate the gradient at each iteration. This can be combined with the algorithms \code{GD}, \code{PGD} and \code{APGD}. In many cases, the preconditioning step multiplies the gradient by a preconditioning array $P_k$, see for example step 4 in Algorithm~\ref{alg:PGD}.

For the SAGA-1, 
we select a preconditioning array $P_{k}$ given by 
the inverse of the negative row-sum of the Hessian of the Poisson log-likelihood~\cite{Tsai2017}: 
\begin{equation}
    P_{k} = \left ( \bA_i^\top \left \{\frac{\bb_i}{(\bA_i \bx_{k} + \bm{\eta}_i)^2}\right \} \bA_i \bold{1}  \right )^{-1} = (\kappa^2 +\varepsilon)^{-1}.
\end{equation}

The code block below demonstrates how to define the above preconditioner by subclassing the \code{Preconditioner} class and overriding the \texttt{apply} method. This allows users to implement custom preconditioning strategies for their specific optimisation problems.
\begin{center}
    \begin{tcolorbox}[
        enhanced,
        attach boxed title to top center={yshift=-2mm},
        colback=darkspringgreen!20,
        colframe=darkspringgreen,
        colbacktitle=darkspringgreen,
        title=SAGA-1: $P_{k}$,
        text width = 9.5cm,
        fonttitle=\bfseries\color{white},
        boxed title style={size=small,colframe=darkspringgreen,sharp corners},
        sharp corners,
    ]
\begin{minted}{python}
class TsaiPreconditioner(Preconditioner):

    def __init__(self, kappa, epsilon=1e-6):
        self.kappasq = kappa*kappa + epsilon

    def apply(self, algorithm, gradient, out=None):
        return gradient.divide(self.kappasq, out=out)

preconditioner = TsaiPreconditioner(kappa)        
\end{minted}
\label{code:pet-preconditioner}
\end{tcolorbox}
\end{center}

    \begin{center}
\begin{tcolorbox}[
    enhanced,
    attach boxed title to top center={yshift=-2mm},
    colback=darkspringgreen!20,
    colframe=darkspringgreen,
    colbacktitle=darkspringgreen,
    title=SAGA-1 with Decreasing Stepsize and Preconditioning,
    text width = 0.7\linewidth,
    fonttitle=\bfseries\color{white},
    boxed title style={size=small,colframe=darkspringgreen,sharp corners},
    sharp corners,
]
    \begin{minted}{python}
sampler = Sampler.random_with_replacement(len(list_fi), seed=40)    
F_saga = -SAGAFunction(list_fi, sampler=sampler)
G = IndicatorBox(lower=0) 
saga_1 = ISTA(initial = OSEM_recon, f=F_saga, g=G, 
            update_objective_interval = N, 
            step_size = step_size_rule,
            preconditioner = preconditioner)
saga_1.run(iterations=K, callbacks=[metric_callback])
\end{minted}
\end{tcolorbox}
\end{center}

Note that in the above code snippet, the leading minus sign in the \code{SAGAFunction} is introduced to adjust to a minimisation problem, see \ref{subsec:PET-partitioning}.
    
The second algorithm, SAGA-SOS, is the SAGA entry to the PETRIC challenge~\cite{PETRIC_arxiv} submitted by the Stochastic Optimisation Submission (SOS) team. The algorithm features a decaying step size rule and  a modified BSREM-type preconditioner~\cite{Ahn2003}, defined as:  

\begin{equation}
    P_k = \frac{\bx^{k}}{\bA^\top \mathbf{1}}
\end{equation} 
where $\bA^\top \mathbf{1}$ represents the sensitivity image, obtained by applying a backprojection to uniform data. 
The step size is determined using the Armijo rule: starting from an initial user-defined step size, it is reduced iteratively by a factor of 2 until the condition \begin{equation}
    f(\bx_{k} - \gamma_{k} d_{k}) \leq f(\bx_{k}) + \sigma \gamma_{k} \nabla f(\bx_{k})^\top d_{k}
\end{equation} 
is satisfied. Here, $\sigma$ is a user-defined tolerance (set to $0.2$ for SAGA-SOS), and $d_{k} = P_{k}\nabla f(\bx_{k})$ denotes the preconditioned search direction~\cite{armijo1966minimization}. 
The SAGA-SOS algorithm begins by estimating a suitable step size through 5 iterations of PGD. Note that under the non-negativity constraint, this reduces to a Projected Gradient Descent with an Armijo line search.
The smallest step is defined as $\gamma_0 = \min\{\alpha_1, \alpha_2, \dots, \alpha_5\}$ where $\alpha_l$ is the step size computed at the lth iteration. This $\gamma_0$ is then used as the initial step size in the stochastic reconstruction phase. During this phase, the step size is updated iteratively according to $\gamma_k = \gamma_0 / (1 + \beta k)$, where $k$ is the iteration index and $\beta = 0.01$ is a user-defined decay parameter.

\subsection{Results}

\begin{figure}
    \centering
    \includegraphics[width=0.9\linewidth]{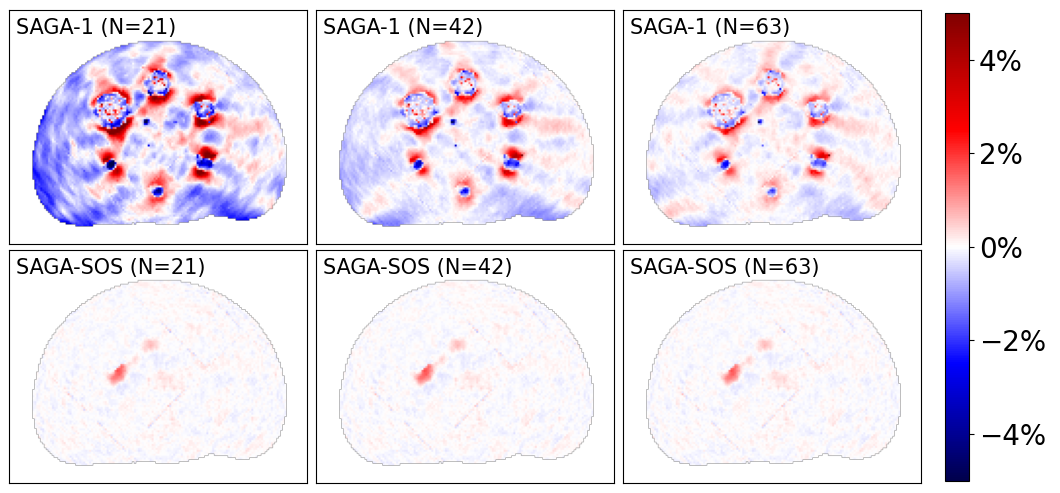}
    \caption{Relative difference of the reconstructions at 30 data passes of the SAGA-SOS and SAGA-1 algorithms with respect to the reference solution for the whole volume as defined in fig~\ref{fig:nema_pet}. The first row relates to the SAGA-1 algorithm, while the second relates to SAGA-SOS. The columns relate to the number of subsets: 21, 42, 63 from left to right. The diverging colour map goes from -5\% (blue) to +5\% (red), with white being an exact match. Red indicates that the algorithm has overestimated the value, whilst blue indicates underestimation. We can see that the SAGA-SOS algorithm has converged more closely than SAGA-1 to the reference solution without a clear dependence on the number of subsets. 
    }
    \label{fig:PET-rel-diff}
\end{figure}

\begin{figure}
    \centering
    \includegraphics[width=0.48\linewidth]{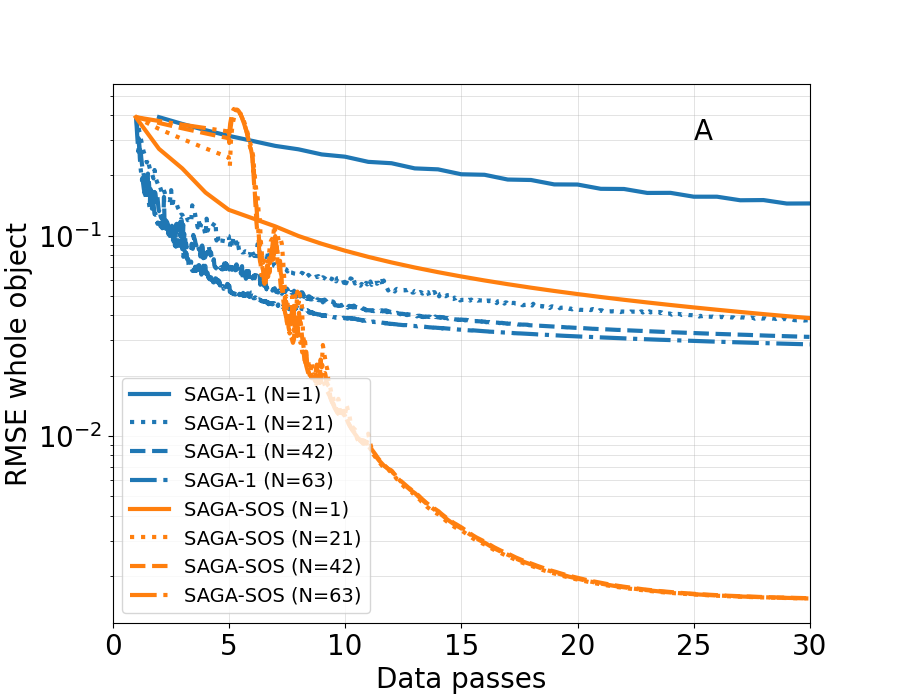}
    \includegraphics[width=0.48\linewidth]{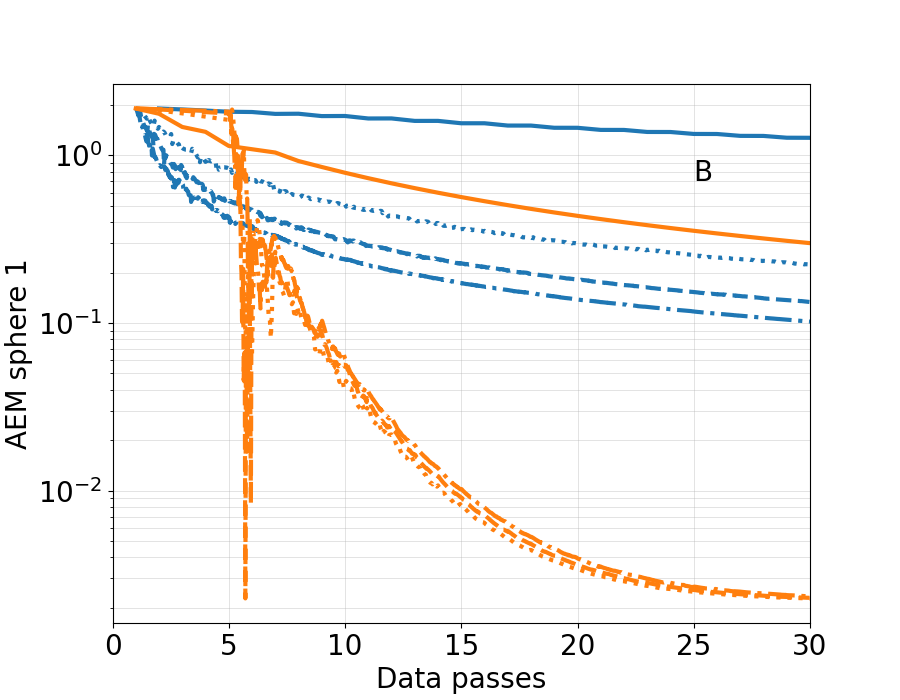}
    \caption{Comparison of the metrics used in the PETRIC challenge as a function of data passes: the RMSE whole object (A) and AEM for the sphere 1 (B) for SAGA-1 and SAGA-SOS submission with 1, 21, 42, and 63 subsets. One subset is a deterministic algorithm, added here as a reference to demonstrate the potential of stochastic methods. These metrics compare the reconstructions obtained by the algorithms with a reference image obtained by the organisers running a different algorithm to "convergence" on the same optimisation problem. The lower the metric, the closer the solution is to the reference solution. It can be seen that SAGA-SOS converges faster than SAGA-1 to the reference solution in both metrics.
    }
    \label{fig:PET-PETRIC-metrics}
\end{figure}

We compared reconstructions obtained by SAGA-1 and SAGA-SOS of the low-count acquisition of the NEMA phantom shown in Figure~\ref{fig:nema_pet}, with 3 different configurations with increasing number of subsets, 21, 42 and 63. We also ran SAGA-SOS and SAGA-1 with one single subset, $N=1$, resulting in a deterministic version of the respective algorithms. In Figure~\ref{fig:PET-rel-diff}, we show the relative difference of the reconstructions obtained by the SAGA-SOS and SAGA-1 algorithms at 30 data passes with respect to the reference solution for the whole volume shown in Figure~\ref{fig:nema_pet}.
In Figure~\ref{fig:PET-PETRIC-metrics}, we show as a function of data passes, the results of a selection of the metrics used during the evaluation of the PETRIC challenge, i.e. the root mean square error RMSE of the whole object and the absolute error of the mean (AEM) for the smallest sphere. These metrics compare the solution reached by the algorithm to a reference solution (solving the same optimisation problem) provided by the organisers of the PETRIC challenge. The exact definition of the metrics and reference algorithm can be found on the PETRIC website,~\cite{PETRIC_arxiv}, but in both cases, a lower value corresponds to a solution that is closer to the reference solution. From all these figures, we can see that SAGA-SOS achieves faster convergence to the reference solution than SAGA-1. One can also observe that SAGA-SOS does not show a strong dependence on the number of subsets, whilst SAGA-1 seems to improve slightly with increasing subsets. We can also see that the deterministic version of the algorithm is always outperformed by its stochastic counterpart, see Figure \ref{fig:PET-PETRIC-metrics}.

\subsubsection{Conclusion}
This case study demonstrates the flexibility of the CIL stochastic framework to tackle a PET problem using SIRF in conjunction; we demonstrated how to use and define custom preconditioners and step size rules with CIL's stochastic framework; we observed that the performance of the SAGA algorithms is not an intrinsic property, but it is deeply influenced by the global configuration, such as the choice of step size and preconditioning. Similar results have been observed during the PETRIC challenge for other algorithms such as SVRG, see ~\cite{Matthiasetal2025}.

\subsection{Code efficiency and speed up}

So far, we compared the convergence speed of different algorithms with respect to data passes, i.e. the computational complexity. However, other parameters have an effect on the actual run time of these algorithms: memory access and transfer between RAM and CPU and to the GPU, parallelisation of the underlying code, to name a few. We report here on the speed of these algorithms with respect to wall clock time.

\begin{figure}
\centering
    \includegraphics[width=0.7\linewidth]{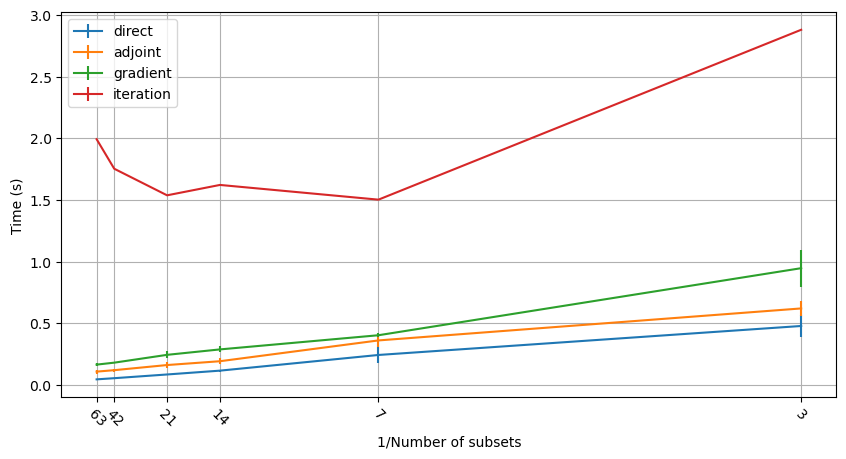}
    \caption{Time to calculate the gradient, direct, adjoint and time per iteration of the SAGA algorithm without prior.
    }
    \label{fig:PET-timings}
\end{figure}

In Figure~\ref{fig:PET-timings} we plot the (wall-clock) time for calling the PET operator's \code{direct} and \code{adjoint}, the gradient of the data discrepancy function and finally the wall-clock time per data pass of the SAGA-1 algorithm, versus the inverse of the number of subsets in the problem. Both the PET operator and gradient of the data discrepancy function are dominated by data computation, and the expected linear behaviour of the execution time vs inverse of the number of subsets is visible in the curves. However, the whole runtime of a single iteration of the SAGA-1 algorithm contains 3 algebraic operations of the type \code{AX+BY} on gradient images (3 times the image data size), one memory copy of the gradient and the evaluation of the proximal of the prior; in this case it is visible that these additional tasks dominate the runtime (in SIRF 3.8) and since the prior depends only on the reconstructed volume, there is no reduced computational cost by increasing the number of subsets. We observe that the minimum runtime is between 7 and 21 subsets, and it is larger for either larger or smaller numbers of subsets. In particular, running with 21 subsets is 25\% faster in terms of wall clock time per iteration than with 63.

%% file: version3/06_discussion.tex
\section{Discussion}\label{sec:discussion}

In this paper, we have demonstrated the CIL stochastic framework: partitioning, sampling, and defining stochastic gradient estimators, which are then plugged into deterministic algorithms. We also illustrated how these algorithms can be further customised through choices of step size, preconditioners, and momentum methods.  In the CT case, staggered partitioning, where each subset contains information from across the projections, was shown to be beneficial compared to sequential sampling. Additionally, increasing the number of subsets improved the convergence rate when measured in data passes. The stochastic algorithms consistently converged more quickly, in terms of data passes, to an ``optimal solution" compared to their deterministic counterparts. In the PET case, we highlighted the impact of selecting appropriate preconditioners and step size strategies, demonstrating how these choices can significantly enhance algorithm performance. The PET example also showcased the ability of CIL and the stochastic framework to integrate with other open-source software and applications.

This paper was primarily written to demonstrate the stochastic framework implemented in CIL, with a focus on its plug-and-play nature, flexibility, and extensibility. Rather than conducting an exhaustive comparison of algorithms, our aim was to showcase how the framework can be applied to real-world data, giving users a sense of the kinds of results they might expect on their own datasets. As such, we did not perform extensive parameter sweeps or fine-tune algorithmic performance. However, the framework itself provides the necessary infrastructure to support future studies dedicated to systematic comparisons and benchmarking of stochastic optimisation methods.

The design of the framework prioritises ease of use and rapid prototyping over computational efficiency. While this makes it highly adaptable across different applications and methods, it also means that algorithms tailored and optimised for specific tasks may outperform those implemented within the framework in terms of speed. Ongoing development efforts include integration of CIL data containers with array APIs, including PyTorch tensors, and will increase the speed of prototyping on small data.  

Although the computational burden of each iteration can be reduced by using only a subset of the data, the benefit in terms of time to convergence needs to account for other effects, such as consistency between the subsets, step size limitations for many stochastic algorithms, and data transfer times from RAM to CPU/GPU. These costs can take a toll on the advantages of using subsets and need to be evaluated in each case.

Another important consideration is memory usage. Stochastic gradient estimators, particularly in large-scale imaging problems, can be memory-intensive, posing challenges for real-data applications. Future work could integrate more memory-efficient algorithms into CIL to help mitigate this issue. 

%% file: version3/07_conclusion.tex
\section{Conclusion} \label{sec:conclusion}

This work presents a flexible framework for implementing a wide range of stochastic algorithms by combining stochastic gradient estimators with deterministic solvers, and enhancing them through preconditioning, momentum, and adaptive step size strategies. Demonstrated on realistic PET and CT datasets, the approach proves effective not only for prototyping and algorithm development but also for deployment in real-world scenarios. This modular framework simplifies both software development and code maintenance, while offering a flexible interface for user-defined experimentation. It enables users to construct their own stochastic estimators and extend the base CIL algorithms to implement new stochastic optimisation methods. Such adaptability is crucial, especially given the rapid pace at which these algorithms are evolving.  We hope that it will serve as an invaluable tool for researchers and practitioners in imaging inverse problems, and that it will also encourage broader adoption of stochastic methods across imaging and inverse problems, fostering innovation and practical impact in the field.

%% file: refs.bib
@article{Nikolova2004,
  author = {Nikolova, Mila},
  description = {SpringerLink - Journal Article},
  journal = {Journal of Mathematical Imaging and Vision},
  month = {jan},
  number = 1,
  pages = {99--120},
  title = {A Variational Approach to Remove Outliers and Impulse Noise},
  url = {http://dx.doi.org/10.1023/B:JMIV.0000011326.88682.e5},
  volume = 20,
  year = 2004
}

@article{Le2007,
  title = {A Variational Approach to Reconstructing Images Corrupted by Poisson Noise},
  volume = {27},
  ISSN = {1573-7683},
  url = {http://dx.doi.org/10.1007/s10851-007-0652-y},
  DOI = {10.1007/s10851-007-0652-y},
  number = {3},
  journal = {Journal of Mathematical Imaging and Vision},
  publisher = {Springer Science and Business Media LLC},
  author = {Le,  Triet and Chartrand,  Rick and Asaki,  Thomas J.},
  year = {2007},
  month = mar,
  pages = {257–263}
}

@article{Rudin1992,
  author ={L. Rudin and S. Osher and E. Fatemi},
  title ={Nonlinear total variation based noise removal algorithms},
  journal = {Physica D: Nonlinear Phenomena},
  year = {1992},
  volume ={60},
  pages = {259-268}, 
  doi={10.1016/0167-2789(92)90242-F},
  url={http://dx.doi.org/10.1016/0167-2789(92)90242-F}
}

@book{tikhonov1977solutions,
  address = {Washington, D.C.: John Wiley \& Sons, New York},
  author = {Tikhonov, Andrey N. and Arsenin, Vasiliy Y.},
  description = {MR: Publications results for "MR Number=(455365)"},
  mrclass = {65J05 (65R05 65NXX)},
  mrnumber = {0455365 (56 \#13604)},
  mrreviewer = {M. Z. Nashed},
  note = {Translated from the Russian, Preface by translation editor Fritz John, Scripta Series in Mathematics},
  pages = {xiii+258},
  publisher = {V. H. Winston \& Sons},
  title = {Solutions of ill-posed problems},
  year = 1977
}

@article{Bredies2010,
  doi = {10.1137/090769521},
  url = {https://doi.org/10.1137/090769521},
  year = {2010},
  month = jan,
  publisher = {Society for Industrial {\&} Applied Mathematics ({SIAM})},
  volume = {3},
  number = {3},
  pages = {492--526},
  author = {Kristian Bredies and Karl Kunisch and Thomas Pock},
  title = {Total Generalized Variation},
  journal = {{SIAM} Journal on Imaging Sciences}
}

@ARTICLE{Holt2014,
  author={Holt, Kevin M.},
  journal={IEEE Trans. Image Process.}, 
  title={Total Nuclear Variation and {J}acobian Extensions of Total Variation for Vector Fields}, 
  year={2014},
  volume={23},
  number={9},
  pages={3975-3989},
  url = {http://dx.doi.org/10.1109/TIP.2014.2332397},
  DOI = {10.1109/tip.2014.2332397}
}

@InProceedings{Lefkimmiatis2013,
author="Lefkimmiatis, Stamatios
and Roussos, Anastasios
and Unser, Michael
and Maragos, Petros",
title="Convex Generalizations of Total Variation Based on the Structure Tensor with Applications to Inverse Problems",
booktitle="Scale Space and Variational Methods in Computer Vision",
year="2013",
pages="48--60",
url = {http://dx.doi.org/10.1007/978-3-642-38267-3_5},
  DOI = {10.1007/978-3-642-38267-3_5},
}

@article{Condat2013,
author = {Condat, Laurent},
doi = {10.1007/s10957-012-0245-9},
isbn = {00223239 (ISSN)},
issn = {00223239},
journal = {Journal of Optimization Theory and Applications},
number = {2},
pages = {460--479},
title = {{A Primal-Dual Splitting Method for Convex Optimization Involving Lipschitzian, Proximable and Linear Composite Terms}},
volume = {158},
year = {2013},
url = {https://doi.org/10.1007/s10957-012-0245-9},
}

@article{Vu2011,
  title = {A splitting algorithm for dual monotone inclusions involving cocoercive operators},
  volume = {38},
  ISSN = {1572-9044},
  url = {http://dx.doi.org/10.1007/s10444-011-9254-8},
  DOI = {10.1007/s10444-011-9254-8},
  number = {3},
  journal = {Advances in Computational Mathematics},
  publisher = {Springer Science and Business Media LLC},
  author = {V\~u, Bang C{\^{o}}ng},
  year = {2011},
  month = nov,
  pages = {667–681}
}

@article{ChambollePock2016, 
title={An introduction to continuous optimization for imaging},
volume={25}, 
journal={Acta Numer.}, 
author={Chambolle, Antonin and Pock, Thomas}, 
year={2016}, 
pages={161-319}}

@article{BenningBurger2018, 
title={Modern regularization methods for inverse problems}, 
volume={27}, 
journal={Acta Numer.}, author={Benning, Martin and Burger, Martin},
year={2018}, pages={1–111}}

@article{Bredies2020,
  title = {Higher-order total variation approaches and generalisations},
  volume = {36},
  ISSN = {1361-6420},
  url = {http://dx.doi.org/10.1088/1361-6420/ab8f80},
  DOI = {10.1088/1361-6420/ab8f80},
  number = {12},
  journal = {Inverse Problems},
  publisher = {IOP Publishing},
  author = {Bredies,  Kristian and Holler,  Martin},
  year = {2020},
  month = dec,
  pages = {123001}
}

@book{Mueller2012,
  title = {Linear and Nonlinear Inverse Problems with Practical Applications},
  ISBN = {9781611972344},
  url = {http://dx.doi.org/10.1137/1.9781611972344},
  DOI = {10.1137/1.9781611972344},
  publisher = {Society for Industrial and Applied Mathematics},
  author = {Mueller,  Jennifer L. and Siltanen,  Samuli},
  year = {2012},
  month = oct 
}

@book{Ryu2022,
  title = {Large-Scale Convex Optimization: Algorithms {\&}; Analyses via Monotone Operators},
  ISBN = {9781009160858},
  url = {http://dx.doi.org/10.1017/9781009160865},
  DOI = {10.1017/9781009160865},
  publisher = {Cambridge University Press},
  author = {Ryu,  Ernest K. and Yin,  Wotao},
  year = {2022},
  month = nov 
}

@article{Yan2018,
  title = {A New Primal–Dual Algorithm for Minimizing the Sum of Three Functions with a Linear Operator},
  volume = {76},
  ISSN = {1573-7691},
  url = {http://dx.doi.org/10.1007/s10915-018-0680-3},
  DOI = {10.1007/s10915-018-0680-3},
  number = {3},
  journal = {Journal of Scientific Computing},
  publisher = {Springer Science and Business Media LLC},
  author = {Yan,  Ming},
  year = {2018},
  month = mar,
  pages = {1698–1717}
}

@article{Yan2024,
  title = {On the Improved Conditions for Some Primal-Dual Algorithms},
  volume = {99},
  ISSN = {1573-7691},
  url = {http://dx.doi.org/10.1007/s10915-024-02537-x},
  DOI = {10.1007/s10915-024-02537-x},
  number = {3},
  journal = {Journal of Scientific Computing},
  publisher = {Springer Science and Business Media LLC},
  author = {Yan,  Ming and Li,  Yao},
  year = {2024},
  month = may 
}

@article{Chen2016,
  title = {A primal-dual fixed point algorithm for minimization of the sum of three convex separable functions},
  volume = {2016},
  ISSN = {1687-1812},
  url = {http://dx.doi.org/10.1186/s13663-016-0543-2},
  DOI = {10.1186/s13663-016-0543-2},
  number = {1},
  journal = {Fixed Point Theory and Applications},
  publisher = {Springer Science and Business Media LLC},
  author = {Chen,  Peijun and Huang,  Jianguo and Zhang,  Xiaoqun},
  year = {2016},
  month = apr 
}

@article{Salim2022,
  title = {Dualize,  Split,  Randomize: Toward Fast Nonsmooth Optimization Algorithms},
  volume = {195},
  ISSN = {1573-2878},
  url = {http://dx.doi.org/10.1007/s10957-022-02061-8},
  DOI = {10.1007/s10957-022-02061-8},
  number = {1},
  journal = {Journal of Optimization Theory and Applications},
  publisher = {Springer Science and Business Media LLC},
  author = {Salim,  Adil and Condat,  Laurent and Mishchenko,  Konstantin and Richtárik,  Peter},
  year = {2022},
  month = jul,
  pages = {102–130}
}

@article{Qian2021,
  author  = {Xun Qian and Zheng Qu and Peter Richt{\'a}rik},
  title   = {{L-SVRG and L-Katyusha with Arbitrary Sampling}},
  journal = {Journal of Machine Learning Research},
  year    = {2021},
  volume  = {22},
  number  = {112},
  pages   = {1--47},
  url     = {http://jmlr.org/papers/v22/20-156.html}
}

@article{Esser2010,
  title = {A General Framework for a Class of First Order Primal-Dual Algorithms for Convex Optimization in Imaging Science},
  volume = {3},
  ISSN = {1936-4954},
  url = {http://dx.doi.org/10.1137/09076934X},
  DOI = {10.1137/09076934x},
  number = {4},
  journal = {SIAM Journal on Imaging Sciences},
  publisher = {Society for Industrial & Applied Mathematics (SIAM)},
  author = {Esser,  Ernie and Zhang,  Xiaoqun and Chan,  Tony F.},
  year = {2010},
  month = jan,
  pages = {1015–1046}
}

@article{Jiang2022,
  title = {Bregman Three-Operator Splitting Methods},
  volume = {196},
  ISSN = {1573-2878},
  url = {http://dx.doi.org/10.1007/s10957-022-02125-9},
  DOI = {10.1007/s10957-022-02125-9},
  number = {3},
  journal = {Journal of Optimization Theory and Applications},
  publisher = {Springer Science and Business Media LLC},
  author = {Jiang,  Xin and Vandenberghe,  Lieven},
  year = {2022},
  month = nov,
  pages = {936–972}
}

@ARTICLE{Fadili2011,
  author={Fadili, Jalal M. and Peyré, Gabriel},
  journal={IEEE Transactions on Image Processing}, 
  title={Total Variation Projection With First Order Schemes}, 
  year={2011},
  volume={20},
  number={3},
  pages={657-669},
  doi={10.1109/TIP.2010.2072512},
  url = {https://doi.org/10.1109/TIP.2010.2072512},
}

@article{Latafat2017,
  title = {Asymmetric forward–backward–adjoint splitting for solving monotone inclusions involving three operators},
  volume = {68},
  ISSN = {1573-2894},
  url = {http://dx.doi.org/10.1007/s10589-017-9909-6},
  DOI = {10.1007/s10589-017-9909-6},
  number = {1},
  journal = {Computational Optimization and Applications},
  publisher = {Springer Science and Business Media LLC},
  author = {Latafat,  Puya and Patrinos,  Panagiotis},
  year = {2017},
  month = apr,
  pages = {57–93}
}

@article{Daubechies2004,
  title = {An iterative thresholding algorithm for linear inverse problems with a sparsity constraint},
  volume = {57},
  ISSN = {1097-0312},
  url = {http://dx.doi.org/10.1002/cpa.20042},
  DOI = {10.1002/cpa.20042},
  number = {11},
  journal = {Communications on Pure and Applied Mathematics},
  publisher = {Wiley},
  author = {Daubechies,  I. and Defrise,  M. and De Mol,  C.},
  year = {2004},
  month = aug,
  pages = {1413–1457}
}

@article{Nesterov1983,
  title = {A method of solving a convex programming problem with convergence rate $O\bigl(\frac1{k^2}\bigr)$},
  volume = {269},
  number = {3},
  journal = {Dokl. Akad. Nauk SSSR},
  author = {Yu.~E.~Nesterov},
  year = {1983},
  pages = {543--547}
}

@article{Nesterov2004,
  title = {Smooth minimization of non-smooth functions},
  volume = {103},
  ISSN = {1436-4646},
  url = {http://dx.doi.org/10.1007/s10107-004-0552-5},
  DOI = {10.1007/s10107-004-0552-5},
  number = {1},
  journal = {Mathematical Programming},
  publisher = {Springer Science and Business Media LLC},
  author = {Nesterov,  Yu.},
  year = {2004},
  month = dec,
  pages = {127–152}
}

@article{Nesterov2012,
  title = {Gradient methods for minimizing composite functions},
  volume = {140},
  ISSN = {1436-4646},
  url = {http://dx.doi.org/10.1007/s10107-012-0629-5},
  DOI = {10.1007/s10107-012-0629-5},
  number = {1},
  journal = {Mathematical Programming},
  publisher = {Springer Science and Business Media LLC},
  author = {Nesterov,  Yu.},
  year = {2012},
  month = dec,
  pages = {125–161}
}

@article{Taylor2017,
  title = {Exact Worst-Case Performance of First-Order Methods for Composite Convex Optimization},
  volume = {27},
  ISSN = {1095-7189},
  url = {http://dx.doi.org/10.1137/16M108104X},
  DOI = {10.1137/16m108104x},
  number = {3},
  journal = {SIAM Journal on Optimization},
  publisher = {Society for Industrial & Applied Mathematics (SIAM)},
  author = {Taylor,  Adrien B. and Hendrickx,  Julien M. and Glineur,  Fran\c{c}ois},
  year = {2017},
  month = jan,
  pages = {1283–1313}
}

@article{Aspremont2021,
author = {Alexandre d'Aspremont and Damien Scieur and Adrien Taylor},
title = {Acceleration Methods},
Year = {2021},
Eprint = {arXiv:2101.09545},
journal = {Foundations and Trends in Optimization: Vol. 5: No. 1-2, pp 1-245 (2021)},
Doi = {10.1561/2400000036},
  url = {https://doi.org/10.1561/2400000036},
}

@article{Kaufman1993,
  title = {Maximum likelihood,  least squares,  and penalized least squares for PET},
  volume = {12},
  ISSN = {0278-0062},
  url = {http://dx.doi.org/10.1109/42.232249},
  DOI = {10.1109/42.232249},
  number = {2},
  journal = {IEEE Transactions on Medical Imaging},
  publisher = {Institute of Electrical and Electronics Engineers (IEEE)},
  author = {Kaufman,  L.},
  year = {1993},
  month = jun,
  pages = {200–214}
}

@article{Clinthorne1993,
  title = {Preconditioning methods for improved convergence rates in iterative reconstructions},
  volume = {12},
  ISSN = {0278-0062},
  url = {http://dx.doi.org/10.1109/42.222670},
  DOI = {10.1109/42.222670},
  number = {1},
  journal = {IEEE Transactions on Medical Imaging},
  publisher = {Institute of Electrical and Electronics Engineers (IEEE)},
  author = {Clinthorne,  N.H. and Pan,  T.-S. and Chiao,  P.-C. and Rogers,  W.L. and Stamos,  J.A.},
  year = {1993},
  month = mar,
  pages = {78–83}
}

@article{Qi2006,
  title = {Iterative reconstruction techniques in emission computed tomography},
  volume = {51},
  ISSN = {1361-6560},
  url = {http://dx.doi.org/10.1088/0031-9155/51/15/R01},
  DOI = {10.1088/0031-9155/51/15/r01},
  number = {15},
  journal = {Physics in Medicine and Biology},
  publisher = {IOP Publishing},
  author = {Qi,  Jinyi and Leahy,  Richard M},
  year = {2006},
  month = jul,
  pages = {R541–R578}
}

@Article{Chambolle2004,
author={Chambolle, Antonin},
title={An Algorithm for Total Variation Minimization and Applications},
journal={J. Math. Imaging Vis.},
year={2004},
day={01},
volume={20},
number={1},
pages={89-97}
}

@article{Robbins1951,
  title = {A Stochastic Approximation Method},
  volume = {22},
  ISSN = {0003-4851},
  url = {http://dx.doi.org/10.1214/aoms/1177729586},
  DOI = {10.1214/aoms/1177729586},
  number = {3},
  journal = {The Annals of Mathematical Statistics},
  publisher = {Institute of Mathematical Statistics},
  author = {Robbins,  Herbert and Monro,  Sutton},
  year = {1951},
  month = sep,
  pages = {400–407}
}

@article{Bottou2018,
  title = {Optimization Methods for Large-Scale Machine Learning},
  volume = {60},
  ISSN = {1095-7200},
  url = {http://dx.doi.org/10.1137/16M1080173},
  DOI = {10.1137/16m1080173},
  number = {2},
  journal = {SIAM Review},
  publisher = {Society for Industrial & Applied Mathematics (SIAM)},
  author = {Bottou,  Léon and Curtis,  Frank E. and Nocedal,  Jorge},
  year = {2018},
  month = jan,
  pages = {223–311}
}

@article{Schmidt2017,
author={Schmidt, Mark
and Le Roux, Nicolas
and Bach, Francis},
title={Minimizing finite sums with the stochastic average gradient},
journal={Mathematical Programming},
year={2017},
month={Mar},
day={01},
volume={162},
number={1},
pages={83-112},
issn={1436-4646},
doi={10.1007/s10107-016-1030-6},
url={https://doi.org/10.1007/s10107-016-1030-6}
}

@article{Blatt2007,
  title = {A Convergent Incremental Gradient Method with a Constant Step Size},
  volume = {18},
  ISSN = {1095-7189},
  url = {http://dx.doi.org/10.1137/040615961},
  DOI = {10.1137/040615961},
  number = {1},
  journal = {SIAM Journal on Optimization},
  publisher = {Society for Industrial & Applied Mathematics (SIAM)},
  author = {Blatt,  Doron and Hero,  Alfred O. and Gauchman,  Hillel},
  year = {2007},
  month = jan,
  pages = {29–51}
}

@article{Driggs2022,
  author  = {Derek Driggs and Jingwei Liang and Carola-Bibiane Schonlieb},
  title   = {On Biased Stochastic Gradient Estimation},
  journal = {Journal of Machine Learning Research},
  year    = {2022},
  volume  = {23},
  number  = {24},
  pages   = {1--43},
  url     = {http://jmlr.org/papers/v23/20-316.html}
}

@inproceedings{Defazio2014,
   author = {Aaron Defazio and Francis Bach and Simon Lacoste-Julien},
   journal = {Advances in Neural Information Processing Systems},
   pages = {1646-1654},
   title = {SAGA: A fast incremental gradient method with support for non-strongly convex composite objectives},
   volume = {2},
  booktitle={Advances in Neural Information Processing Systems},
   year = {2014},
}

@article{Karimi2016,
  title = {A hybrid stochastic-deterministic gradient descent algorithm for image reconstruction in cone-beam computed tomography},
  volume = {2},
  ISSN = {2057-1976},
  url = {http://dx.doi.org/10.1088/2057-1976/2/1/015008},
  DOI = {10.1088/2057-1976/2/1/015008},
  number = {1},
  journal = {Biomedical Physics {\&}; Engineering Express},
  publisher = {IOP Publishing},
  author = {Karimi,  Davood and Ward,  Rabab K},
  year = {2016},
  month = feb,
  pages = {015008}
}

@inproceedings{Johnson2013,
 author = {Johnson, Rie and Zhang, Tong},
 booktitle = {Advances in Neural Information Processing Systems},
 editor = {C.J. Burges and L. Bottou and M. Welling and Z. Ghahramani and K.Q. Weinberger},
 pages = {},
 publisher = {Curran Associates, Inc.},
 title = {Accelerating Stochastic Gradient Descent using Predictive Variance Reduction},
 url = {https://proceedings.neurips.cc/paper_files/paper/2013/file/ac1dd209cbcc5e5d1c6e28598e8cbbe8-Paper.pdf},
 volume = {26},
 year = {2013}
}

@article{Xiao2014,
  title = {A Proximal Stochastic Gradient Method with Progressive Variance Reduction},
  volume = {24},
  ISSN = {1095-7189},
  url = {http://dx.doi.org/10.1137/140961791},
  DOI = {10.1137/140961791},
  number = {4},
  journal = {SIAM Journal on Optimization},
  publisher = {Society for Industrial & Applied Mathematics (SIAM)},
  author = {Xiao,  Lin and Zhang,  Tong},
  year = {2014},
  month = jan,
  pages = {2057–2075}
}

@article{Kovalev2020,
  title = 	 {Don't Jump Through Hoops and Remove Those Loops:  SVRG and Katyusha are Better Without the Outer Loop},
  author =       {Kovalev, Dmitry and Horv{\'a}th, Samuel and Richt{\'a}rik, Peter},
  journal = 	 {Proceedings of the 31st International Conference  on Algorithmic Learning Theory},
  pages = 	 {451--467},
  year = 	 {2020},
  editor = 	 {Kontorovich, Aryeh and Neu, Gergely},
  volume = 	 {117},
  series = 	 {Proceedings of Machine Learning Research},
  month = 	 {08 Feb--11 Feb},
  publisher =    {PMLR},
  url = 	 {https://proceedings.mlr.press/v117/kovalev20a.html}
}

@InProceedings{pmlr-v124-zhou20a,
  title = 	 {{Amortized Nesterov’s Momentum: A Robust Momentum and Its Application to Deep Learning}},
  author =       {Zhou, Kaiwen and Jin, Yanghua and Ding, Qinghua and Cheng, James},
  booktitle = 	 {Proceedings of the 36th Conference on Uncertainty in Artificial Intelligence (UAI)},
  pages = 	 {211--220},
  year = 	 {2020},
  editor = 	 {Peters, Jonas and Sontag, David},
  volume = 	 {124},
  series = 	 {Proceedings of Machine Learning Research},
  month = 	 {03--06 Aug},
  publisher =    {PMLR},
  url = 	 {https://proceedings.mlr.press/v124/zhou20a.html}
}

@article{JMLR:v10:duchi09a,
  author  = {John Duchi and Yoram Singer},
  title   = {Efficient Online and Batch Learning Using Forward Backward Splitting},
  journal = {Journal of Machine Learning Research},
  year    = {2009},
  volume  = {10},
  number  = {99},
  pages   = {2899--2934},
  url     = {http://jmlr.org/papers/v10/duchi09a.html}
}

@inproceedings{Hofmann2015,
 author = {Hofmann, Thomas and Lucchi, Aurelien and Lacoste-Julien, Simon and McWilliams, Brian},
 booktitle = {Advances in Neural Information Processing Systems},
 editor = {C. Cortes and N. Lawrence and D. Lee and M. Sugiyama and R. Garnett},
 pages = {},
 publisher = {Curran Associates, Inc.},
 title = {Variance Reduced Stochastic Gradient Descent with Neighbors},
 url = {https://proceedings.neurips.cc/paper_files/paper/2015/file/effc299a1addb07e7089f9b269c31f2f-Paper.pdf},
 volume = {28},
 year = {2015}
}

@InProceedings{Bingcong2019,
author = {Bingcong Li and Meng Ma and Georgios B. Giannakis},
title = {{On the Convergence of SARAH and Beyond}},
year = {2019},
eprint = {arXiv:1906.02351},
url = {https://arxiv.org/abs/1906.02351},
booktitle = {Proceedings of Machine Learning Research},
}

@InProceedings{Nhan2019,
author = {Nhan H. Pham and Lam M. Nguyen and Dzung T. Phan and Quoc Tran-Dinh},
title = {ProxSARAH: An Efficient Algorithmic Framework for Stochastic Composite Nonconvex Optimization},
year = {2019},
eprint = {arXiv:1902.05679},
url = {https://arxiv.org/abs/1902.05679},
booktitle = {Journal of Machine Learning Research},
}

@book{Nesterov2004book,
  title = {Introductory Lectures on Convex Optimization},
  ISBN = {9781441988539},
  ISSN = {1384-6485},
  url = {http://dx.doi.org/10.1007/978-1-4419-8853-9},
  DOI = {10.1007/978-1-4419-8853-9},
  journal = {Applied Optimization},
  publisher = {Springer US},
  author = {Nesterov,  Yurii},
  year = {2004}
}

@article{Driggs2020,
  title = {Accelerating variance-reduced stochastic gradient methods},
  volume = {191},
  ISSN = {1436-4646},
  url = {http://dx.doi.org/10.1007/s10107-020-01566-2},
  DOI = {10.1007/s10107-020-01566-2},
  number = {2},
  journal = {Mathematical Programming},
  publisher = {Springer Science and Business Media LLC},
  author = {Driggs,  Derek and Ehrhardt,  Matthias J. and Sch\"{o}nlieb,  Carola-Bibiane},
  year = {2020},
  month = sep,
  pages = {671–715}
}

@inproceedings{Yurtsever2016,
 author = {Yurtsever, Alp and Vu, Bang Cong and Cevher, Volkan},
 booktitle = {Advances in Neural Information Processing Systems},
 editor = {D. Lee and M. Sugiyama and U. Luxburg and I. Guyon and R. Garnett},
 pages = {},
 publisher = {Curran Associates, Inc.},
 title = {Stochastic Three-Composite Convex Minimization},
 url = {https://proceedings.neurips.cc/paper_files/paper/2016/file/5d6646aad9bcc0be55b2c82f69750387-Paper.pdf},
 volume = {29},
 year = {2016}
}

@InProceedings{Zhao2019,
  title = 	 {An Optimal Algorithm for Stochastic Three-Composite Optimization},
  author =       {Zhao, Renbo and Haskell, William B. and Tan, Vincent Y. F.},
  booktitle = 	 {Proceedings of the Twenty-Second International Conference on Artificial Intelligence and Statistics},
  pages = 	 {428--437},
  year = 	 {2019},
  editor = 	 {Chaudhuri, Kamalika and Sugiyama, Masashi},
  volume = 	 {89},
  series = 	 {Proceedings of Machine Learning Research},
  month = 	 {16--18 Apr},
  publisher =    {PMLR},
  url = 	 {https://proceedings.mlr.press/v89/zhao19a.html}
}

@article{Chambolle_2018_SPDHG,
  doi = {10.1137/17m1134834},
  url = {https://doi.org/10.1137/17m1134834},
  month = jan,
  year = {2018},
  publisher = {Society for Industrial {\&} Applied Mathematics ({SIAM})},
  volume = {28},
  number = {4},
  pages = {2783--2808},
  author = {Antonin Chambolle and Matthias J. Ehrhardt and Peter Richt{\'{a}}rik and Carola-Bibiane Sch\"{o}nlieb},
  title = {Stochastic Primal-Dual Hybrid Gradient Algorithm with Arbitrary Sampling and Imaging Applications},
  journal = {{SIAM} Journal on Optimization}
}

@article{Herman1993,
  title = {Algebraic reconstruction techniques can be made computationally efficient (positron emission tomography application)},
  volume = {12},
  ISSN = {0278-0062},
  url = {http://dx.doi.org/10.1109/42.241889},
  DOI = {10.1109/42.241889},
  number = {3},
  journal = {IEEE Transactions on Medical Imaging},
  publisher = {Institute of Electrical and Electronics Engineers (IEEE)},
  author = {Herman,  G.T. and Meyer,  L.B.},
  year = {1993},
  pages = {600–609}
}

@article{Blondel2021,
  title={Efficient and Modular Implicit Differentiation},
  author={Blondel, Mathieu and Berthet, Quentin and Cuturi, Marco and Frostig, Roy 
    and Hoyer, Stephan and Llinares-L{\'o}pez, Felipe and Pedregosa, Fabian 
    and Vert, Jean-Philippe},
  journal={arXiv preprint arXiv:2105.15183},
  year={2021}
}

@article{BARZILAI1988,
  title = {Two-Point Step Size Gradient Methods},
  volume = {8},
  ISSN = {1464-3642},
  url = {http://dx.doi.org/10.1093/imanum/8.1.141},
  DOI = {10.1093/imanum/8.1.141},
  number = {1},
  journal = {IMA Journal of Numerical Analysis},
  publisher = {Oxford University Press (OUP)},
  author = {Barzilai, Jonathan and Borwein,  Jonathan M.},
  year = {1988},
  pages = {141–148}
}

@article{Vamvakeros2021,
  title = {Cycling Rate‐Induced Spatially‐Resolved Heterogeneities in Commercial Cylindrical Li‐Ion Batteries},
  volume = {5},
  ISSN = {2366-9608},
  url = {http://dx.doi.org/10.1002/smtd.202100512},
  DOI = {10.1002/smtd.202100512},
  number = {9},
  journal = {Small Methods},
  publisher = {Wiley},
  author = {Vamvakeros,  Antonis and Matras,  Dorota and Ashton,  Thomas E. and Coelho,  Alan A. and Dong,  Hongyang and Bauer,  Dustin and Odarchenko,  Yaroslav and Price,  Stephen W. T. and Butler,  Keith T. and Gutowski,  Olof and Dippel,  Ann‐Christin and Zimmerman,  Martin von and Darr,  Jawwad A. and Jacques,  Simon D. M. and Beale,  Andrew M.},
  year = {2021},
  month = aug 
}

@article{Drakopoulos2015,
  title = {I12: the Joint Engineering,  Environment and Processing (JEEP) beamline at Diamond Light Source},
  volume = {22},
  ISSN = {1600-5775},
  url = {http://dx.doi.org/10.1107/S1600577515003513},
  DOI = {10.1107/s1600577515003513},
  number = {3},
  journal = {Journal of Synchrotron Radiation},
  publisher = {International Union of Crystallography (IUCr)},
  author = {Drakopoulos,  Michael and Connolley,  Thomas and Reinhard,  Christina and Atwood,  Robert and Magdysyuk,  Oxana and Vo,  Nghia and Hart,  Michael and Connor,  Leigh and Humphreys,  Bob and Howell,  George and Davies,  Steve and Hill,  Tim and Wilkin,  Guy and Pedersen,  Ulrik and Foster,  Andrew and De Maio,  Nicoletta and Basham,  Mark and Yuan,  Fajin and Wanelik,  Kaz},
  year = {2015},
  month = apr,
  pages = {828–838}
}

@misc{Docherty2025,
  title = {{Battery Imaging Library: Multi-length scale and multi-modal synchrotron and laboratory battery imaging data for all}},
  url = {https://chemrxiv.org/doi/full/10.26434/chemrxiv-2025-sbp73/v2},
  publisher = {American Chemical Society (ACS)},
  author = {Docherty,  Ronan and Riley,  Sam and Morley,  John D. and Papoutsellis,  Evangelos and Bobitan,  Antonia Diana and Donoghue,  Jack and Rankin,  Katy and Alvarez Borges,  Fernando and Salinas-Farran,  Luis and Michalik,  Stefan and Liptak,  Alexander and Burca,  Genoveva and Autran,  Pierre-Olivier and Wright,  Jon and Kockelmann,  Winfried and Gutowski,  Olof and Dippel,  Ann-Christin and von Zimmermann,  Martin and Matras,  Dorota and Winiarski,  Bartlomiej and Mangum,  John and Finegan,  Donal and Gott,  James A. and Proprentner,  Daniela and Dong,  Hongyang and George,  Chandramohan and Beale,  Andrew M. and Jacques,  Simon D.M. and Cooper,  Samuel J. and Vamvakeros,  Antonis},
  year = {2025},
  month = sep 
}

@article{Vo2021,
author = {Nghia T. Vo and Robert C. Atwood and Michael Drakopoulos and Thomas Connolley},
journal = {Opt. Express},
number = {12},
pages = {17849--17874},
publisher = {Optica Publishing Group},
title = {Data processing methods and data acquisition for samples larger than the field of view in parallel-beam tomography},
volume = {29},
month = {Jun},
year = {2021},
url = {https://opg.optica.org/oe/abstract.cfm?URI=oe-29-12-17849},
doi = {10.1364/OE.418448},
}

@article{Chambolle2010,
  doi = {10.1007/s10851-010-0251-1},
  url = {https://doi.org/10.1007/s10851-010-0251-1},
  year = {2010},
  month = dec,
  publisher = {Springer Science and Business Media {LLC}},
  volume = {40},
  number = {1},
  pages = {120--145},
  author = {Antonin Chambolle and Thomas Pock},
  title = {A First-Order Primal-Dual Algorithm for Convex Problems with~Applications to Imaging},
  journal = {Journal of Mathematical Imaging and Vision}
}

@article{Beck2009,
  doi = {10.1109/tip.2009.2028250},
  url = {https://doi.org/10.1109/tip.2009.2028250},
  year = {2009},
  month = nov,
  publisher = {Institute of Electrical and Electronics Engineers ({IEEE})},
  volume = {18},
  number = {11},
  pages = {2419--2434},
  author = {A. Beck and M. Teboulle},
  title = {Fast Gradient-Based Algorithms for Constrained Total Variation Image Denoising and Deblurring Problems},
  journal = {{IEEE} Transactions on Image Processing}
}

@article{OConnor2018,
  title = {On the equivalence of the primal-dual hybrid gradient method and Douglas–Rachford splitting},
  volume = {179},
  ISSN = {1436-4646},
  url = {http://dx.doi.org/10.1007/s10107-018-1321-1},
  DOI = {10.1007/s10107-018-1321-1},
  number = {1–2},
  journal = {Mathematical Programming},
  publisher = {Springer Science and Business Media LLC},
  author = {O’Connor,  Daniel and Vandenberghe,  Lieven},
  year = {2018},
  month = aug,
  pages = {85–108}
}

@article{Jorgensen2021,
Author = {Jakob S. J{\o}rgensen and Evelina Ametova and Genoveva Burca and Gemma Fardell and Evangelos Papoutsellis and Edoardo Pasca and Kris Thielemans and Martin Turner and Ryan Warr and William R. B. Lionheart and Philip J. Withers},
Title = {{Core Imaging Library -- Part I: a versatile Python framework for tomographic imaging}},
Year = {2021},
journal = {Phil. Trans. R. Soc. A 20200192},
doi = {10.1098/rsta.2020.0192},
url = {https://doi.org/10.1098/rsta.2020.0192}
}

@article{Papoutsellis2021,
  doi = {10.1098/rsta.2020.0193},
  url = {https://doi.org/10.1098/rsta.2020.0193},
  year = {2021},
  month = jul,
  publisher = {The Royal Society},
  volume = {379},
  number = {2204},
  pages = {20200193},
  author = {Evangelos Papoutsellis and Evelina Ametova and Claire Delplancke and Gemma Fardell and Jakob S. J{\o}rgensen and others},
  title = {{Core Imaging Library} - Part {II}: multichannel reconstruction for dynamic and spectral tomography},
  journal = {Philos Trans A Math Phys Eng Sci}
}

@article{Evelina2021,
Author = {Evelina Ametova and Genoveva Burca and Suren Chilingaryan and Gemma Fardell and Jakob S. J{\o}rgensen and Evangelos Papoutsellis and Edoardo Pasca and Ryan Warr and Martin Turner and William R. B. Lionheart and Philip J. Withers},
Title = {Crystalline phase discriminating neutron tomography using advanced reconstruction methods},
Year = {2021},
journal = {Journal of Physics D: Applied Physics }, 
doi = {10.1088/1361-6463/ac02f9},
url = {https://doi.org/10.1088/1361-6463/ac02f9},
}

@article{SimmendefeldtSchmidt2024,
  title = {Anisotropic regularization for inversion of fast-ion loss detector measurements},
  volume = {64},
  ISSN = {1741-4326},
  url = {http://dx.doi.org/10.1088/1741-4326/ad75a5},
  DOI = {10.1088/1741-4326/ad75a5},
  number = {10},
  journal = {Nuclear Fusion},
  publisher = {IOP Publishing},
  author = {Simmendefeldt Schmidt,  Bo and Sauer Jørgensen,  Jakob and Rueda-Rueda,  José and Galdon-Quíroga,  Joaquín and García-Muñoz,  Manuel and Salewski,  Mirko},
  year = {2024},
  month = sep,
  pages = {106053}
}

@article{Warr2021,
  title = {Enhanced hyperspectral tomography for bioimaging by spatiospectral reconstruction},
  volume = {11},
  ISSN = {2045-2322},
  url = {http://dx.doi.org/10.1038/s41598-021-00146-4},
  DOI = {10.1038/s41598-021-00146-4},
  number = {1},
  journal = {Scientific Reports},
  publisher = {Springer Science and Business Media LLC},
  author = {Warr,  Ryan and Ametova,  Evelina and Cernik,  Robert J. and Fardell,  Gemma and Handschuh,  Stephan and Jørgensen,  Jakob S. and Papoutsellis,  Evangelos and Pasca,  Edoardo and Withers,  Philip J.},
  year = {2021},
  month = oct 
}

@article{Jorgensen2023,
  title = {A directional regularization method for the limited-angle Helsinki Tomography Challenge using the Core Imaging Library (CIL)},
  volume = {1},
  ISSN = {2994-7669},
  url = {http://dx.doi.org/10.3934/ammc.2023011},
  DOI = {10.3934/ammc.2023011},
  number = {2},
  journal = {Applied Mathematics for Modern Challenges},
  publisher = {American Institute of Mathematical Sciences (AIMS)},
  author = {Jørgensen,  Jakob Sauer and Papoutsellis,  Evangelos and Murgatroyd,  Laura and Fardell,  Gemma and Pasca,  Edoardo},
  year = {2023},
  pages = {143–169}
}

@article{Brown2021,
  title={{Motion estimation and correction for simultaneous PET/MR using SIRF and CIL}},
  author={Richard Brown and Christoph Kolbitsch and Claire Delplancke and Evangelos Papoutsellis and Johannes Mayer and Evgueni Ovtchinnikov and Edoardo Pasca and Radhouene Neji and Casper da Costa-Luis and Ashley G. Gillman and  Matthias J. Ehrhardt and Jamie McClelland and Bjoern Eiben and Kris Thielemans},
  year={2021},
  journal={Phil. Trans. R. Soc. A 20200208},
  doi = {10.1098/rsta.2020.0208},
  url = {https://doi.org/10.1098/rsta.2020.0208},
}

@article{Watson2024,
  title = {Resolving full-wave through-wall transmission effects in multi-static synthetic aperture radar},
  volume = {40},
  ISSN = {1361-6420},
  url = {http://dx.doi.org/10.1088/1361-6420/ad5b83},
  DOI = {10.1088/1361-6420/ad5b83},
  number = {8},
  journal = {Inverse Problems},
  publisher = {IOP Publishing},
  author = {Watson,  F M and Andre,  D and Lionheart,  W R B},
  year = {2024},
  month = jul,
  pages = {085009}
}

@article{NEMAstandard,
    author = {National Electrical Manufacturers Association, NEMA Standards Publication NU},
    title = {Performance measurements of positron emission tomographs},
    journal = {Rosslyn, VA},
    year = {2007}
}

@misc{PETRIC_arxiv,
Author = {Casper da Costa-Luis and Matthias J. Ehrhardt and Christoph Kolbitsch and Evgueni Ovtchinnikov and Edoardo Pasca and Kris Thielemans and Charalampos Tsoumpas},
Title = {{PET Rapid Image Reconstruction Challenge (PETRIC)}},
Year = {2025},
Eprint = {arXiv:2511.22566},

  url = {https://arxiv.org/abs/2511.22566},
}

@article{Kazantsev,
  doi = {10.1016/j.softx.2019.04.003},
  url = {https://doi.org/10.1016/j.softx.2019.04.003},
  year = {2019},
  month = jan,
  publisher = {Elsevier {BV}},
  volume = {9},
  pages = {317--323},
  author = {Daniil Kazantsev and Edoardo Pasca and Martin J. Turner and Philip J. Withers},
  title = {{CCPi}-Regularisation toolkit for computed tomographic image reconstruction with proximal splitting algorithms},
  journal = {{SoftwareX}}
}

@article{Ovtchinnikov2020,
  doi = {10.1016/j.cpc.2019.107087},
  url = {https://doi.org/10.1016/j.cpc.2019.107087},
  year = {2020},
  month = apr,
  publisher = {Elsevier {BV}},
  volume = {249},
  pages = {107087},
  author = {E. Ovtchinnikov and R. Brown and C. Kolbitsch and E. Pasca and C. da Costa-Luis and A. G. Gillman and B. A. Thomas and N. Efthimiou and J. Mayer and P. Wadhwa and M. J. Ehrhardt and S. Ellis and J. S. J{\o}rgensen and J. Matthews and C. Prieto and A. J. Reader and C. Tsoumpas and M. J. Turner and D. Atkinson and K. Thielemans},
  title = {{SIRF}: Synergistic Image Reconstruction Framework},
  journal = {Computer Physics Communications}
}

@article{beck2009fast,
doi= {10.1137/080716542},
url={https://doi.org/10.1137/080716542},
  title={A fast iterative shrinkage-thresholding algorithm for linear inverse problems},
  author={Beck, Amir and Teboulle, Marc},
  journal={SIAM journal on imaging sciences},
  volume={2},
  number={1},
  pages={183--202},
  year={2009},
  publisher={SIAM}
}

@article{armijo1966minimization,
  title={Minimization of functions having Lipschitz continuous first partial derivatives},
  author={Armijo, Larry},
  journal={Pacific Journal of mathematics},
  volume={16},
  number={1},
  pages={1--3},
  year={1966},
  publisher={Mathematical Sciences Publishers},
doi={10.2140/pjm.1966.16.1},
url={https://msp.org/pjm/1966/16-1/p01.xhtml}
}

@article{Nuyts2002,
  author={Nuyts J. and Beque D. and Dupont P. and Mortelmans L.},
  journal={IEEE Transactions on Nuclear Science}, 
  title={A concave prior penalizing relative differences for maximum-a-posteriori reconstruction in emission tomography}, 
  year={2002},
  volume={49},
  number={1},
  pages={56-60},
  doi={10.1109/TNS.2002.998681},
  url = {https://doi.org/10.1109/TNS.2002.998681},
}

@article{Ahn2003,
  title={Globally convergent image reconstruction for emission tomography using relaxed ordered subsets algorithms},
  author={Ahn, Sangtae and Fessler, Jeffrey A},
  journal={IEEE transactions on medical imaging},
  volume={22},
  number={5},
  pages={613--626},
  year={2003},
  publisher={IEEE},
  doi={10.1109/TMI.2003.812251}
}

@article{yan2018new,
  title={A new primal--dual algorithm for minimizing the sum of three functions with a linear operator},
  author={Yan, Ming},
  journal={Journal of Scientific Computing},
  volume={76},
  pages={1698--1717},
  year={2018},
  publisher={Springer}
}

@article{Tsai2017,
  author={Tsai, Yu-Jung and Schramm, Georg and Ahn, Sangtae and Bousse, Alexandre and Arridge, Simon and Nuyts, Johan and Hutton, Brian F. and Stearns, Charles W. and Thielemans, Kris},
  journal={IEEE Transactions on Medical Imaging}, 
  title={Benefits of Using a Spatially-Variant Penalty Strength With Anatomical Priors in PET Reconstruction}, 
  year={2020},
  volume={39},
  number={1},
  pages={11-22},
  doi={10.1109/TMI.2019.2913889},
  url = {https://doi.org/10.1109/TMI.2019.2913889},
}

@article{Tsai2020,
  author={Tsai, Yu-Jung and Schramm, Georg and Ahn, Sangtae and Bousse, Alexandre and Arridge, Simon and Nuyts, Johan and Hutton, Brian F. and Stearns, Charles W. and Thielemans, Kris},
  journal={IEEE Transactions on Medical Imaging}, 
  title={Benefits of Using a Spatially-Variant Penalty Strength With Anatomical Priors in PET Reconstruction}, 
  year={2020},
  volume={39},
  number={1},
  pages={11-22},
  doi={10.1109/TMI.2019.2913889},
  url = {https://doi.org/10.1109/TMI.2019.2913889},
}

@misc{PET_partitioner,
howpublished={\url{https://github.com/SyneRBI/SIRF-Contribs/blob/master/src/Python/sirf/contrib/partitioner/partitioner.py}}
}

@article{Matthiasetal2025,
  title={{Fast PET reconstruction with variance reduction and prior-aware preconditioning}},
  author={Ehrhardt, Matthias J and Kereta, Zeljko and Schramm, Georg},
  journal={Frontiers in nuclear medicine},
  volume={5},
  pages={1641215},
  year={2025},
  publisher={Frontiers Media SA},
DOI={10.3389/fnume.2025.1641215},
url={https://doi.org/10.3389/fnume.2025.1641215}
}

@misc{PETRIC_NEMAIQlowcountsDataset,
    author = {Thomas BA and Sanderson T.},
    title = {NEMA image quality phantom acquisition on the Siemens
  mMR scanner},
    doi = {10.5281/zenodo.1304454},
    year = {2018},
  url = {https://doi.org/10.5281/zenodo.1304454},
}

@article{mcir,
	author={Protopapa, Letizia and Duff, Margaret A G and Mayer, Johannes and Schulz-Menger, Jeanette and Thielemans, Kris and Kolbitsch, Christoph and Pasca, Edoardo},
	title={Efficient motion-corrected image reconstruction for 3D cardiac MRI through stochastic optimisation},
	journal={Physics in Medicine {\&} Biology},
	url={http://iopscience.iop.org/article/10.1088/1361-6560/adf609},
	year={2025}
}

@article{Ehrhardt_2025,
doi = {10.1088/1361-6420/adc0b7},
url = {https://dx.doi.org/10.1088/1361-6420/adc0b7},
year = {2025},
month = {may},
publisher = {IOP Publishing},
volume = {41},
number = {5},
pages = {053001},
author = {Ehrhardt, Matthias J and Kereta, Zeljko and Liang, Jingwei and Tang, Junqi},
title = {A guide to stochastic optimisation for large-scale inverse problems},
journal = {Inverse Problems}
}

@book{nocedal2006numerical,
  title={Numerical optimization},
  author={Nocedal, Jorge and Wright, Stephen J},
  year={2006},
  publisher={Springer}
}

@article{bbstepsize,
    author = {Barzilai, Jonathan and Borwein, Jonathan M.},
    title = {Two-Point Step Size Gradient Methods},
    journal = {IMA Journal of Numerical Analysis},
    volume = {8},
    number = {1},
    pages = {141-148},
    year = {1988},
    month = {01},
    issn = {0272-4979},
    doi = {10.1093/imanum/8.1.141},
    url = {https://doi.org/10.1093/imanum/8.1.141},
    eprint = {https://academic.oup.com/imajna/article-pdf/8/1/141/2402762/8-1-141.pdf},
}

@article{Balke2022,
  title = {Scientific Computational Imaging Code (SCICO)},
  volume = {7},
  ISSN = {2475-9066},
  url = {http://dx.doi.org/10.21105/joss.04722},
  doi = {10.21105/joss.04722},
  number = {78},
  journal = {Journal of Open Source Software},
  publisher = {The Open Journal},
  author = {Balke,  Thilo and Davis,  Fernando and Garcia-Cardona,  Cristina and Majee,  Soumendu and McCann,  Michael and Pfister,  Luke and Wohlberg,  Brendt},
  year = {2022},
  month = oct,
  pages = {4722}
}

@article{Ravasi2024,
  title = {PyProximal - scalable convex optimization in
Python},
  volume = {9},
  ISSN = {2475-9066},
  url = {http://dx.doi.org/10.21105/joss.06326},
  DOI = {10.21105/joss.06326},
  number = {95},
  journal = {Journal of Open Source Software},
  publisher = {The Open Journal},
  author = {Ravasi,  Matteo and \"{O}rnhag,  Marcus Valtonen and Luiken,  Nick and Leblanc,  Olivier and Uruñuela,  Eneko},
  year = {2024},
  month = mar,
  pages = {6326}
}

@article{Ravasi2020,
  title = {PyLops—A linear-operator Python library for scalable algebra and optimization},
  volume = {11},
  ISSN = {2352-7110},
  url = {http://dx.doi.org/10.1016/j.softx.2019.100361},
  DOI = {10.1016/j.softx.2019.100361},
  journal = {SoftwareX},
  publisher = {Elsevier BV},
  author = {Ravasi,  Matteo and Vasconcelos,  Ivan},
  year = {2020},
  month = jan,
  pages = {100361}
}

@misc{tachella2025deepinverse,
      title={DeepInverse: A Python package for solving imaging inverse problems with deep learning},
      author={Julián Tachella and Matthieu Terris and Samuel Hurault and Andrew Wang and Dongdong Chen and Minh-Hai Nguyen and Maxime Song and Thomas Davies and Leo Davy and Jonathan Dong and Paul Escande and Johannes Hertrich and Zhiyuan Hu and Tobías I. Liaudat and Nils Laurent and Brett Levac and Mathurin Massias and Thomas Moreau and Thibaut Modrzyk and Brayan Monroy and Sebastian Neumayer and Jérémy Scanvic and Florian Sarron and Victor Sechaud and Georg Schramm and Romain Vo and Pierre Weiss},
      year={2025},
      eprint={2505.20160},
      archivePrefix={arXiv},
      primaryClass={eess.IV},
      url={https://arxiv.org/abs/2505.20160},
}

@misc{PyXu,
  author       = {Matthieu Simeoni and
                  Sepand Kashani and
                  Joan Rué-Queralt and
                  Pyxu Developers},
  title        = {pyxu-org/pyxu: pyxu},
  publisher    = {Zenodo},
  doi          = {10.5281/zenodo.4486431},
  url          = {https://doi.org/10.5281/zenodo.4486431}
}

@misc{ODL,
  doi = {10.5281/ZENODO.1442734},
  url = {https://zenodo.org/record/1442734},
  author = {Adler,  Jonas and Kohr,  Holger and Ringh,  Axel and Moosmann,  Julian and {Sbanert} and Ehrhardt,  Matthias J. and Lee,  Gregory R. and {Niinimaki} and {Bgris} and Verdier,  Olivier and Karlsson,  Johan and {Zickert} and Palenstijn,  Willem Jan and \"{O}ktem,  Ozan and Chen,  Chong and Loarca,  Hector Andrade and Lohmann,  Michael},
  title = {odlgroup/odl: ODL 0.7.0},
  publisher = {Zenodo},
  year = {2018},
  copyright = {Open Access}
}

@INPROCEEDINGS{ImrajHerman,
  author={Singh, Imraj RD. and Barbano, Riccardo and Twyman, Robert and Kereta, Željko and Jin, Bangti and Arridge, Simon and Thielemans, Kris},
  booktitle={2022 IEEE Nuclear Science Symposium and Medical Imaging Conference (NSS/MIC)}, 
  title={Deep Image Prior PET Reconstruction using a SIRF-Based Objective}, 
  year={2022},
  volume={},
  number={},
  pages={1-2},
  doi={10.1109/NSS/MIC44845.2022.10399292},
  url = {https://doi.org/10.1109/NSS/MIC44845.2022.10399292},
}
